\documentclass[1p]{elsarticle}

\usepackage{amsmath,amssymb,amsthm,amsfonts,afterpage,hyperref}
\usepackage{amscd,subeqnarray}

\usepackage{graphicx}

\usepackage{array}
\usepackage{wrapfig}
\usepackage{color}
\usepackage{ragged2e}
\usepackage{subfig}
\usepackage{stmaryrd}
\usepackage{siunitx}
\usepackage{multirow}
\usepackage{xcolor,cancel}
\usepackage{boldfonts}
\usepackage{symbols}
\usepackage{footmisc}
\usepackage{float}
\usepackage{setspace}
\usepackage{bm}
\usepackage{booktabs}
\usepackage{tikz, xcolor}
\usepackage{soul}

\usepackage[linesnumbered,ruled,vlined]{algorithm2e}

\usepackage[mathlines]{lineno}

\SetKwInput{KwInput}{Input}                
\SetKwInput{KwOutput}{Output} 

\usepackage{hyperref}

\usepackage{natbib}
\setcitestyle{numbers,sort&compress}

\usetikzlibrary{shapes,arrows, positioning,calc}	
\usetikzlibrary{arrows,snakes,backgrounds}
\usetikzlibrary{plotmarks}
\DeclareMathAlphabet{\mathpzc}{OT1}{pzc}{m}{it}

\newtheorem{remark}{Remark}

\journal{European Journal of Mechanics-B/Fluids}

\begin{document}

\begin{frontmatter}

%\title{A  stabilized finite element method for steady Darcy-Brinkman-Forchheimer model}
%\title{{A stabilized finite element method for steady Darcy-Brinkman-Forchheimer flow model in porous media %in 2D Lid-driven cavity flow 
%and the role of Forchheimer parameter}}
%\title{{A stabilized finite element method for steady Darcy-Brinkman-Forchheimer flow model in porous media}}
%\title{{A stabilized finite element method for steady Darcy-Brinkman-Forchheimer flow model in porous media --\\Comparison study with different viscous and inertial resistances}}% 
%\title{{A comparison study on non-Darcy flow models with different viscous and inertial resistances with a stabilized finite element method for steady Darcy-Brinkman-Forchheimer flow model in porous media}}% 
\title{{A stabilized finite element method for steady Darcy-Brinkman-Forchheimer flow model with different viscous and inertial resistances in porous media}}% 

\author[KIGAM]{Hyun Chul Yoon}
\ead{hyun.yoon@kigam.re.kr}

\author[TAMU]{S. M. Mallikarjunaiah\corref{cor1}}
\ead{M.Muddamallappa@tamucc.edu}

\cortext[cor1]{Corresponding author}

\address[KIGAM]{Marine Geology \& Energy Division,
Korea Institute of Geoscience and Mineral Resources,
124 Gwahak-ro,
Daejeon 34132, Republic of Korea}
\address[TAMU]{Department of Mathematics \& Statistics,
Texas A\&M University-Corpus Christi,
6300 Ocean Drive, Unit 5825
Corpus Christi, Texas X 78412 USA}

\begin{abstract}
We implement a stabilized finite element method for steady  Darcy-Brinkman-Forchheimer model within the continuous Galerkin framework.
The nonlinear fluid model, {a generalized and sophisticated flow model with inertial effect in porous media}, is first linearized using a standard \textit{Newton's method}. The sequence of linear problems is then discretized utilizing a stable \textit{inf-sup} type continuous finite elements based on the  \textit{Taylor-Hood} pair to approximate the 
{primary} variables{: velocity and pressure}. Such a pair is known to be optimal for the approximation of the isotropic Navier-Stokes equation. 
{To overcome {the well-known numerical instability} in the convection-dominated problems, the Grad-Div stabilization is employed with an efficient \textit{augmented Lagrangian-type} penalty method.}
We use the penalty term to develop the \textit{block Schur complement} preconditioner, which is later coupled with a Krylov-space-based iterative linear solver. 
{In addition, the Kelly error estimator for the adaptive mesh refinement is employed 
to achieve better numerical results with less computational cost.}
{Performance of the proposed algorithm is verified for a classical benchmark problem, i.e., the 2D lid-driven cavity flow of an incompressible fluid trapped between two plates.}  {{Extending the 2D lid-driven cavity flow,} a detailed comparative study of the flow nature is presented for a variety of \textit{Reynolds} number and \textit{Darcy} number combinations {for variety of viscous and inertial resistances.}} {Particularly for the Forchheimer parameter,} we present some interesting flow patterns  
with the velocity components and their streamlines along the mid-lines in the computational domain.  
The role of the Forchheimer term is  
{highlighted} for 
different porous medium scenarios. 
{This study {can offer} an attractive setting for discretizing many multi-physics problems along with the fluid flow having inertial effects in porous media.} 
\end{abstract}


\end{frontmatter}

%\linenumbers

\section{Introduction}\label{intro}
{Homogenized and simplified from Stoke's equation, the conventional Darcy's model %homogenized and simplified from Stoke's equation,
can account only for limited instances of flow in porous media.} 
For example, ignoring the inertia effect with low fluid viscosity 
may not be appropriate to describe real-world problems of convection phenomena (e.g., enhanced oil recovery in petroleum engineering, radioactive waste disposal, extraction, and geothermal energy storage). 
{More details on the limitations of Darcy's model can be found {in} \cite{vafai1995limitations,nield1984non,munaf1993boundary,rajagopal2007hierarchy,srinivasan2013flow}}.
{Thus, 
{a sophisticated non-Darcy flow model is required for flow in porous media}. %such as the the Darcy-Brinkman-Forchheimer \cite{joseph1982nonlinear,JunG2015}.} %but the analysis of a nonlinear flow and its modeling are not 
%straightforward in terms of %not only physical aspects but also some numerical aspects.}
{%As for the Darcy-Brinkman-Forchheimer, 
Joseph et al. \cite{joseph1982nonlinear} formulated a nonlinear extension of Brinkman's theory~\cite{Brink1949} for the flow of a viscous fluid through a swarm of spherical particles, and proposed the Darcy-Brinkman-Forchheimer model \cite{joseph1982nonlinear,JunG2015}. The literature on the application of Darcy-Brinkman-Forchheimer model is quite large:  to name a few, experiments on the improvement of heat and mass transfer in porous medium \cite{JinY2017,JouN2016,AvrA2020}, the fuel cell design with the inertial effect of the Forchheimer term \cite{KimJ2017}. 
%).
} 
{%For example, the 
%Due to a variety of reasons, 

%but 
{However, the analysis of a nonlinear flow and its modeling are not straightforward in terms of %not only physical aspects but also 
some numerical aspects.
Regarding the trapped flow of an incompressible fluid between two parallel plates (henceforth referred as the {\textit{lid-driven cavity flow})}}, %has been in the attention {from} computational mechanics community: 
%due to a variety of reasons: 
{the difficulties may lie in} the presence of corner singularities which occur at the points of intersection between the moving lid and other stationary walls; devising a stable numerical discretization method for the model with small viscosity (equivalently, higher \textit{Reynolds} number); {furthermore,} a motivation to design highly stable computational algorithms with remarkable accuracy. 
%It is well known that a \textit{naive Galerkin finite element method} is known to produce spurious oscillations and numerical instabilities near the corners \cite{BrenS2002, ForM2013}. Thus, 
{A more fundamental issue of incompressible fluid flow systems, such as the Navier-Stokes,
{equivalent to the free fluid for the Darcy-Brinkman-Forchheimer %\cite{joseph1982nonlinear,JunG2015} 
when the $Darcy$ number (i.e., permeability) goes to infinite}, 
is the stabilization of numerical solution. %It is well known that a \textit{naive Galerkin finite element method} is known to produce spurious oscillations and numerical instabilities near the corners \cite{BrenS2002, ForM2013}.} 
%And the Navier-Stokes system is 
%{equivalent to the free fluid for the Darcy-Brinkman-Forchheimer %\cite{joseph1982nonlinear,JunG2015} 
%when the $Darcy$ number (i.e., permeability) goes to infinite}.
Interested readers can refer the reviews \cite{gatski1991review,gresho1991incompressible} about 
mathematical issues for the stabilization of the incompressible Navier-Stokes equation. {Furthermore, for the {convection-dominated} incompressible fluid flow problems, it is well known that a \textit{naive Galerkin finite element method} is known to produce spurious oscillations and numerical instabilities near the corners \cite{BrenS2002, ForM2013}.}} %For the {convection-dominated} incompressible fluid flow problems,} various 

{For the lid-driven cavity flow,} various high-order methods have been proposed to control the oscillations and to achieve the optimal accuracy in the Sobolev space norm \cite{vallalaa2014higher,kim2016spectral,deville2002high}. %for the {convection-dominated} incompressible fluid flow problems. 
Some studies have reported the semi-analytical solutions for the flow problems with low \textit{Reynolds} number using the biorthogonal eigenfunctions \cite{shankar1993eddy,srinivasan1995accurate,wang2009recirculating,muddamallappa2009numerical}. For the flow problems with moderate or high  \textit{Reynolds} number in a cavity with a top moving plate, many strategies have been proposed in the literature, some of which 
include regularizing the whole boundary value problem by altering the top-lid tangential velocity condition to a polynomial that vanishes at the endpoints \cite{shen1991hopf,barragy1997stream}.  
{However, such methods are computationally expensive, and particularly in \cite{barragy1997stream}, the underlying finite element test function space approximating the velocity variable still needs to be far smoother.} {The  Streamline Upwind Petro-Galerkin (SUPG) method  is another popular countermeasure to the classical Galerkin finite element formulation that stabilizes the computations for higher  
\textit{Reynolds} number   \cite{brooks1982streamline,wang2016supg,tezduyar2000finite}, and also in compressible flow simulations \cite{tezduyar2006stabilization}.
An open and challenging issue in the SUPG method is analyzing the coupling between velocity and pressure terms in the discrete variational formulation} {(see the works in \cite{franca1992stabilized,hansbo1990velocity} for the issues with SUPG method)}. Meanwhile, the study from Barragy and Carey~\cite{barragy1997stream} used $p$-type finite elements for the vorticity-stream function formulation and reported a highly accurate solution for the steady cavity problem. In Botella and Peyret's work~\cite{botella1998benchmark}, a Chebyshev collocation method is used for the lid-driven cavity flow problem by subtracting the leading order terms from the asymptotic expansion of the solution to the steady Navier-Stokes equation,
where %a 
{very high-order polynomials} are still used in the approximation.    
Further, %in other approaches,  
many researchers have used the velocity-vorticity formulations rather than computing velocity-pressure as %primitive 
{primary} variables \cite{davies2001novel}  
due to the lack of an independent equation for pressure. However,  
the boundary conditions for the vorticity variable at the corner of the top lid could be more straightforward. 

%The primary objective of this paper is two-fold: {first,} to 
{{In this study,} we devise a stabilized finite element method %and illustrate its performance for the lid-driven cavity filled with 
for the stationary Darcy-Brinkman-Forchheimer model with incompressible, isotropic, Newtonian fluid filled in porous media, and illustrate its performance through the lid-driven cavity flow problem.}
%in the Darcy-Brinkman-Forchheimer porous medium}.} %{second, to} 
{Using the stabilized solutions, we also highlight} %present 
some important flow features -- such as the effect of dimensionless constants and the role of the Forchheimer term.  To that end, the \textit{inf-sup} stable finite element approximation %with the %primitive 
%{primary} variable pair of velocity-pressure 
and a Grad-Div stabilization 
is employed with the {primary} variable pair of velocity-pressure. %in the stationary Darcy-Brinkman-Forchheimer model.
The Grad-Div stabilization, suggested in various studies \cite{franca1992stabilized,whiting2001stabilized,braack2007stabilized,olshanskii2002low}, is another popular and robust approach to reducing the velocity error due to the pressure error and  
preserving the mass conservation.  
%\magenta{cite the step-57 correctly} 
%In this study, 
{Particularly motivated by \cite{ZhaL2017}, we employ} the product of a stabilization parameter and the Grad-Div term %is added 
to the momentum equation, {where the} choice of the stabilization parameter is vital for the accuracy of numerical solutions \cite{olshanskii2009grad}. An effecient Schur complement-based preconditioner is develped by utilizing the block-strucutre of the discrete system \cite{ZhaL2017}. 
%This 
{which reduces the number of Newton iterations and makes the whole process independent of the total degrees of freedom.}  Besides the preconditioner for the linear solver, 
the adaptive mesh refinement based on the {\textit{Kelly} error estimator \cite{kelly1983}} is employed 
to achieve better numerical results with less computational cost.

{To validate} our finite element discretization along with the preconditioner and corresponding implementation, 
the computational code is tested  
{with} a manufactured solution for the optimal convergence rate in $H^1$- and $L^2$-norms, {where we report the optimal orders %of 
%convergence rates 
for both velocities and pressure.} %in terms of $H^1$- and $L^2$-norms.} %with a manufactured solution.} 
{Our finite element solution is also compared against a finite difference solution in \cite{Ghia1982} and reports an excellent agreement between the two results.} 
%{From the numerical experiments,} we report the optimal-order convergence rates for both velocities and pressure in terms of $H^1$- and $L^2$-norms with a manufactured solution.
{We then perform further numerical experiments for different viscous and inertial resistances in porous media. Two groups of \textit{Reynolds} and \textit{Darcy} numbers are set, for which \textit{Re} $\times$ \textit{Da} values are less than and greater than 1.0.}
{Regarding the numerical performance, we find the preconditioner-based algebraic penalty term stabilizes the {relatively} high Reynolds number (up to $1000$). %{for the porous media} problem %our method takes \hl{few iterations (less than $10$ particularly)} within Newton's method. 
Furthermore, our method takes {few iterations within Newton's method even for the nonlinear model}. %the preconditioner-based algebraic penalty term stabilizes the {relatively} high Reynolds number (up to $1000$) {for the porous media} problem; 
%In addition, the whole procedure is shown to be effective and independent of the large-scale sparse problem due to local mesh adaptivity.
Due to local mesh adaptivity, the whole procedure is shown to be independent of the degree of freedom and effective for large-scale sparse problem even when \textit{Re} $\times$ \textit{Da} are increased. 
%are then 
Our study further focuses on the comparison, i.e., %of 
%between 
the Brinkman and Darcy-Brinkman models without the Forchheimer term %compared to 
with the full nonlinear Darcy-Brinkman-Forchheimer.} 
For two groups of \textit{Reynolds} and \textit{Darcy} numbers, the primary, secondary, and tertiary vortices of the models are compared. 
%{and} the role of Forchheimer term for the turbulent flow regime in the porous medium is specifically highlighted.
%Therefore, the 
%Numerical advantages of our overall algorithm are three-fold: first, 
%{For the numerical performance, we find our method takes few iterations (less than $10$) within Newton's method; furthermore, the preconditioner-based algebraic penalty term stabilizes the {relatively} high Reynolds number (up to $1000$) {for the porous media} problem; in addition, the whole procedure is shown to be effective and independent of {the large-scale sparse problem due to local mesh adaptivity.}} 
%For the 
In the Darcy-Brinkman-Forchheimer model, we identify distinct streamline patterns and the departures from other models in the location and size of primary, secondary, and tertiary vortices. {We also highlight the role of Forchheimer term for the turbulent flow regime in the porous medium.}
The Forchheimer term is found to play its role as the resistance to the inertial force and convection, delaying the flow regime change {between the laminar and turbulence regime of flow in porous media}. 

{The paper is organized as follows: Section~\ref{sec:GoverningEquations} describes the mathematical model that characterizes the flow of incompressible fluid in a non-Darcy porous medium.  Some essential nations from functional analysis,  a damped Newton's method, a corresponding variational formulation, and a stable finite element discretization are %all 
described in Section~\ref{numerical}. In the same section, a linear solver, preconditioner, and Kelly error indicator-based local mesh refinement are explained in detail.  Section~\ref{num_exp} contains the numerical experiments using the proposed finite element algorithm. The conclusions and some directions for future works are outlined in Section~\ref{conclusion}. }

\section{Governing Equations}\label{sec:GoverningEquations}
\label{sec:background}

%{The literature on the application of Darcy-Brinkman-Forchheimer model is quite large:  to name a few, experiments on the improvement of heat and mass transfer in porous medium \cite{JinY2017,JouN2016,AvrA2020}, the fuel cell design with the inertial effect of the Forchheimer term \cite{KimJ2017}.} Joseph et al. \cite{joseph1982nonlinear} formulate anonlinear extension of Brinkman's theory~\cite{Brink1949} for the flow of a viscous fluid through a swarm of spherical particles. 
For the momentum balance of the steady flow, {the following nonlinear equation, i.e., the Darcy-Brinkman-Forchheimer model} governs the flow in a fully saturated porous
medium when the inertia effects and the fluid incompressibility constraint are {considered}:
\begin{subequations}\label{eq:Darcy-Brinkan-Forch} 
\begin{align}
\rho \mathbf{u} \cdot \nabla\mathbf{u}+\nabla p  &=\mu \Delta
\mathbf{u}-\phi\frac{\mu}{K}\mathbf{u}
-\phi\frac{\rho\,\beta  }{ \sqrt{K} } \, | \mathbf{u} |^{\alpha} \, \mathbf{u} + \widetilde{\mathbf{f}},  \quad {\mbox{in} \quad \Omega}\label{eq:Darcy-Brinkan-Forch-a}\\
\nabla \cdot \mathbf{u} &=0,  \quad {\mbox{in} \quad \Omega}, 
\end{align}
\end{subequations}
where $\mathbf{u} \colon \Omega \to \mathbb{R}^2$ is the velocity, $p \colon \Omega \to \mathbb{R}$ is the pressure, in which $\Omega \subset \mathbb{R}^2$ is a 
smooth, bounded, non-convex domain, {and $\partial \Omega$ is its boundary}. $K$ and $\phi$ are
the intrinsic permeability and porosity of the medium, respectively. Also, $\mu$ is the coefficient of
dynamic viscosity of the fluid, 
$\rho$ is the fluid density. {$\beta$ is the Forchheimer coefficient, which is a positive quantity for the inertial resistance depending on
the geometry of the medium.} %having the unit of reciprocal length.} 
%Moreover, the right hand side term $\widetilde{\mathbf{f}}$ represents an external force (e.g., externally applied magnetic field). 
{The symbol $| \cdot |$ is the Euclidean norm, and $| \mathbf{u} |^2 = \mathbf{u} \cdot \mathbf{u}$}. 
%In this study, intrinsic permeability ($K$) and inertial resistance ($\beta$) are taken as scalar coefficients instead of full tensor expression for an assumption of isotropic medium. {For simple analysis, we set the porosity constant as $\phi=1.0$, unless otherwise noted. Although the porosity is a variable depending upon the fluid pressure or volumetric strain of porous medium, setting it constant can eliminate its effect on the flow, and we can purely investigate numerical aspects of the model and perform comparison study between the models.}
%Also, for simplicity %purposes only, 
%purpose, the porosity of the medium is assumed to be uniform and constant as $\phi=1.0$. 
{Note that the Forchheimer term has a power, $\alpha$, generally regarded as $1.0\leq\alpha\leq2.0$ \cite{FirdM1997,CimF2013}. %Although the exact value may depend on the physical nature of the model, %for illustration purposes, 
In this study, we take  $\alpha=1.0$, resulting in the quadratic form deviating from the linear Darcy's law.} %Although the exact value may depend on the physical nature of the model, the validity of the quadratic form was confirmed in [cite-Chen et al., 2001] through matching with the experimental results [McDonald et al., 1979].} \\
Moreover, the right hand side term $\widetilde{\mathbf{f}}$ represents an external force (e.g., externally applied magnetic field). 

{{Assuming $\phi=1.0$} for simplicity}, we introduce the following dimensionless quantities: 
\begin{equation}
\mathbf{x}^{*} = \frac{\mathbf{x}}{L},  \quad 
\mathbf{u}^{*} = \frac{\mathbf{u}}{U_{0}},  \quad
p^{*} = \frac{p}{\rho \, U_{0}^{2}},   
\end{equation}
where $L$ and $U_{0}$ are the characteristic properties of length and velocity, respectively. %, of which scales depend on a problem. In essence, 
{Their scales depend on a problem, %. In essence,}
and for a porous medium, $L$ is the reference distance (e.g., the mean diameter of grain)} 
and {$U_{0}$ denotes the magnitude of the reference discharge relative to the solid grain, having the unit of velocity.} 
We omit the asterisks, $(\cdot)^{*}$, {of the dimensionless quantities} henceforth for simplicity of the notation. {Under the assumption of no external force,} the dimensionless equations of motion are as follows:
\begin{subequations}\label{eq:Darcy-Brinkan-Forch-1}
\begin{align}
 \mathbf{u} \cdot \nabla\mathbf{u}+\nabla p  &= \underbrace{\frac{1}{Re}}_{\displaystyle (a)} \Delta
\mathbf{u}- \underbrace{\frac{1}{Re \, Da}}_{\displaystyle (b)}  \, \mathbf{u} 
-  \underbrace{\frac{c_F}{\sqrt{Da}}}_{\displaystyle (c)} \,  \| \mathbf{u} \| \, \mathbf{u}, \label{eq:dimless-full-nonlinear} \\
\nabla \cdot \mathbf{u} &=0, \quad \label{eq:imcompressible}
\end{align}
\end{subequations}
where $c_F$ is the dimensionless coefficient for inertial resistance. {$Re$ and $Da$ are characteristic and dimensionless quantities of} the \textit{Reynolds} and \textit{Darcy} numbers, respectively,  
defined as 
\begin{align}
Re &= \frac{ \rho \, L \, U_0}{\mu}, \label{eq:Re} \\
Da &= \frac{K}{L^2}. \label{eq:Da}
\end{align}
As seen in \eqref{eq:Re}, note that {the inverse of $Re$ number ($(a)$ in \eqref{eq:dimless-full-nonlinear}) %, i.e., $(a)$ in \eqref{eq:dimless-full-nonlinear}, 
is the ratio of viscous force to inertial force for a fluid.} {Then, the coefficient} $ (b) \dfrac{1}{Re \, Da}$ in \eqref{eq:dimless-full-nonlinear} plays its role as the viscous resistance of the fluid flow through a porous medium, whereas the coefficient $(c) \dfrac{c_F}{\sqrt{Da}}$ in \eqref{eq:dimless-full-nonlinear} represents the inertial resistance, %, respectively, 
%of the fluid flow through a porous medium, 
respectively.

\begin{remark}
%The Forchheimer term has been criticized for its role in the nonlinear effect regardless of the magnitude of the velocity, especially for an incompressible fluid with a steady uni-directional flow~\cite{CimF2013}. Nonetheless, 
{The Forchheimer term with the power $\alpha$ in \eqref{eq:Darcy-Brinkan-Forch-a} has its basis derived %experimentally 
in the framework of the homogenization theory {of} the pore scales, and it can have various forms (e.g., the cubic form for the weak inertia)~\cite{FirdM1997,ZimR2004}. Although the exact value may depend on the physical nature of the model, the validity of the quadratic form was confirmed in \citep{CheZ2001}, through matching with the experimental results \citep{MacI1979}.} {Regarding the inertial resistance (or the Forchheimer) coefficient,} 
there are several methodologies for its determination 
 based on both theoretical and empirical aspects, which depend on the corresponding properties of the porous medium. For example, it can be theoretically determined by employing the {Ergun} equation~\cite{ErgS1957, YanD2012}, while empirically by using a correlation between the coefficient and the permeability, which can result in {certain range of} values due to the scale difference of environments between the lab and reality.   
In this study, 
neither aspects of the theoretical nor empirical {backgrounds are utilized}, but %instead, 
we fix its value within a {typical} range, {according to} the literature.   
 %instead. 
\end{remark}
\begin{remark}
{Setting the porosity $(\phi)$ value as unity} is based on the simplification process, and this framework %can be categorized as a 
{is the same in its form as the} penalization method~\cite{BrunC2008,CimF2013} %{Compared to 
{for {fully}} coupling method 
{between the Navier-Stokes equation and the} {{full} nonlinear model}{, i.e.,} the Darcy-Brinkman-Forchheimer. {Note the resistance terms can also imply 
porous medium since the resistance can be induced by porosity, i.e., the flow channel, and its property.} %{without explicitly applying the domain .} 
It can also be viewed
as a fictitious domain approach~\cite{CimF2013,AngP1999,KhadK2000} for modeling fluid flow resistance in porous media. 
\end{remark}
\begin{remark}\label{rmk2-3}
As for the {full} nonlinear Darcy-Brinkman-Forchheimer {model} \eqref{eq:dimless-full-nonlinear}, some {studies} (e.g., see \cite{GutR2015}) have omitted the Forchheimer term, while other literature (e.g., see \cite{JunG2015}) has no convection term. In this study, note that 
{we name an equation as the} (linearized) Brinkman {model}, {when}
\begin{equation}\label{eq:Brinkman}
\nabla p  = \frac{1}{Re} \Delta
\mathbf{u}- \frac{1}{Re \, Da}  \, \mathbf{u},
\end{equation}
whereas the Darcy-Brinkman {model}, 
{being} expressed as 
\begin{equation}\label{eq:Darcy-Brinkman}
\mathbf{u} \cdot \nabla\mathbf{u}+\nabla p  = \frac{1}{Re} \Delta\mathbf{u}- \frac{1}{Re \, Da}  \, 
\mathbf{u}. 
\end{equation}
The only difference between the two models above lies in the nonlinear convective term, $\mathbf{u} \cdot \nabla\mathbf{u}$. 
Note that the Darcy-Brinkman(-Forchheimer) model is
{convergent} to the Navier-Stokes model {when} $Da$ goes to {large number,}
{ultimately being identical to it, if} $Da\rightarrow\infty$. 
 \end{remark}

\section{Numerical Methods}\label{numerical}
{In this section, numerical methods employed 
for the rigorous numerical results of this study 
are explained. In Section~\ref{sec:prelim}, 
some notations about the functional spaces are defined and introduced. The nonlinear governing equations, the Darcy-Brinkman-Forchheimer with the incompressibility constraint, are linearized first using Newton's method in Section~\ref{sec:Newton}. The linearized equations are turned into the variational formulation in Section~\ref{sec:variational_form}. In Section~\ref{sec:FEM}, the stabilized finite element method is introduced, and some preconditioners for the linear solver employed in this study are explained in Section~\ref{sec:preconditioner}. Finally, the local adaptive mesh refinement method based on the {{Kelly} error estimator \cite{kelly1983}} is described in Section~\ref{sec:LGR}.}
\subsection{Notations and preliminaries}\label{sec:prelim}
{Throughout this 
{study,} we assume $\Omega \subset \mathbb{R}^{2}$ a {non-convex}, bounded polygonal domain  with boundary $\partial  \Omega$}. 
Let $X$ be a normed linear space defined
 for $\Omega$ and $L^{p}\left( X    \right)$ is a space of all {\textit{Leb$\grave{e}$sgue}} integrable functions defined on $\Omega$,  in particular 
\begin{equation}
L^{2}\left( \Omega  \right) := \left\{ q \colon \Omega \to \mathbb{R} \, | \, \int_{\Omega} |q|^{2} d\boldsymbol{x} < \infty  \right\}
\end{equation} 
is the classical space of \textit{square-integrable functions}, and the corresponding inner-product and norm are defined as
\begin{equation}
\left(p, \, q  \right):= \int_{\Omega} p \, q \, d\boldsymbol{x} \quad \mbox{and} \quad|| q || := \left( \int_{\Omega} |q|^{2} d\boldsymbol{x}  \right)^{1/2}.
\end{equation}
Also, the norm $|| \cdot ||_{p}$ is defined as
\begin{equation}
|| u ||_{p} = \left( \int_{\Omega} || u ||^{p}_{X} \; d\boldsymbol{x} \right)^{1/p}, \quad p \in [1, \; \infty).
\end{equation}
 
 {For the problem on hand,} the natural Hilbert space is defined as follows:
\begin{align} 
 \mathbf{X}_{0} &=  H_{0}^{1}\left( \Omega   \right)^{2} := \left\{ {\mathbf{v}} \colon \Omega \to \mathbb{R}^2 \;  | \; {\mathbf{v}} \in  \left( L^{2}\left( \Omega \right) \right)^{2}, \;  \nabla {\mathbf{v}} \in  \left(  L^{2}\left( \Omega \right) \right)^{2 \times 2}, \; {\mathbf{v}} = \mathbf{0} \;\; \mbox{on} \;\; \partial \Omega \right\},  \\
 Q_0 &= L_{0}^{2} \left( \Omega \right)  := \left\{ {q} \in L^{2} \left( \Omega \right) \colon \int_{\Omega} {q} \; d\boldsymbol{x} = 0  \right\},
\end{align}
and further $\mathbf{X} = H^{1}\left( \Omega   \right)^{2}$, $Q=L^{2}\left( \Omega \right)$. The $H^{1}$-norm is then defined as 
\begin{equation}
\| \mathbf{u} \|_{H^{1}(\Omega)} := \left( \| \mathbf{u} \|_{2}^2 + \| \nabla \mathbf{u} \|_{2}^2  \right)^{1/2}.
\end{equation}

\subsection{Linearization using Newton's method}\label{sec:Newton} 
Our main goal in this paper is to propose a stable finite element discretization for the nonlinear {partial differential equations (PDEs) system for} the Darcy-Brinkman-Forchheimer {model} 
described in Section~\ref{sec:GoverningEquations}. A common strategy for the numerical solution of nonlinear differential equations is to first obtain a sequence of linear problems by using Picard's or Newton's method at the PDE level, which is then solved by an appropriate numerical technique \cite{charnyi2019efficient,kelley2003solving,ortega2000iterative}. Alternatively, one can discretize the nonlinear partial differential equations and obtain a system of nonlinear algebraic equations, which are then solved by any befitting linear algebra techniques. Both methods are mathematically identical, and 
there 
{is less likely to be}
{any} difference {in computational perspective}. {One can also linearize the continuous variational formulation and then solve the resulting formulation by an appropriate numerical method for convergence.  In this work, we linearize at the differential equation level, and then the linear problems are solved using the continuous Galerkin finite element method.} 

Let $\mathbf{x}^{n} = \left( \mathbf{u}^n, \; p^n \right)$ be a numerical approximation to the primary variables, $\mathbf{u}$ and $p$,
and $n$ denotes \textit{iterative step} number for the Newton's method.
We then write the system \eqref{eq:dimless-full-nonlinear}-\eqref{eq:imcompressible} by moving all the terms to the left-hand side and obtain
\begin{equation}
F(\mathbf{x}):=F(\mathbf{u}, \, p) = 
\begin{pmatrix}
 \mathbf{u} \cdot \nabla\mathbf{u}+\nabla p  - \dfrac{1}{Re} \, \Delta
\mathbf{u} + \dfrac{1}{Re \, Da} \, \mathbf{u}  
+ \dfrac{c_F}{\sqrt{Da}} \| \mathbf{u} \|  \mathbf{u}  \\
- \nabla \cdot \mathbf{u}
\end{pmatrix}.
\end{equation}
{{We then} see 
{an} approximation to $\mathbf{x}^{n}$  of the form} 
\begin{equation}
\mathbf{x}^{n+1} = \mathbf{x}^{n} + {\widehat{\alpha}} \, \delta \mathbf{x}^{n}. \label{NMethod-1}
\end{equation}
In \eqref{NMethod-1}, the term $\delta \mathbf{x}^{n}$ is known as the \textit{update} or the \textit{correction} term, and ${\widehat{\alpha}}$ is a line search parameter that helps to achieve the quadratic convergence property of Newton's method, {for which we use $\alpha=1.0$ in this study}. The updates are obtained  using 
\begin{equation}
\delta \mathbf{x}^{n} = - \left( \nabla F(\mathbf{x}^n) \right)^{-1} F(\mathbf{x}^n),  \label{NMethod-2}
\end{equation}
and $\nabla F(\mathbf{x}^n) $ denotes the \textit{Jacobian} of $ F(\mathbf{x})$ evaluated at $\mathbf{x}^n$. Finally, Newton's iterative algorithm for the current mixed {fields} of velocity and pressure can be formulated as a combination of \eqref{NMethod-1} and \eqref{NMethod-2}. The equation \eqref{NMethod-2} can also be written as 
\begin{equation}
\nabla F(\mathbf{x}^n) \, \delta \mathbf{x}^{n} = - \, F(\mathbf{x}^n).  \label{NMethod-3}
\end{equation}
The term on the left-hand side of \eqref{NMethod-3} is the \textit{directional derivative} of $F(\mathbf{x}^n)$ in the direction of $\delta \mathbf{x}^{n}$. We compute the {directional derivative}  by using the basic definition 
\begin{align}
& \nabla F(\mathbf{x}^n) \, \delta \mathbf{x}^n =                      \nabla F(\mathbf{u}^{n}, \, p^{n})(\delta \mathbf{u}^n, \, \delta p^n) \\
& = \lim_{\epsilon \to 0} \frac{F(\mathbf{u}^{n} + \epsilon  \delta \mathbf{u}^n, \, p^n + {\epsilon \delta p^n}) - F(\mathbf{u}^n, \, p^n) }{\epsilon}  \\
&=\begin{pmatrix}
  \mathbf{u}^{n}\cdot \nabla \delta \mathbf{u}^{n} + \delta \mathbf{u}^{n}\cdot \nabla \mathbf{u}^{n}  + \nabla \delta p^n   -  \dfrac{1}{Re} \, \Delta \delta \mathbf{u}^{n} + \dfrac{1}{Re \, Da} \,  \delta  \mathbf{u}^n  \\
  +  \dfrac{c_F}{\sqrt{Da}}  \left(  \| \mathbf{u}^n \|    \delta \mathbf{u}^n  + \dfrac{\mathbf{u}^n \cdot \delta \mathbf{u}^n }{\| \mathbf{u}^n \|} \mathbf{u}^n  \right)   \\
  \\
- \nabla \cdot \delta \mathbf{u}^n
\end{pmatrix}.\label{eq:direc_Jacobian}
\end{align}
The resultant linearized {PDEs} system that
is to be solved at every \textit{Newton step} for $\delta \mathbf{u}^n$ and {$\delta p^n$} reads
\begin{multline}\label{lin1}
\begin{pmatrix}
  \mathbf{u}^{n}\cdot \nabla \delta \mathbf{u}^{n} + \delta \mathbf{u}^{n}\cdot \nabla \mathbf{u}^{n}  + \nabla \delta p^n   -  \dfrac{1}{Re} \, \Delta \delta \mathbf{u}^{n} + \dfrac{1}{Re \, Da} \,  \delta  \mathbf{u}^n  \\
  +  \dfrac{c_F}{\sqrt{Da}}  \left(  \| \mathbf{u}^n \|    \delta \mathbf{u}^n  + \dfrac{\mathbf{u}^n \cdot \delta \mathbf{u}^n }{\| \mathbf{u}^n \|} \mathbf{u}^n  \right)   \\
  \\
- \nabla \cdot \delta \mathbf{u}^n
\end{pmatrix} \\ = - F(\mathbf{u}^n, \, p^n),
\end{multline}
and the right-hand side is given by
\begin{equation}\label{rhs}
F(\mathbf{u}^n, \, p^n) = 
\begin{pmatrix}
\mathbf{f}^n  \\
- \nabla \cdot \mathbf{u}^n
\end{pmatrix},
\end{equation}
where $\mathbf{f}^n$ denotes
\begin{equation}
\mathbf{f}^n = \mathbf{u}^n \cdot \nabla\mathbf{u}^n+\nabla p^n  - \dfrac{1}{Re} \, \Delta
\mathbf{u}^n + \dfrac{1}{Re \, Da} \, \mathbf{u}^n  
+ \dfrac{c_F}{\sqrt{Da}} \| \mathbf{u}^n \|  \mathbf{u}^n.
\end{equation}
Note that $\mathbf{u}^n$ and $p^n$ are the {discrete} solutions obtained in the previous iterations, and the unknowns in \eqref{lin1} are the \textit{Newton updates}, $\delta \mathbf{u}^n$ and $\delta p^n${, i.e., the components of $\delta \mathbf{x}^n$ in \eqref{NMethod-1}}. {To} obtain the linear algebraic system for the updates in this study, we will use the standard Galerkin conforming finite element method with appropriate function spaces. Meanwhile, the linearized problem \eqref{lin1}-\eqref{rhs} must be supplied with an appropriate boundary condition {for well-posedness.} 
{Suppose} 
the non-homogenous Dirichlet boundary conditions are given in the form
\begin{equation}
\mathbf{u} = \mathbf{g} \quad \mbox{on} \quad \partial \Omega,
\end{equation} 
{we then} impose the {(\textit{known}) non-zero values 
on} the initial guess and zero 
values on the \textit{Newton updates} {for boundary conditions} in the subsequent iterations, i.e.,
\begin{align}
\mathbf{u}^{0} &=  \mathbf{g},  \quad \mbox{on} \quad {\partial\Omega}\\
\delta \mathbf{u}^n &=0, \quad \mbox{for} \quad n > 0 \quad  \mbox{ %and 
{on}} \quad {\partial \Omega}. 
\end{align}
Note that $\mathbf{u}^{n+1}=\mathbf{u}^{n} + \delta \mathbf{u}^{n}$ always guarantees $\mathbf{u}^{n+1} = \mathbf{g}$ on $\partial \Omega$.

It is well understood that the \textit{Newton's method} %needs 
necessitates a good initial guess for better convergence, and %the 
a particular initial guess needs to be ``close'' {enough} to the exact solution. For a better initial guess {in this study}, {we solve \eqref{eq:Darcy-Brinkman} (namely, what we call as the {Darcy-Brinkman} model) for the full {Darcy-Brinkman-Forchheimer model}
on a coarser mesh and then interpolate the solution to the newly refined mesh. %, {of which approach, i.e., 
This \textit{local grid refinement} approach is described in detail in Section~\ref{sec:LGR}}. 

\subsection{Variational formulation}\label{sec:variational_form}
{In this section, {for the standard Galerkin conforming finite element method with appropriate function spaces,}}
{a} continuous variational formulation of the problem is first introduced within the suitable Sobolev spaces for velocity and pressure solution fields.  Let {$\mathbf{X}, \; \mathbf{X}_0, \;Q$, and $Q_0$} are the functional spaces defined in  
{Section}~\ref{sec:prelim}  
 and let us introduce the following bilinear, trilinear, and semilinear forms: 
\begin{subequations}\label{lforms}
\begin{align}\label{sub-a}
&m \colon \mathbf{X} \times \mathbf{X} \to \mathbb{R}, \quad m(\delta \mathbf{u}^n, \; \mathbf{v}) = \dfrac{1}{Re} \, \left(\nabla \delta \mathbf{u}^n, \; \nabla \mathbf{v} \right), \\
&c \colon \mathbf{X} \times \mathbf{X} \to \mathbb{R}, \quad c(\delta \mathbf{u}^n, \; \mathbf{v}) = \dfrac{1}{Re \, Da} \, \left(\delta \mathbf{u}^n, \;  \mathbf{v} \right), 
  \\
&n \colon \mathbf{X} \times  \mathbf{X} \times \mathbf{X} \to \mathbb{R}, \quad n(\mathbf{u}^n, \; \delta \mathbf{u}^n, \; \mathbf{v}) = \left( \mathbf{u}^n \cdot \nabla \delta \mathbf{u}^n, \;  \mathbf{v} \right) + \left( \delta \mathbf{u}^n \cdot \nabla  \mathbf{u}^n, \;  \mathbf{v} \right), \\
&s \colon \mathbf{X} \times  \mathbf{X} \times \mathbf{X} \to \mathbb{R}, \quad s(\mathbf{u}^n; \; \delta \mathbf{u}^n, \; \mathbf{v}) = \dfrac{c_F}{\sqrt{Da}} \left[ \left(  \| \mathbf{u}^n \| \, \delta \mathbf{u}^n, \;  \mathbf{v} \right) +
 \left( \dfrac{\mathbf{u}^n \cdot \delta \mathbf{u}^n}{\| \mathbf{u}^n \|} \; \mathbf{u}^n, \; \mathbf{v} \right) \right],\label{sub-d}
 \\
&b \colon \mathbf{X} \times Q \to \mathbb{R}, \quad b(\delta \mathbf{u}^n, \; q) =  \left( \nabla \cdot  \delta \mathbf{u}^n, \;  q \right).
\end{align} 
\end{subequations}
Furthermore for \eqref{sub-a}-\eqref{sub-d}, we define ${\widetilde{a} \colon \mathbf{X} \times  \mathbf{X} \times \mathbf{X} \to \mathbb{R}}$, as 
\begin{equation}\label{eq:a_trilinear}
{\widetilde{a}(\mathbf{u}^n; \; \delta \mathbf{u}^n, \; \mathbf{v})} := m(\delta \mathbf{u}^n, \; \mathbf{v})  + c(\delta \mathbf{u}^n, \; \mathbf{v}) + n(\mathbf{u}^n, \; \delta \mathbf{u}^n, \; \mathbf{v}) + s(\mathbf{u}^n; \; \delta \mathbf{u}^n, \; \mathbf{v}). 
\end{equation}
For a fixed $\mathbf{u}^n$ {in \eqref{eq:a_trilinear}, {we obtain {$\widetilde{a}(\mathbf{u}^n; \; \cdot, \, \cdot)$} as bilinear}. Also, the solution $\mathbf{u}^n$ is \textit{divergence-free}. Then, to derive a weak formulation of the continuous  problem, {multiplying} the linearized version of {momentum and mass balance equations} in \eqref{lin1} by test functions $\mathbf{v} \in \mathbf{X}_0$ and $q \in Q_0$ with integrating by parts implies the following weak formulation: \\

\noindent \textbf{Continous weak formulation: } \textit{Find $\mathbf{u}^{n+1} \in \mathbf{X}$ and {$p^{n+1} \in Q$} with  {$\mathbf{u}^{n+1}|_{\partial\Omega} = \mathbf{g}$},  given $\left( \mathbf{u}^0, \; p^0 \right) \in {\mathbf{X}} \times Q$, and for each $n = 0, \, 1, \, 2, \ldots$, 
such that}
\begin{subequations}\label{f1}
\begin{align}
&\widetilde{a}(\mathbf{u}^n; \; \delta \mathbf{u}^n, \; \mathbf{v})  - b(\mathbf{v},  {\delta p^n}) = \left( \mathbf{f}^n, \; \mathbf{v}  \right),  \label{f1-a} \\
&b(\delta \mathbf{u}^n, \, q) = 0, \label{f1-b}
\end{align}
\end{subequations}
where $\forall \; \mathbf{v} \in \mathbf{X}_0$, $\forall \;  q \in Q_0$. \\

%Also notice that since $\Omega$ is bounded and connected, there exists a constant $c$ satisfying the continuous \textit{inf-sup} condition:
%\begin{equation}\label{infsup1}
%\sup_{{\mathbf{v}\neq\mathbf{0},}\;{\mathbf{v}} \in \mathbf{X}} \;\; \dfrac{\left( p, \; \nabla \cdot {\mathbf{v}}     \right)}{  \| p \|_{ Q} \; \| {\mathbf{v}} \|_{ \mathbf{X}}} \geq c \quad \forall \, p \in Q_0.
%\end{equation}
%Additionally, since $\mathbf{u}^n$ is fixed and \textit{divergence-free}, the blinear form $\widetilde{a}(\mathbf{u}^n; \; \cdot, \, \cdot)$ is both coercive and continuous over $\mathbf{X}$:
%\begin{align}
%\widetilde{a}(\mathbf{u}^n; \; \mathbf{v}, \, \mathbf{v}) &\geq \beta_1 \, \| \nabla \mathbf{v} \|_{\mathbf{X}} \quad \forall \mathbf{v} \; \in \mathbf{X}, \\
%| \widetilde{a}(\mathbf{u}^n; \; \mathbf{u}, \, \mathbf{v}) | &\leq \beta_2 \, \| \nabla \mathbf{u} \|_{ \mathbf{X}} \,  \| \nabla \mathbf{v} \|_{ \mathbf{X}} \quad  \forall \mathbf{u} \; \in \mathbf{X}, \; \forall \mathbf{v} \; \in \mathbf{X}.
%\end{align}
%Where the constants $\beta_1$ and $\beta_2$ are, respectively, the coercivity and continuity constants. The existence and uniqueness of a solution to \eqref{f1} follows from {the} \textit{Lax-Milgram lemma} \cite{girault2012finite,ern2013theory}.

\begin{remark}
{It is easy to see that the solution $\left( \delta \mathbf{u}^{n}, \;  \delta p^{n}     \right)$, satisfying \eqref{lin1}, also satisfies \eqref{f1} along with proper boundary conditions.} 
However, when interpolated from the coarse mesh, the initial guess might not have the divergence-free property. In that case, \eqref{f1-b}  
should then be replaced by 
\begin{equation}\label{eq:b}
b(\delta \mathbf{u}^n, \, q) =  - b( \mathbf{u}^n, \, q).
\end{equation}
\end{remark}
The solution to the linearized problem \eqref{f1} can be obtained by a suitable discretization using {a stabilized finite element method} combined with a convenient stopping criteria for the \textit{Newton's method}.

\subsection{Finite element discretization}\label{sec:FEM}
{We introduce a {stabilized finite element discretization} of the formulation described in {\eqref{f1-a} and \eqref{eq:b}}.} 
We discretize the domain $\Omega$ as follows; 
let $\mathpzc{T}_{h}$ be a subdivision of $\Omega$ that possesses \emph{quadrilateral} elements $K \in  \mathpzc{T}_{h}$ such that $\overline{\Omega} = \bigcup_{K \in \mathpzc{T}_{h}} \overline{K}$ and $K_{i} \cap K_{j} = \emptyset$ for $K_{i}, \; K_{j} \, \in \mathpzc{T}_{h}$ with $i \neq j$. 
The boundary of $\Omega$ is discretized as the union of faces from  $\mathpzc{T}_{h}$ {and} we assume that the subdivision of $\Omega$ is conforming. {Then, the intersection of the closure of any two elements $K_{i}$ and $K_{j}$ is either empty {or} along a face.} We also assume that this discretization is \textit{shape-regular} (or regular) and {\textit{simplicial}} 
 in the sense of Ciarlet \cite{ciarlet2002finite}. 
 
 The finite element space considered in this work is  $\mathcal{Q}_{k}, \; k \in \mathbb{N}$, is defined by
\begin{equation}
\mathcal{Q}_{k} := \left\{ \mathbf{v} \in C(\Omega) \,  \colon \; \mathbf{v} |_{K} \in \mathbb{Q}_k, \; K \in \mathpzc{T}_{h}  \right\},
\end{equation}
{where $\mathbb{Q}_k$ is the tensor-product of polynomials up to order $k$ on each element.}  

We then describe discrete spaces to approximate velocity and pressure variables. Any traditional {equal-order} spaces may result in spurious oscillations in the solution {\cite{BrenS2002, ForM2013, YooH2018}}.   {We denote the discrete spaces $\mathbf{X}_h = \left[\mathcal{Q}_k\right]^2 \, \cap \, \mathbf{X}$ and ${Q_h} = \mathcal{Q}_k \, \cap \, Q$}. The chosen discrete space pair are {the} Taylor-Hood pair, $\mathbf{X}_h  \times Q_h$, {which} satisfies the \textit{discrete inf-sup} or Ladyzhenskaya–Babu$\check{s}$ka–Brezzi (LBB) compatibility condition
\begin{equation}\label{lbb}
\inf_{{p_{h}\neq0},\;p_{h} \in Q_{h}} \; \sup_{{\mathbf{u}_h\neq\mathbf{0}},\;\mathbf{u}_h \in  \mathbf{X}_{h}} \; \frac{\left( p_{h}, \; \nabla \cdot \mathbf{u}_{h}  \right)}{|| p_{h} ||_Q \, ||\nabla \mathbf{u}_{h} ||_\mathbf{X}} \geq c_0 >0,
\end{equation}
where constant $c_0$ is independent of $h$. A more detailed discussion and proof of the LBB condition and its proof can be found in \citep{ern2013theory,girault2012finite}. \\

The discrete finite element formulation for our problems is \\

\noindent {\textbf{Discrete weak formulation:}\textit{Find $\mathbf{u}_h^{n+1} \in \mathbf{X}_h$ and $p_h^{n+1} \in Q_h$ with $\mathbf{u}_h^{n+1} |_{\partial\Omega} = \widetilde{\mathbf{g}}$,  given $\left( \mathbf{u}_h^0, \; p_h^0 \right) \in \mathbf{X}_h \times Q_h$, and for each $n = 0, \, 1, \, 2, \ldots$, 
such that}
\begin{subequations}\label{fe_formulation0}
\begin{align}
&\widetilde{a}(\mathbf{u}_h^n; \; \delta \mathbf{u}_h^n, \; \mathbf{v}_h)  - b(\mathbf{v}_h, p_h^n) = \left( \mathbf{f}_h^n, \; \mathbf{v}_h  \right), \\
&b(\delta \mathbf{u}_h^n, \, q_h) = - b( \mathbf{u}_h^n, \, q_h),
\end{align}
\end{subequations}
\textit{where $\forall \; [ \mathbf{v}_h, \; q_h ]  \in \mathbf{X}_{0, \, h} \times Q_{0, \, h}$.  Also, $\widetilde{\mathbf{g}}$ is the projection of $\mathbf{g}$ onto finite element mesh.}}\\ 

{For the above} finite element discretization that is {the} inf-sup stable, we consider a stabilization similar to the one studied in \cite{olshanskii2009grad,heister2013efficient,ahmed2017grad,zhao2018numerical,benzi2006augmented}, 
resulting in a modified bilinear form
\begin{equation}\label{stablization_eq}
{a}(\mathbf{u}_h^n; \; \delta \mathbf{u}_h^n, \; \mathbf{v}_h) := \widetilde{a}(\mathbf{u}_h^n; \; \delta \mathbf{u}_h^n, \; \mathbf{v}_h)  + \left( \gamma \, \nabla \cdot \delta \mathbf{u}_h^n, \;  \nabla \cdot  \mathbf{v}_h^n\right).
\end{equation}
Here, $\gamma \geq 0$ is a stabilization parameter, 
{which can also be considered} as a piece-wise constant non-negative function concerning the finite element partition $\mathpzc{T}_{h}$. {Since 
we may have $\nabla \cdot \mathbf{u}_h^n \neq 0$ in many finite element discretizations; it can also be seen as adding a constant term to the momentum equation, which helps the penalty term in the mass conservation.} {Note that $\gamma$ can be taken as a constant per each mesh element, but in this work, it is taken as a global constant}.
\begin{remark}
{The optimal choice for $\gamma$ depends on the behavior of the solution on the finite element mesh $\mathpzc{T}_{h}$.  Many of the studies concerning the implementation of \textit{inf-sup} stable finite elements for the Oseen problem indicated a stable choice as $\gamma = \mathcal{O}(1)$ \cite{jenkins2014parameter,matthies2009some}. For sufficiently smooth solutions, $\gamma \sim 1$ should work to stabilize the solution. 
In our work, {we simply 
tested with $\gamma =0$ and $\gamma =1$, and from which 
we figured out no convergence obtained from $\gamma =0$} and selected $\gamma =1$. %An exhaustive 
A comprehensive study of the effect of this parameter is not the aim of the present contribution but that of the future one. }
\end{remark}

\noindent Along with the stabilization  term, the final discrete finite element problem is as follows: 
\\{\bf Discrete weak formulation: }\textit{Find $\mathbf{u}_h^{n+1} \in \mathbf{X}_h$ and $p_h^{n+1} \in Q_h$ with $\mathbf{u}_h^{n+1} |_{\partial\Omega} = \widetilde{\mathbf{g}}$,  given $\left( \mathbf{u}_h^0, \; p_h^0 \right) \in \mathbf{X}_h \times Q_h$, and for each $n = 0, \, 1, \, 2, \ldots$, 
such that}
\begin{subequations}\label{fe_formulation1}
\begin{align}
&{a}(\mathbf{u}_h^n; \; \delta \mathbf{u}_h^n, \; \mathbf{v}_h)  - b(\mathbf{v}_h,  \delta p_h^n) = \left( \mathbf{f}_h^n, \; \mathbf{v}_h  \right), \\
&b(\delta \mathbf{u}_h^n, \, q_h) = - b( \mathbf{u}_h^n, \, q_h),
\end{align}
\end{subequations}
and  {$\forall \; [ \mathbf{v}_h, \; q_h ]  \in \mathbf{X}_{0, \, h} \times Q_{0, \, h}$.} The bilinear form and the linear form are as follows:  
\begin{align}\label{FLHS-a}
\begin{split}
a\left( \mathbf{u}_h^n; \; \delta \mathbf{u}_h^n, \; \mathbf{v}_h \right) := &  \left(\mathbf{u}_h^{n}\cdot \nabla \delta \mathbf{u}_h^{n}, \, \mathbf{v}_{h}  \right) +  \left(  \delta \mathbf{u}_h^{n}\cdot \nabla \mathbf{u}_h^{n}, \, \mathbf{v}_{h}  \right)  + \frac{1}{Re} \left( \nabla  \delta \mathbf{u}_h^{n}, \,  \nabla   \mathbf{v}_h \right)   \\
&  + \frac{1}{Re \, Da} \left( \delta \mathbf{u}^{n}_{h}, \,  \mathbf{v}_{h}\right)  \\
&+  \frac{c_F}{\sqrt{Da}} \left[ \left( \| \mathbf{u}_h^n  \| \, \delta u_{h}^{n}, \, \mathbf{v}_h \right) +  \left( \left( \frac{ \mathbf{u}_{h}^{n} \cdot \delta \mathbf{u}_{h}^{n}}{\| \mathbf{u}_{h}^{n} \|} \right) \mathbf{u}_{h}^{n} , \, \mathbf{v}_h  \right)  \right] \\
& +\gamma\left( \nabla \cdot \delta \mathbf{u}_h^{n} , \, \nabla \cdot \mathbf{v}_h  \right),
\end{split}
\\
\label{FLHS-b}
b\left(\mathbf{v}_h,  \; \delta  p_h^n\right) := & \left( \delta p_{h}^{n}, \, \nabla \cdot   \mathbf{v}_{h} \right), &&
\end{align}
and 
\begin{align}\label{FRHS}
(\mathbf{f}_h^n, \, \mathbf{v}_h) := & -  \left( \mathbf{u}_h^n \cdot  \nabla \,\mathbf{u}_h^n, \,   \mathbf{v}_h  \right) + \left( p_{h}^{n}, \, \nabla \cdot   \mathbf{v}_h   \right)  -  \frac{1}{Re} \left(  \nabla \mathbf{u}_h^n, \, \nabla \mathbf{v}_h \right) - \frac{1}{Re \, Da}  \left(  \mathbf{u}_h^n, \, \mathbf{v}_h\right) && \\ \nonumber
& - \frac{c_F}{\sqrt{Da}} \left(  \| \mathbf{u}_h^n  \| \, \mathbf{u}_h^n, \, \mathbf{v}_h \right) 
   - \gamma \left( \nabla \cdot \mathbf{u}_h^n, \,  \nabla \cdot \mathbf{v}_h \right).&&
\end{align}
%The existence and uniqueness of a solution to the discrete problem \eqref{fe_formulation1} is guaranteed if the discretization satisfies the \textit{discrete inf-sup} (or LBB) compatibility condition with the existence of a positive constant $c_0$ as in \eqref{lbb}.  
\begin{remark}
Regarding the choice of the linear solver for the above linear system of equations, an important point to emphasize 
is that the old solution $\mathbf{u}_h^n$ represents the \textit{interpolant} of $\mathbf{u}_h$ in the discrete space $\mathbf{X}_h$, hence this term will not be point-wise divergence-free. 
The corresponding matrix obtained during the assembly of the integral $\left(\mathbf{u}_h^{n}\cdot \nabla \delta \mathbf{u}_h^{n}, \, \mathbf{v}_{h}  \right)$ will not be symmetric; therefore, {the matrix obtained for \eqref{FLHS-a} overall} is not symmetric. The issue of choosing a suitable solver will be addressed in the following section. 
\end{remark}

\subsection{Linear solver and preconditioner}\label{sec:preconditioner}
{In this paper, one of our main goals} is to 
establish {an efficient linear solver} for the stable solution of the nonlinear {Darcy-Brinkman-Forchheimer} model. 
{In establishing competent linear algebra tools}, it is convenient to ascertain the discrete problem \eqref{fe_formulation1} as a matrix 
form in terms of actual finite element matrices. We need to write the variables in terms of {finite element basis functions} to obtain such matrices. Let the basis define the finite-dimensional functional space sets as
\begin{equation}
\mathbf{X}_h := \mbox{span}\left\{  \boldsymbol{\phi}_i            \right\}_{i=1}^{n_\mathbf{u}}, \quad
Q_h := \mbox{span}\left\{  {\psi}_j            \right\}_{j=1}^{n_p}, 
\end{equation}
where $\phi$ and $\psi$ are scalar basis functions for velocity and pressure, respectively. Also, $n_\mathbf{u}$ and $n_p$ are the number of degrees of freedom for velocity and pressure variables, respectively. The associated discrete solutions, {$\delta \mathbf{u}_h$ and $\delta p_h$ for $\mathbf{u}_h^{n+1} \in \left( \mathbb{R}^{n_\mathbf{u}}\right)^2$ and $p_h^{n+1} \in \mathbb{R}^{n_p }$ expressed with $\mathbf{u}_h^{n+1}=\mathbf{u}_h^{n} + \delta \mathbf{u}_h^{n}$ and ${p}_h^{n+1}={p}_h^{n} + \delta {p}_h^{n}$,} are obtained from a linear combination of basis functions on a finite element mesh as
\begin{equation}
\delta \mathbf{u}_h := \sum_{i=1}^{n_\mathbf{u}} \mathbf{u}_i \, \boldsymbol{\phi}_i, \quad \mbox{and} \quad \delta p_h :=  \sum_{j=1}^{n_p} {p}_j \, {\psi}_j.
\end{equation}
The entries of {the associated system matrix composed of sub-matrices (expressed in \textit{capital letter} below) corresponding to the ones defined in} 
\eqref{lforms} respectively are given by
\begin{subequations}\label{lforms2}
\begin{align}
&M_{ij} := \left( m (\nabla \delta \mathbf{u}_{h}^{n}, \; \nabla \mathbf{v}_h) \right)_{ij} =  \dfrac{1}{Re}\left( \nabla \boldsymbol{\phi}_i, \;   \nabla \boldsymbol{\phi}_j   \right),\\
&C_{ij} := \left( c ( \delta \mathbf{u}_{h}^{n}, \;  \mathbf{v}_h) \right)_{ij} =  \dfrac{1}{Re \, Da}\left(  \boldsymbol{\phi}_i, \;    \boldsymbol{\phi}_j   \right), \\
&N_{ij} :=  \left( n (\mathbf{u}_h^n; \, \delta \mathbf{u}_{h}^{n}, \;  \mathbf{v}_h) \right)_{ij} =  \left( \mathbf{u}_{h}^{n} \cdot  \nabla \boldsymbol{\phi}_i, \;    \boldsymbol{\phi}_j   \right)   + \left(\boldsymbol{\phi}_i \cdot  \nabla \mathbf{u}_{h}^{n},   \;    \boldsymbol{\phi}_j   \right), \\
&S_{ij} :=  \left( s (\mathbf{u}_h^n; \, \delta \mathbf{u}_{h}^{n}, \;  \mathbf{v}_h) \right)_{ij} = \dfrac{c_F}{\sqrt{Da}} \left[ \left( \| \mathbf{u}_{h}^{n} \|   \boldsymbol{\phi}_i, \;    \boldsymbol{\phi}_j   \right)   + \left( \dfrac{\mathbf{u}_h^n \cdot \boldsymbol{\phi}_i}{\| \mathbf{u}_h^n \|} \; {\mathbf{u}_h^n},   \;    \boldsymbol{\phi}_j   \right) \right],
\end{align}
\end{subequations}
and the entries of the divergence matrix are given by
\begin{equation}
B_{ij} := -\left( b(\mathbf{v}_h, p_h^n) \right)_{ij} =  - \left( \psi_i, \; \nabla \cdot  \boldsymbol{\phi}_j \right),
\end{equation}
{and finally, the matrix entries from the added ``stabilized'' Augmented Lagrangian term are given by
\begin{equation}\label{eq:L}
{L_{ij} := \left( l( \mathbf{u}_h^n,\,  \mathbf{v}_h) \right)_{ij}  =\left(  \nabla \cdot  \boldsymbol{\phi}_i, \,  \nabla \cdot  \boldsymbol{\phi}_j \right).}
\end{equation}
}
{Then, the assembled linear system of the discretized system results in} 
\begin{equation}\label{linsystem1}
{\mathcal{A}} \mathbf{x} = \mathcal{G},
\end{equation}
where $\mathbf{x} = \left(\mathbf{U}, \; \mathbf{P} \right)^{T}$, {and the vectors $\mathbf{U}$ and $\mathbf{P}$} are the discrete version {of vectors} that belong to the finite element degrees of freedom of the variables, i.e., velocity $\mathbf{u}_h$ and pressure $p_h$. {{Right-hand side vector is $\mathcal{G} =  \left(\mathbf{F}, \; \mathbf{0} \right)^{T}$}, where the vector $\mathbf{F}$ is obtained by assembling the Equation~\ref{FRHS}}. The 
structure of the left-hand side block matrix is
\begin{equation}\label{bmatrix1}
\mathcal{A} = \begin{pmatrix}
\mathbf{A} &  & &\mathbf{B}^{T} \\
\mathbf{B} & &  &\mathbf{0}
\end{pmatrix},
\end{equation}
{where} the block matrices {$\mathbf{A} =\left(\mathbf{M} + \mathbf{C} +\mathbf{N} +\mathbf{S}\right)$} and $\mathbf{B}$, which correspond to the bilinear forms $\widetilde{a}$ and $b$ in \eqref{fe_formulation0}, respectively, are based on the sub-matrices described in \eqref{lforms2}. The sparse linear system \eqref{linsystem1}-\eqref{bmatrix1} are known in the literature as \textit{saddle point {problem\;(system)}}, and a great deal of work can be found about developing stable and efficient solvers for such a system (e.g., {%including 
the classical {Navier-Stokes} system) via using $\mathbb{Q}_2\mathbb{Q}_1$ pair, so-called the \textit{Taylor-Hood} %$\mathbb{Q}_2\mathbb{Q}_1$ pair, 
finite elements.}

To stabilize the linear solver for the nonlinear model \eqref{eq:Darcy-Brinkan-Forch} or \eqref{eq:Darcy-Brinkan-Forch-1} in this work, we then use the classical \textit{Augmented Lagrangian} approach {as aforementioned in \eqref{stablization_eq}}. In this method, we replace the above linear matrix system \eqref{linsystem1} with the equivalent {system}
{
\begin{equation}\label{bmatrix2}
\widetilde{\mathcal{A}} \mathbf{x} = \mathcal{G},
\end{equation}
{where  
\begin{equation}\label{bmatrix3}
\widetilde{\mathcal{A}} =
\begin{pmatrix}
\mathbf{A}_{\gamma} &  & &\mathbf{B}^{T} \\
\mathbf{B} & &  &\mathbf{0}
\end{pmatrix}.
\end{equation}}
{For the regularized block matrix $\mathbf{A}_{\gamma} = \mathbf{A}+\gamma \, \mathbf{B}^T \, \mathbf{W}^{-1} \,\mathbf{B}$, i.e., the $(1, 1)$ block of $\widetilde{\mathcal{A}} $ in \eqref{bmatrix3},} $\mathbf{W}$ is {a positive definite matrix which is invertible.}} 
 {{However,} a notable issue {lies in} 
the choice of the parameter $\gamma$ 
{such that the block matrix  $\mathbf{A}_{\gamma}$ could become} ill-conditioned for large $\gamma$.} Thus, it is desirable
 to keep  the value of $\gamma$ balanced such that it satisfies $\gamma  \approx O(1)$ or specifically $\gamma \approx O( \| \mathbf{u}_{h} \|)$ \cite{benzi2006augmented}; therefore devising an effective block-solver 
becomes challenging. 

{In this work,} we solve the system by applying a {suitable preconditioning as follows.} {With an operator $\mathcal{P}$}, which is a block triangular matrix as described in \cite{saad2003iterative,olshanskii2009grad,heister2013efficient}, {we} obtain the solution by solving 
\begin{equation}
\widetilde{\mathcal{A}} \mathcal{P}^{-1} \mathbf{y} = \mathcal{G},
\end{equation} 
where 
\begin{equation}
\mathbf{x} = \mathcal{P}^{-1} \mathbf{y}.
\end{equation}
Here, we define the block triangular preconditioner operator $\mathcal{P}$ as
\begin{equation}
\mathcal{P} =  \begin{pmatrix}
\mathbf{A}_{\gamma}  &  & &\mathbf{B}^{T} \\
\mathbf{0} & &  & \widetilde{\mathbf{S}}
\end{pmatrix}, 
\end{equation}
then the corresponding  inverse block preconditioner is
\begin{equation} \label{bpreconditioner}
\mathcal{P}^{-1} :=  \begin{pmatrix}
\mathbf{A}_{\gamma}  &  & &\mathbf{B}^{T} \\
\mathbf{0} & &  &\widetilde{\mathbf{S}}
\end{pmatrix}^{-1} 
=  \begin{pmatrix}
\mathbf{A}_{\gamma}^{-1}  &  & & \mathbf{0} \\
\mathbf{0} & &  & \mathbf{I}
\end{pmatrix} \begin{pmatrix}
\mathbf{I}  &  & & \mathbf{B}^{T}  \\
\mathbf{0} & &  & -\mathbf{I}
\end{pmatrix} \begin{pmatrix}
\mathbf{I}  &  & & \mathbf{0}  \\
\mathbf{0} & &  & - \widetilde{\mathbf{S}}^{-1}
\end{pmatrix},
\end{equation}
where {$\mathbf{I}$ is the rank-2 identity tensor,} $\widetilde{\mathbf{S}}$ is a Schur complement, i.e.,  
 $\widetilde{\mathbf{S}} = \mathbf{B} \, \mathbf{A}_{\gamma}^{-1} \, \mathbf{B}^{T}$, and it is implicitly defined by
\begin{equation}\label{scomplement1}
\widetilde{\mathbf{S}}^{-1} := - \nu \; \mathbf{M}_{p}^{-1} - \gamma \; \mathbf{W}^{-1},
\end{equation}
where $\mathbf{M}_p$ denotes the ``pressure mass matrix'' and its components are
\begin{equation}
\left( \mathbf{M}_p \right)_{ij} = \int \int \psi_i \, \psi_j \; dx \, dy.
\end{equation}
{The term $\nu$ in Equation~\ref{scomplement1} is the viscosity for which we have assumed %it 
as globally constant in this work.} A main motivation  for the definition  of $\widetilde{\mathbf{S}}$ is largely {due to} the fact that 
\begin{align}
\widetilde{\mathbf{S}}^{-1} &=\left[\mathbf{B} \left(\underbrace{{\mathbf{A} + \gamma \mathbf{B}^T \mathbf{W}^{-1} \mathbf{B}}}_{\displaystyle \mathbf{A}_{\gamma}} \right)^{-1} \mathbf{B}^T \right]^{-1}  \notag \\
&= \left( \mathbf{B} \mathbf{A}^{-1} \mathbf{B}^T\right)^{-1} + \gamma \; \mathbf{W}^{-1},
\end{align}
{for which a similar argument can be found in Lemma $4.1$ {in} \cite{benzi2006augmented}}, 
and we let $\mathbf{W} = \mathbf{M}_p$ in 
{the} computations. To build the Schur complement, which is the inverse of the velocity block, 
$\mathbf{A}^{-1} $ needs to be solved  and then a matrix-vector multiplication with $\mathbf{B}^T$ is performed. {We use UMFPACK \cite{davis2006direct,davis2004algorithm} direct solver to compute $\mathbf{A}_{\gamma}^{-1}$} and CG with ILU {preconditioner} to build $\widetilde{\mathbf{S}}^{-1}$ \cite{dealII91}. Finally, to solve the linear matrix system, $\widetilde{\mathcal{A}}\mathbf{x} = \mathcal{G}$,  we use the {generalized} minimal residual method with flexible preconditioning (flexible GMRES or FGMRES) \cite{saad1993}. 

\begin{remark}
{Particularly when computing the Schur complement, globally constant viscosity $\nu$ and the Augmented Lagrangian parameter $\gamma$ {have been chosen in this work.} However, these parameters could be made variable, and meaningful modified pressure mass matrix elements can also be computed via:
\begin{equation}\label{eq:Schur}
\left( \mathbf{M}_p(\nu, \; \gamma)\right)_{ij} = \int_{\Omega} \left( \nu(\mathbf{x}) + \gamma(\mathbf{x}) \right)^{-1} \phi_i  \, \phi_j \; d\Omega.
\end{equation}} 
\end{remark}

\subsection{{Local adaptive mesh refinement}}\label{sec:LGR}
It is well-known that the finite element solvers require finer mesh for accurate solution \cite{gresho1981,verfurth1994posteriori}. {For the local grid refinement, {the \textit{h-version} mesh refinement indicators} developed by Kelly et al. in \cite{kelly1983} and De S.R. Gago et al. in  \cite{de1983posteriori} are employed in this study,  which were are originally derived for the second-order, advection-type, linear problems.} The \textit{Kelly error estimator} is based on {the integration of 
the gradient of the solution} along the faces of each cell. It has been proven effective in detecting the regions of high gradients and jumps in scalar flux between elements. In a nutshell, the computation of this indicator is simple {as}
\begin{equation}\label{kelly}
\eta_{K} = \frac{h}{24} \int_{\partial K}  \llbracket \partial_{n} {q}_{h}   \rrbracket^{2} \, ds,
\end{equation}
where $\llbracket \partial_{n} {q}_{h}   \rrbracket$ is the jump in the normal derivative of the finite element solution of ${q}_{h}$ across each element boundary $\partial K$. In our actual implementation, after computing the indicator $\eta_{K}$ for each cell, 
it is resolved which cells must be refined or coarsened via assigning either a refinement or a coarsening flag to each cell. For the computational efficiency,  we refine  
$30~\%$ of cells 
{from} the highest error indicator and coarsen those cells within $3~\%$ of elements  
{from} the lowest error using \eqref{kelly}. 
\begin{remark}\label{rmk:local_refinement}
 {The main} purpose of the mesh refinement indicators is to provide the solution error in a norm that can be computable {on} a specified finite-dimensional function space if the solution and all the necessary boundary data are available. However, in many instances, no methods are available to find the exact solution to the nonlinear {partial differential equations}; hence, the actual error is still being determined. A ``good'' error estimator should be accurate, because as the mesh length $h$ tends to zero, the error estimate should also tend to zero at the same rate. The error estimator should also be robust for practical nonlinear problems, and implementation should be possible. {In this sense, the use of \textit{h-adaptivity}, without remeshing the whole domain and without the necessity of doing a trail-and-error procedure, facilitates the solution to be smooth at the ``wiggly''  places  
{related} to high gradients.} 
{By  
generating a ``suitable'' mesh using \textit{a posteriori} error estimate, 
the advantages {of local \textit{h-adaptivity}} is clear: low computational cost with less effort on the regions of the vanishing gradient.}
 Meanwhile, a significant bottleneck in the quest for a mesh refinement indicator is that such an ideal estimator may only be available for some situations. For more information on finite element error estimators and additional prospects and challenges, interested readers can see \cite{zienkiewicz1991adaptivity,usmani1999solution,rios2011h,holst2016convergence,bangerth2013adaptive,ainsworth2011posteriori}, and the references therein.  
\end{remark}

{Finally}, we summarize our approach for {numerical solutions of the} nonlinear Darcy-Brinkman-Forchheimer model with a pseudo-algorithm, and 
the steps  
involved are as {follows:}

\begin{itemize}
\item[(i)] \textit{\textbf{{Domain and boundary conditions}}} \quad {{ 
For a  
computational domain of our interest, set the boundary conditions for velocity} (e.g., the moving top lid and zero for the remaining parts).}
{\item[(ii)] \textit{\textbf{Mesh Adaptation}}\quad Refine the mesh using {Kelly error indicator} \eqref{kelly}, and interpolate the old solution to the new mesh as the mesh refinement adds new nodes to the old mesh.} %{\hl{mesh adpatation before Newton's}}
  
\item[(iii)] \textit{\textbf{{Assembly}}} \quad Set up the matrices %$\widetilde{\mathbf{A}}$ 
$\mathbf{A}$, $\mathbf{B}$ and the right-hand side vector $\mathbf{F}$ as defined in \eqref{bmatrix1} from the weak formulation \eqref{fe_formulation1}, and create a linear system $\widetilde{\mathcal{A}}\mathbf{x} ={\mathcal{G}}$ \eqref{bmatrix2}, {where the stablitzation term $\mathbf{A}_{\gamma}$ is included}.
{\item[(iv)] \textit{\textbf{Initial Guess}} {Find a pair of initial guess $(\mathbf{u}_{h}^{0}, \; p_{h}^{0})$ with the Darcy-Brinkman {model}~\eqref{eq:Darcy-Brinkman} using the nonlinear solver {based on the Newton's method as follows}.}}  
\item[(v)] \textit{\textbf{Preconditioner and Solver}}\quad Set up the block preconditioner {as in} \eqref{bpreconditioner} and corresponding Schur complement as in \eqref{scomplement1}. {There are three solvers in total}: first for block matrix with {the Augmented Lagrangian term} $A_{\gamma}^{-1}$, second for the pressure mass matrix $M_p^{-1}$, and third one is for solving the linear system $\widetilde{\mathcal{A}}\mathbf{x} = {\mathcal{G}}$. The first two solvers are used inside the preconditioner 
and the third one {(i.e., FGMRES)} is invoked at the final step of the solution procedure. 
%{\item[(iv)] \textit{\textbf{Mesh Adaptation}}\quad Refine the mesh using {Kelly error indicator} \eqref{kelly}, and interpolate the old solution to the new mesh as the mesh refinement adds new nodes to the old mesh.} %{\hl{mesh adpatation before Newton's}}
%\item[(v)] \textit{\textbf{Initial Guess}} {Find a pair of initial guess $(\mathbf{u}_{h}^{0}, \; p_{h}^{0})$ with the Darcy-Brinkman {model}~\eqref{eq:Darcy-Brinkman} using the nonlinear solver {based on the Newton's method as follows}.} 
\item[(vi)] \textit{\textbf{{Newton's Method}}} \quad 
Compute the \textit{Newton updates} for velocity $\delta \mathbf{u}_{h}^n$ and the pressure $\delta p_{h}^n$, and finally find the solution updates using

\[
\mathbf{u}_{h}^{n+1} = \mathbf{u}_{h}^n +{\delta \mathbf{u}_{h}^n}, \quad \mbox{and} \quad  p_h^{n+1} = p_{h}^n + {\delta p_{h}^n}.
\]
%\item[(vii)] 
%\quad{\textbf{$\circ$ Stopping Criterion}}\quad 
For stopping criterion, compute the $L^2$-norm of the right hand-side vector $F(\mathbf{u}_{h}^{n+1}, \; p_h^{n+1})$ as in \eqref{rhs}, and let 
\[
R^{n+1} = \| F(\mathbf{u}_{h}^{n+1}, \; p_h^{n+1}) \|,
\]
and finally check with the tolerance, $\varepsilon$, for the \textit{Newton's method} {with FGMRES} as follows:
    \begin{itemize}
      \item[a:] if $R^{n+1}  \leq \varepsilon$, then STOP;
      \item[b:]  if $R^{n+1}  \geq \varepsilon$, then go to step {(v)}.
    \end{itemize}
\end{itemize}

\section{Numerical Experiments}\label{num_exp}

{In this section}, we perform several numerical tests %to substantiate 
for the efficacy of %the 
proposed numerical discretization  and solver of the 
 Darcy-Brinkman-Forchheimer nonlinear fluid model. 
{We solve the Darcy-Brinkman-Forchheimer model using the Taylor-Hood  $\mathbb{Q}_2\mathbb{Q}_1$ elements for velocity and pressure, respectively. %Meanwhile, we also aim to validate the theoretical approach of Grad-Div stabilization and convergence of the finite element solution.  
For our computations, we follow the very step of the pseudo-algorithm right above in the previous section.
The stopping criterion, $\varepsilon$, of the Newton method for {outer linear solver using FGMRES} is set as $\varepsilon=1.0e$-$12$.} 

{For the numerical tests, we start with the $h$-convergence study  
%is first %performed  
with a manufactured solution, 
%This 
which is followed by a benchmark problem in the Navier-Stokes system, %and a comparison study between the Darcy-Brinkman and Darcy-Brinkman-Forchheimer models, 
all for the verification of the code.
The code is developed  
by authors based on an open source finite element library, called the \texttt{deal.II}~\cite{dealII91}, where the tutorial \texttt{step-57} \cite{ZhaL2017} is the foundation for the current implementation.} 

%Finally, 
%We then perform several tests with parameter study %are presented 
%to demonstrate some distinct features of the model, highlighting the role of Forchheimer term.}  %{compared to models without it}.}
{
After validating the theoretical approach of Grad-Div stabilization and convergence of the finite element solution, 
we choose a standard  
problem of the lid-driven cavity flow inside a unit square under specific Dirichlet boundary conditions. We then perform several tests with parameter study %are presented 
to demonstrate some distinct features of the model, highlighting the role of Forchheimer term.}  %{compared to models without it}.}

\subsection{{\textit{h}-convergence study}}
%\hl{In this subsection,} 
{We} first perform a $h$-convergence study to verify the proposed algorithm and its corresponding implementation. To this end, we choose a ``manufactured'' solution: 
\begin{equation}\label{m_solution}
\mathbf{u}=(\sin({\pi x}), -\pi y\cos({\pi x})), \quad  \mbox{and} \quad p=\sin({\pi x})\cos({\pi y}).
\end{equation}
Note that the above choice for $\mathbf{u}$ is divergence-free, automatically satisfying the continuity equation. 
For simplicity, all other parameters in the {Darcy-Brinkman-Forchheimer model} model are taken as $1.0$. 
For the numerical implementation, we choose a square domain and the boundary conditions using the velocity vector given in  \eqref{m_solution}. 
%A total %of 
Total $5$ %in 
global refinements %in total 
are performed for the domain, and the numerical solution with  corresponding errors in  {$L^2$-} and $H^1$-norms for the velocity and pressure are recorded. 
\begin{table}[htbp]
\centering
\caption{{Errors in {$L^2$-} and $H^1$-norms for the velocity ($\mathbf{u}$) and in $L^2$-norm for the pressure ($p$) with the manufactured solution and their convergence rates.}} 
\label{tab:h-conv}
\begin{tabular}{|c|| l|l|| l|l|| l|l|}
\hline
\multirow{2}{*}{DoFs}  & \multicolumn{2}{c||}{velocity $\mathbf{u}$}    & \multicolumn{2}{c||}{velocity $\mathbf{u}$}    & \multicolumn{2}{c|}{pressure $p$}                            \\ \cline{2-7} 
&  \multicolumn{1}{c|}{$L^2$ Error}  & \multicolumn{1}{c||}{Rate} & \multicolumn{1}{c|}{$H^1$ Error} &  \multicolumn{1}{c||}{Rate} &  \multicolumn{1}{c|}{$L^2$ Error}  & \multicolumn{1}{c|}{Rate} \\ \hline\hline
59    & %0.027440240607  
2.744e-2 & {-}   & %0.415317177773 
4.153e-1 & {-}    & %0.105965606868  
1.059e-1 & {-}           \\ \hline
187    &  %0.003405715572  
3.405e-3 & 8.0571   & %0.105247445405 
1.052e-3 & 3.9461  & %0.017805930227    
1.780e-2 & 5.9511   \\ \hline
659   & %0.000426246057  
4.262e-4 & 7.9900  & %0.026409855112 
2.640e-2 & 3.9852  & %0.004143296741  
4.143e-3 & 4.2975             \\ \hline
2467  &  %0.000053322387  
5.332e-5 & 7.9938   & %0.006608675234    
6.608e-3 & 3.9962  & %0.001020679949   
1.020e-3 & 4.0593 \\ \hline
9539  & %0.000006666931   
6.666e-6 & 7.9980  &  %0.001652557054 
1.652e-3 & 3.9991 &  %0.000254284008    
2.542e-4 & 4.0139   \\ \hline
\end{tabular}
\end{table}

The convergence rates in $H^1$- and $L^2$-norms for the discretization using Taylor-Hood  $\mathbb{Q}_2\mathbb{Q}_1$ elements for velocity and pressure are  
presented in Table~\ref{tab:h-conv}. It is clear from the above table that for bi-quadratic basis functions, the velocity has a rate of $8$ in  $L^2$-norm and a rate of $4$ in $H^1$-norm. For the linear elements, the pressure has a rate of $4$ in $L^2$-norm. Such convergence rates are optimal for the  $\mathbb{Q}_2\mathbb{Q}_1$ pair of finite elements used in the approximation \cite{girault2012finite}. 

\subsection{{Benchmarking problem of lid-driven cavity flow}}\label{sec:Benchmark}
This section further verifies our algorithm and its implementation by %comparing {two 
{a benchmarking problem with the steady lid-driven cavity flow}. %issues. 
The computational domain for %these problems 
the lid-driven cavity flow is the same 
unit square, 
and the top boundary {($\Gamma_{D_{T}}$)} only has the horizontal velocity  
of unity {(Figure \ref{Fig:CompDomain})}. {Meanwhile, the left ($\Gamma_{D_L}$),  right ($\Gamma_{D_R}$) and  bottom ($\Gamma_{D_B}$) boundaries, we have homogeneous boundaries with zero velocity.} 
Via obtaining stable and convergent numerical solutions for these examples, we verify our stabilized finite element discretization and other numerical methods employed in this study. %for the steady lid-driven cavity flow filled with the Darcy-Brinkman-Forchheimer porous medium.  

%\subsubsection{{Navier-Stokes system}}\label{sec:BM-NS}

{Before} addressing the full nonlinear Darcy-Brinkman-Forchheimer 
porous medium, we first test the lid-driven cavity flow within a pure Navier-Stokes system.  The Navier-Stokes  is equivalent to the Darcy-Brinkman-Forchheimer model when the permeability of Darcy-Brinkman-Forchheimer reaches infinity. To this end, we compare our results with the work of Ghia et al.~\cite{Ghia1982}.
Numerical results for the NS system will be further compared in the following section (Section~\ref{sec:parameter_study}). 

Figure~\ref{Fig:EX_Bench_Velocity_Center} compares the numerical results obtained from our algorithm with the results of Ghia et al.'s work for selected values of Reynolds number: $Re=1000$ and $Re=3200$. It is clear %from the figure 
that our numerical results are in excellent agreement. %with the results obtained in \cite{Ghia1982}. 

\begin{figure}[!h]
\centering
\includegraphics[width=0.3\textwidth]{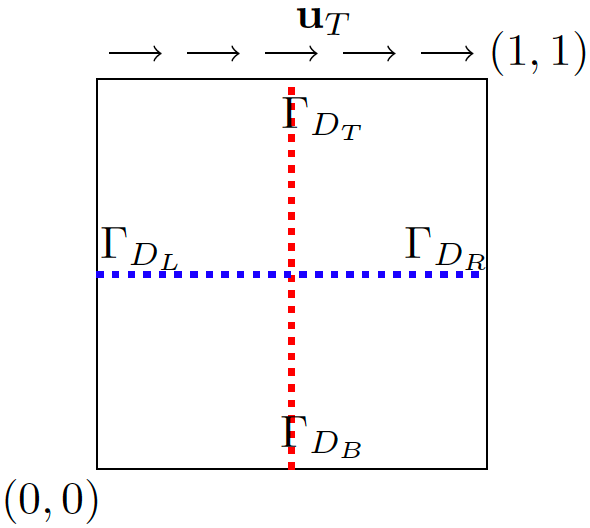}
\caption{An unit square and the boundary conditions for the computational domain. Homogeneous boundary conditions are imposed except for $\mathbf{u}_{T}=(1,0)$. Two dotted lines in color are the reference lines for monitoring velocity: the red for horizontal velocity ($u_x$) and the blue for vertical velocity ($u_y$). }
\label{Fig:CompDomain}
\end{figure}

\begin{figure}[!h]
\centering
\includegraphics[width=0.75\textwidth]{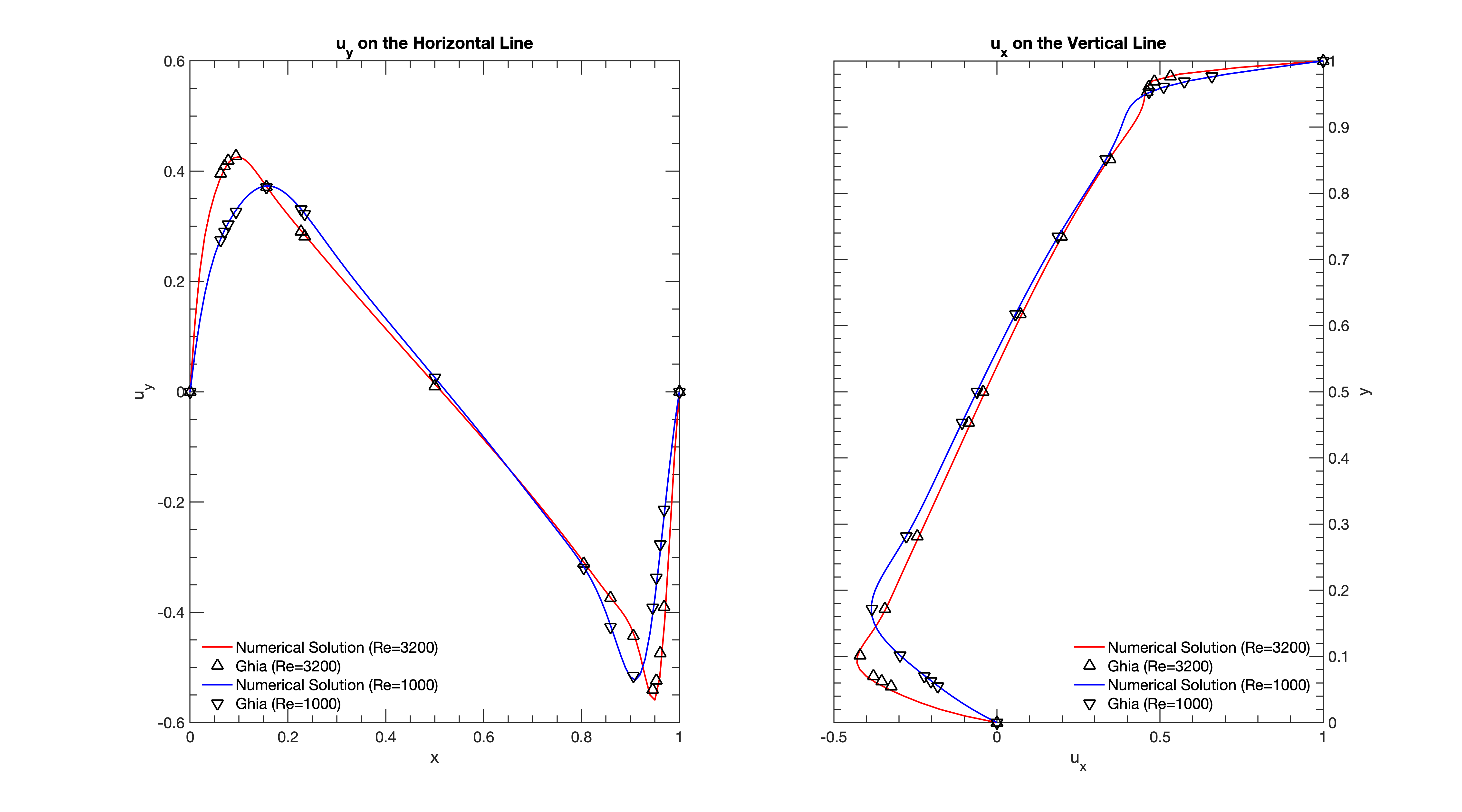}
\caption{\textbf{Benchmark {Example} 1}. Velocity on the reference lines for $Re = 1000$ (blue solid line) and $Re = 3200$ (red solid line): (left) $u_y$ on the horizontal line and (right) $u_x$ on the vertical line. Numerical results are compared with the solutions from \cite{Ghia1982}.}
\label{Fig:EX_Bench_Velocity_Center}
\end{figure}

\subsection{Model comparisons and parameter study}\label{sec:parameter_study}
{In this section,} we perform several numerical tests 
for comparison study between the models of linearized Brinkman \eqref{eq:Brinkman}, Darcy-Brinkman \eqref{eq:Darcy-Brinkman}, and Darcy-Brinkman-Forchheimer \eqref{eq:dimless-full-nonlinear}.  
Particularly for parameter study, the role of Forchheimer term in \eqref{eq:dimless-full-nonlinear} is highlighted.  
%As mentioned, the 
Furthermore, the Navier-Stokes model is also compared if necessary. {Henceforth,} these models may be abbreviated as the D-B (Darcy-Brinkman), D-B-F (Darcy-Brinkman-Forchheimer), and N-S (Navier-Stokes).

The computational domain and boundary conditions for the lid-driven cavity filled with porous medium are depicted in Figure~\ref{Fig:CompDomain}. % the same as the previous one 
%and the details are what follows. 
%{As depicted in Figure}~\ref{Fig:CompDomain}, 
{For the steady lid-driven cavity flow,} we have the Dirichlet boundary conditions in the absence of any sliding factor 
on the top boundary ($\Gamma_{D_{T}}$), {but we have $\mathbf{u}_{T}=(1,0)$ on $\Gamma_{D_{T}}$.}
For the rest boundaries, {the same boundary conditions are applied as in Section~\ref{sec:Benchmark}}, i.e., the homogeneous boundaries. %the left ($\Gamma_{D_L}$),  right ($\Gamma_{D_R}$) and  bottom ($\Gamma_{D_B}$) boundaries, we have homogeneous boundaries with zero velocity. 
Note that the dotted reference lines passing the center are for %monitoring 
comparing the velocity as the benchmarking example in the previous section 
(see Figure~\ref{Fig:EX_Bench_Velocity_Center}); $u_x$ is monitored from the vertical red, 
while $u_y$ is monitored from the horizontal blue.
{The domain starts with total 5 refinements globally, and perform its refinement cycle using the {Kelly} error estimator upto the max refinement number, for which we set 4 in this study.}

{Regarding the properties of porous medium, %For simple analysis, 
we fix porosity constant as $\phi=1.0$. Although the porosity is a variable depending upon the fluid pressure or volumetric strain of porous medium, %it is not the scope of this study, 
%we set it constant as unity for simple analysis so that 
by setting it constant as unity for simple analysis, %so that 
%and by setting it constant, %can eliminate its effect on the flow, and 
%we can purely 
%in order to 
we can purely investigate numerical aspects of the model and perform comparison study between the models. In the same manner, with an assumption of isotropic medium, intrinsic permeability ($K$) and inertial resistance ($\beta$) are taken as scalar coefficients instead of full tensor expressions.}
In order to test with the dimensionless parameters, i.e., $Re$ and $Da$ (see \eqref{eq:Re} and \eqref{eq:Da}) {primarily}, the reference (or the fundamental scale) parameters are set as Table~\ref{tab:Parameters2}, for which $\phi=1.0$ is presumed.
{Note that} these reference values are problems- or environments-specific for porous medium, and they can be obtained from the experimental data such as the grain size and seepage rate from lab experiments. For the inertial resistance coefficient ($c_F$), it can be determined either empirically (from the relation of permeability) or theoretically (such as using Ergun's equation \cite{YanD2012}). In this numerical study, we fix it as $0.5$ for our reference value, {unless otherwise noted}.

\begin{table}[!h]
\begin{center}
\caption{Scale parameters and inertial resistance coefficient that are selected in this study.}
\label{tab:Parameters2}
\begin{tabular}{cll}
\toprule
\textbf{Parameter} & \textbf{Value} & \textbf{Unit}\\
\midrule
Reference distance ($L$) & $0.001$ & $m$\\
Reference discharge ($U_{0}$) & $0.1$ & $m/s$\\ 
Inertial resistance ($c_F$) & {{0.5}} & -- \\
\bottomrule
\end{tabular}
\end{center}
\end{table}

{{For the parameter study}, three different values of Reynolds number are selected, i.e., $Re = 10, 100,$ and $1000$.} 
Following a typical classification of the flow regime in porous media~\cite{BearJ1972}, the selected $Re$ values covers the \textit{laminar} ($Re\leq10$), the \textit{transient} (or the \textit{nonlinear laminar}, $10<Re\leq100$), and the \textit{turbulent} ($100<Re\leq1000$) flow regimes. 
On the other hand, 
the intrinsic permeability $(K)$ 
is set with next six values: $2.5e$-$12$, $2.5e$-$11$, $2.5e$-$10$, $2.5e$-$7$, $2.5e$-$6$, and $2.5e$-$5$, with the unit of $m^2$.
%{In terms of `Darcy' unit, $2.5e$-$12$ is about 2.5 Darcy.} 
The former three can be categorized as relatively low permeability group compared to the latter three of relatively high permeability group. 
%and 
{Using \eqref{eq:Da} with Table~\ref{tab:Parameters2}, $Da$ can then be obtained.}  
With the combinations of these $Re$ and $Da$ numbers, there are 18 cases in total (see Table~\ref{tab:Examples}). For comparison purposes in this study, 
we categorize $Re\times Da$ values into two groups: {\texttt{GROUP I}  and \texttt{GROUP II}} with nine cases for each. 
We note \texttt{GROUP I} is 
for $Re\times Da < 1.0$ and \texttt{GROUP II} is for $Re\times Da \geq 1.0$. 
Accordingly, each $Re$ number owns six cases in total, where each group owns three examples. %: \texttt{GROUP I} and \texttt{GROUP II}. 
In Table~\ref{tab:Examples}, we denote
the numbers for \textbf{\textbf{\textsl{Test}}}  
corresponding to cases in  each group, 
following the order of  
$Re\times Da$ values from small to large. 

\begin{table}[!h]
\begin{center}
\caption{{Numerical tests  based on the values of  $Re$ and $Da$ values. The left three columns of \boldsymbol{$Da$} belong to \texttt{GROUP I}, i.e., \boldsymbol{$Re\times Da$} $(<1.0)$, while the right three of it belong to \texttt{GROUP II}, \boldsymbol{$Re\times Da$} $(>1.0)$.}}
\label{tab:Examples} %\texttt{GROUP I: }\boldsymbol{$Re\times Da$} $(<1.0)$ \texttt{GROUP II: }\boldsymbol{$Re\times Da$} $(\geq 1.0)$
\begin{tabular}{llcccccc}\toprule
\multicolumn{2}{c}{\multirow{2}{*}{\boldsymbol{$Re\times Da$}}} & \multicolumn{6}{c}{\boldsymbol{$Da$}} %& \multicolumn{3}{c}{\boldsymbol{$Da$}} 
\\
\cmidrule(lr){3-5}\cmidrule(lr){6-8}
%&& \multicolumn{3}{c}{} & \multicolumn{3}{c}{} \\
 & &  %$0.0000025$  
 2.5e-6 & %$0.000025$ 
 2.5e-5 & %$0.00025$    
 2.5e-4 & %$0.25$  
 2.5e-1 & %$2.5$ 
 2.5e+0 & %$25.0$
 2.5e+1 \\\midrule
 \multirow{6}{*}{\boldsymbol{$Re$}} & \multirow{2}{*}{$10.0$} & \textbf{\textbf{\textsl{Test}}} \textbf{1} & \textbf{\textbf{\textsl{Test}}} \textbf{2} & \textbf{\textbf{\textsl{Test}}} \textbf{3} & \textbf{\textbf{\textsl{Test}}} \textbf{1} & \textbf{\textbf{\textsl{Test}}} \textbf{2} & \textbf{\textbf{\textsl{Test}}} \textbf{3} \\
 & & %$0.000025$ 
 2.5e-5 & %$0.00025$ 
 2.5e-4 & %$0.0025$ 
 2.5e-3 & %$2.5$ 
 2.5e+0 & %$25.0$ 
 2.5e+1 & %$250.0$ 
 2.5e+2\\
 &  \multirow{2}{*}{$100.0$} & \textbf{\textbf{\textsl{Test}}} \textbf{4} & \textbf{\textbf{\textsl{Test}}} \textbf{5} & \textbf{\textbf{\textsl{Test}}} \textbf{6} & \textbf{\textbf{\textsl{Test}}} \textbf{4} & \textbf{\textbf{\textsl{Test}}} \textbf{5} & \textbf{\textbf{\textsl{Test}}} \textbf{6} \\
 & & %$0.00025$ 
 2.5e-4 & %$0.0025$ 
 2.5e-3 & %$0.025$ 
 2.5e-2 & %$25.0$ 
 2.5e+1 & %$250.0$ 
 2.5e+2 & %$2500.0$
 2.5e+3\\
 & \multirow{2}{*}{$1000.0$} & \textbf{\textbf{\textsl{Test}}} \textbf{7} & \textbf{\textbf{\textsl{Test}}} \textbf{8} & \textbf{\textbf{\textsl{Test}}} \textbf{9} & \textbf{\textbf{\textsl{Test}}} \textbf{7} & \textbf{\textbf{\textsl{Test}}} \textbf{8} & \textbf{\textbf{\textsl{Test}}} \textbf{9}\\
 & & %$0.0025$ 
 2.5e-3 & %$0.025$ 
 2.5e-2 & %$0.25$ 
 2.5e-1 & %$250.0$ 
 2.5e+2 & %$2500.0$ 
 2.5e+3 & %$25000.0$
 2.5e+4\\ \bottomrule
\end{tabular}
\end{center}
\end{table}

\subsubsection{\texttt{GROUP I}: $Re\times Da < 1.0$}\label{sec:Group_1}
%\hl{For} \texttt{GROUP I} of $Re\times Da < 1.0$,  {the}  parameters are as follows; \textbf{\textbf{\textsl{Test}}} \textbf{1} has  $Re\times Da = 0.000025$, where (a) $\dfrac{1}{Re}=0.1$, (b) $\dfrac{1}{ReDa}=40000$, and (c) $\dfrac{c_F}{\sqrt{Da}}=316.23$. For \textbf{\textbf{\textsl{Test}}} \textbf{5}, the corresponding values are $0.0025$, (a) $0.01$, (b) $400$, (c) $100$, respectively, and for \textbf{\textbf{\textsl{Test}}} \textbf{9}, those are $0.25$, (a) $0.001$, (b) $4$, (c) $31.62$, respectively. Note that \textbf{\textbf{\textsl{Test}}} \textbf{1}, \textbf{5}, and \textbf{9} are selected for succinctness.

{As seen in Figure~\ref{Fig:G1_EX159_CASE_1}, %Note that 
\textbf{\textbf{\textsl{Test}}} \textbf{1}, \textbf{5}, and \textbf{9} are selected for succinct demonstration.} 
The (left) column ((a), (d), (g)) is for the linearized Brinkman model, the middle column ((b), (e), (h)) is for {the Darcy-Brinkman model and the (right) column ((c), (f), (i)) is for the full nonlinear Darcy-Brinkman-Forchheimer model.}  The cross mark in red and yellow 
is located at the exact location, the center of each domain. The three models are far from the Navier-Stokes model for the values of Reynolds number.  
This is due to the small $Da$ values resulting in $Re\times Da<1.0$ and higher viscous resistance. {Regarding (b) $\dfrac{1}{ReDa}$ in Equation~\eqref{eq:dimless-full-nonlinear} for \textbf{\textbf{\textsl{Test}}} \textbf{1}, \textbf{5}, and \textbf{9},  their values are $40000$, $400$, and $4$, respectively.} 
 In the case of the full D-B-F model,  higher inertial resistance also works.

The overall shape of streamlines of the three models are very similar, but \textbf{\textsl{Test} 9} exhibits a slight difference 
in  
the location of the center of the eddy between the models, implying that the differentiation between the models has started.  
{Figure}~\ref{Fig:G1_EX159_Velocity_Center} confirms that {a deviation} can be found distinctively with $u_y$ in \textbf{\textsl{Test} 9}. Meanwhile, it is shown that the D-B and the Brinkman models are almost identical in their values for the center line velocities. {It indicates that for $Da$ and $Re$ values of this range, the convection does not make much difference to the momentum balance compared to the inertial resistance from the Forchheimer term.} 
\begin{figure}[tbhp]
\centering
\subfloat[Brinkman \textbf{\textsl{Test} 1}]
{\includegraphics[width=0.3\textwidth]{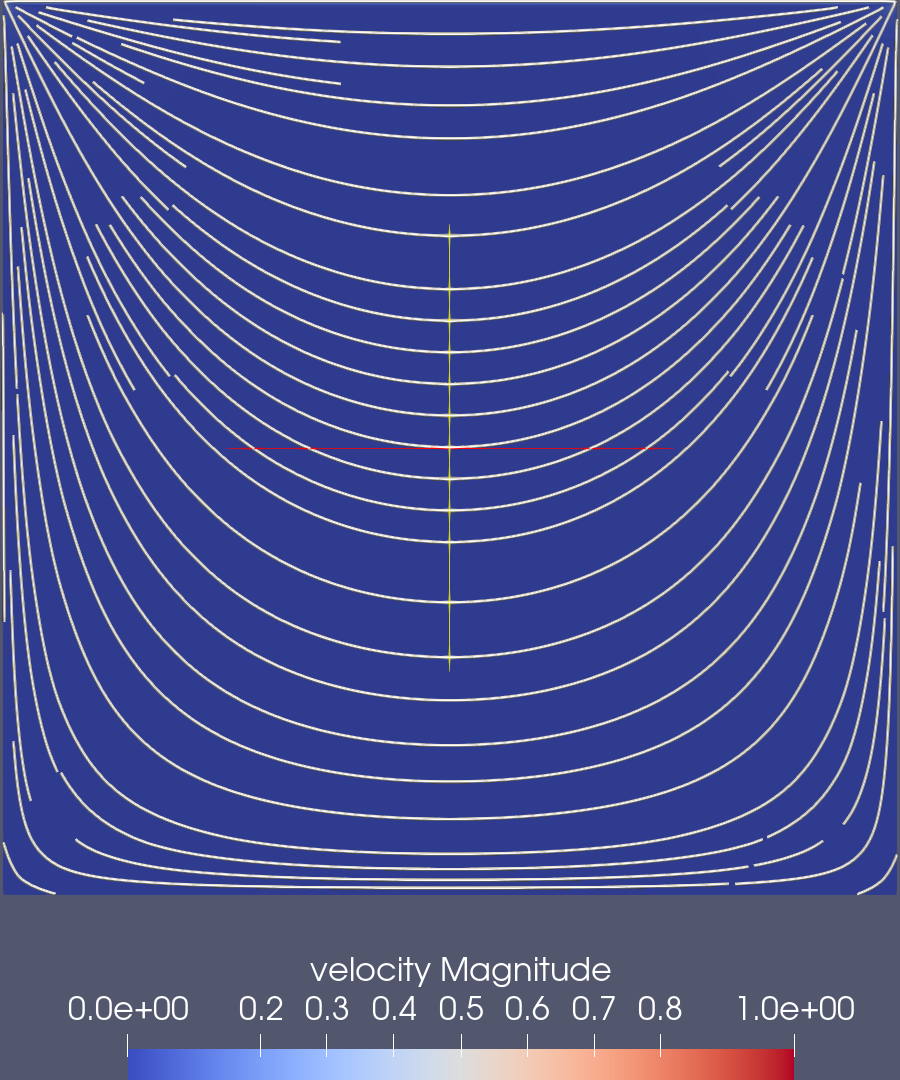}}\hspace*{0.05em}
\subfloat[D-B \textbf{{\textsl{Test} 1}}]
{\includegraphics[width=0.3\textwidth]{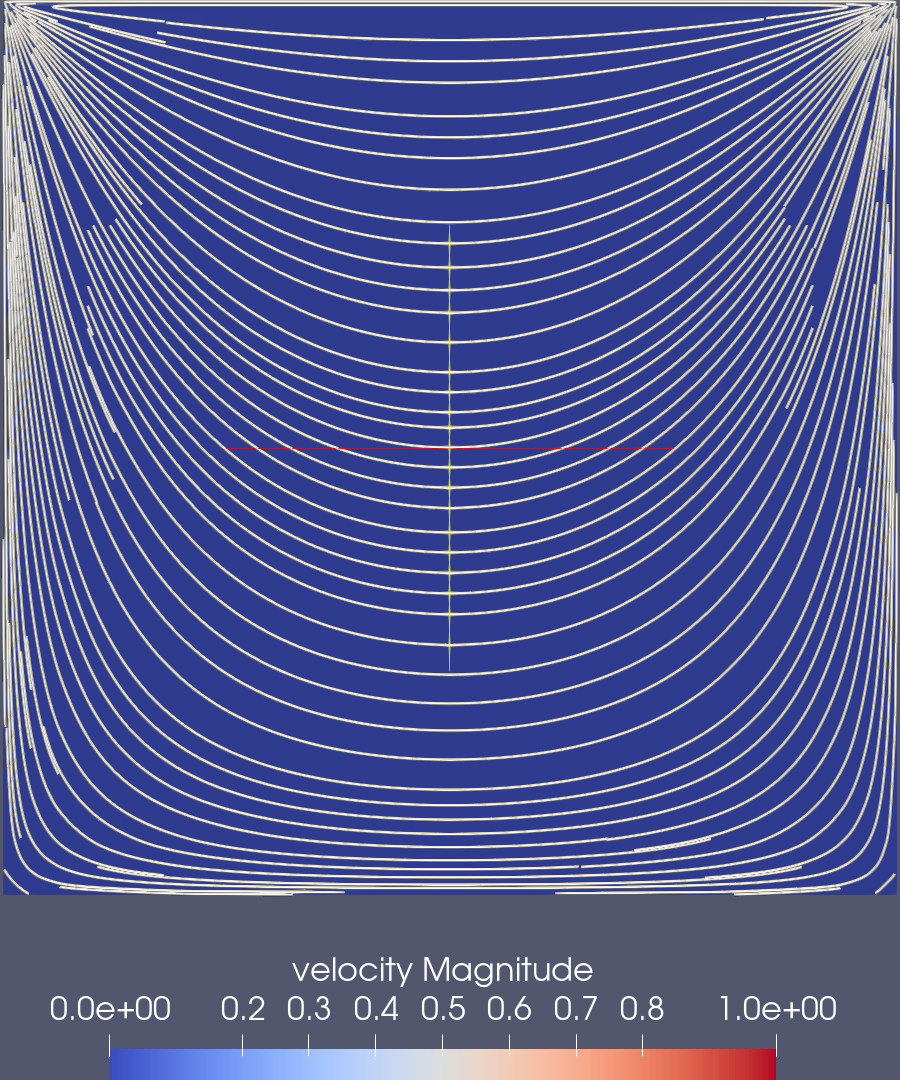}}\hspace*{0.05em}
\subfloat[D-B-F \textbf{\textsl{Test} 1}]
{\includegraphics[width=0.3\textwidth]{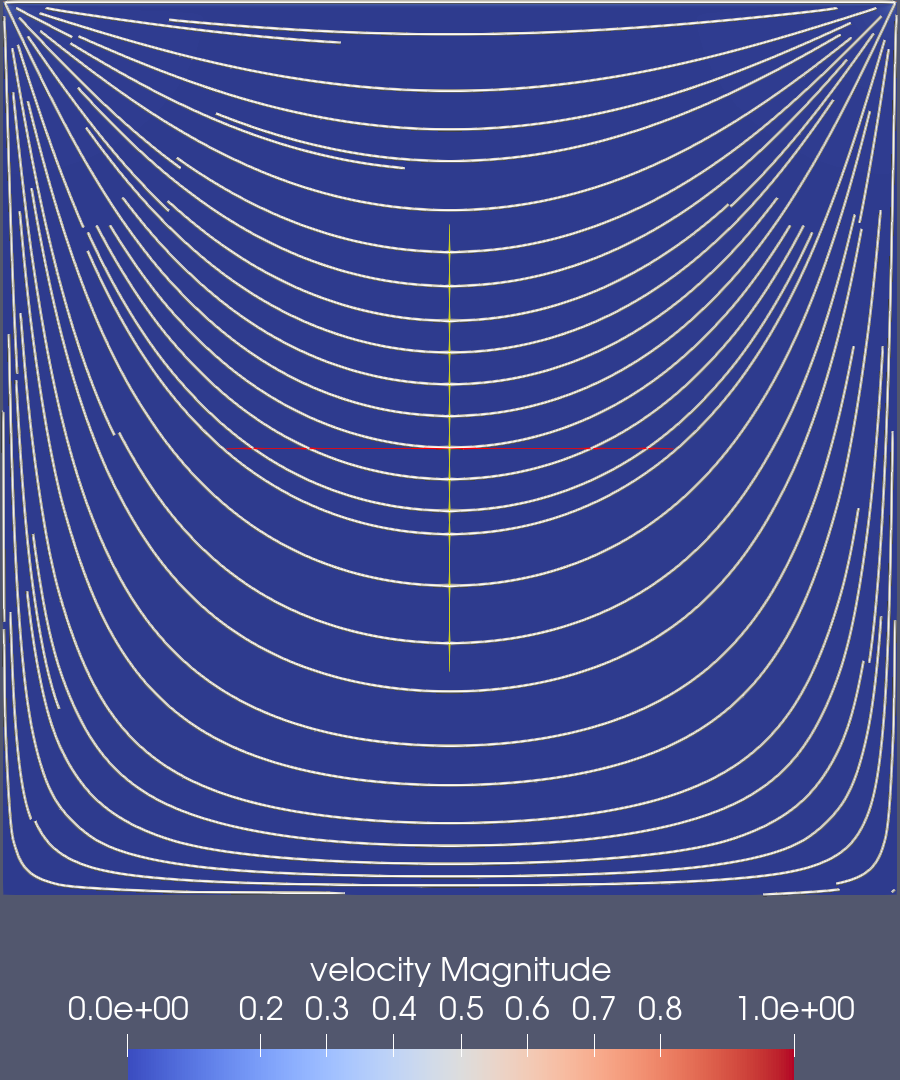}}
\\
\subfloat[Brinkman \textbf{\textsl{Test} 5}]
{\includegraphics[width=0.3\textwidth]{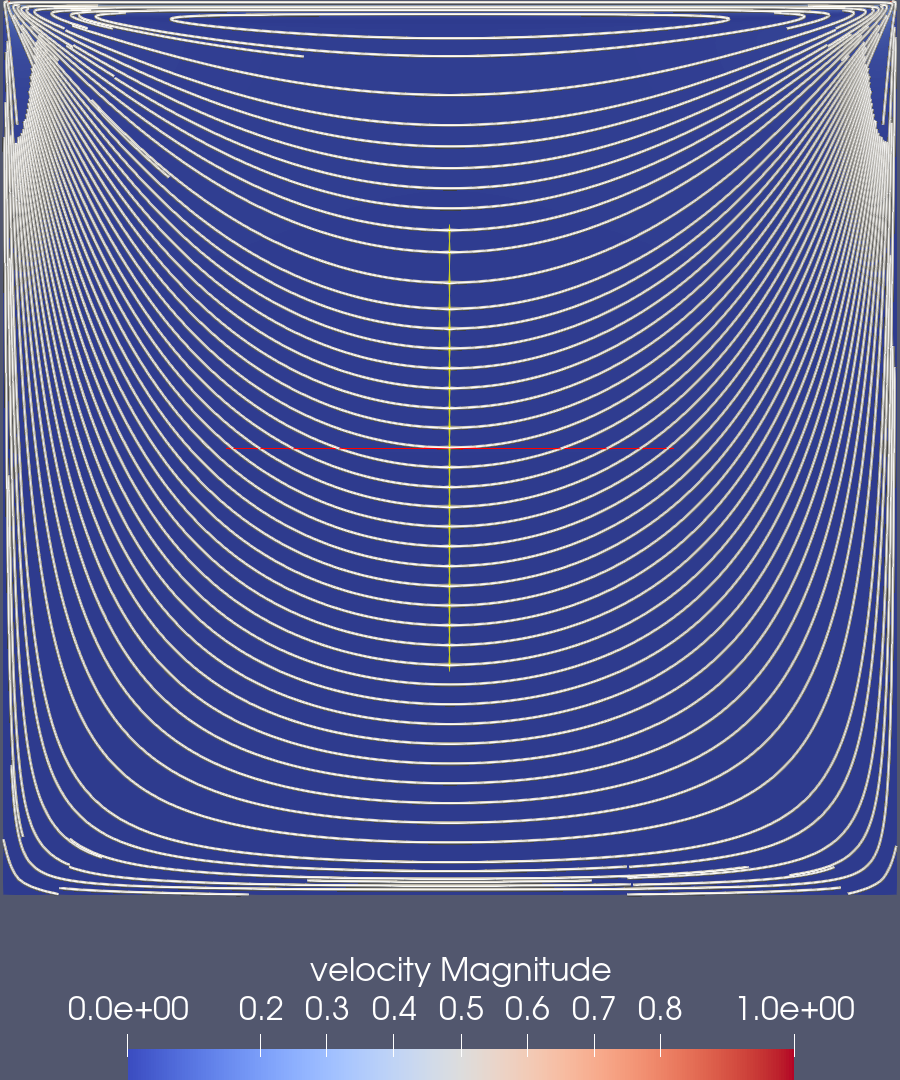}}\hspace*{0.05em}
\subfloat[D-B \textbf{\textsl{Test} 5}]
{\includegraphics[width=0.3\textwidth]{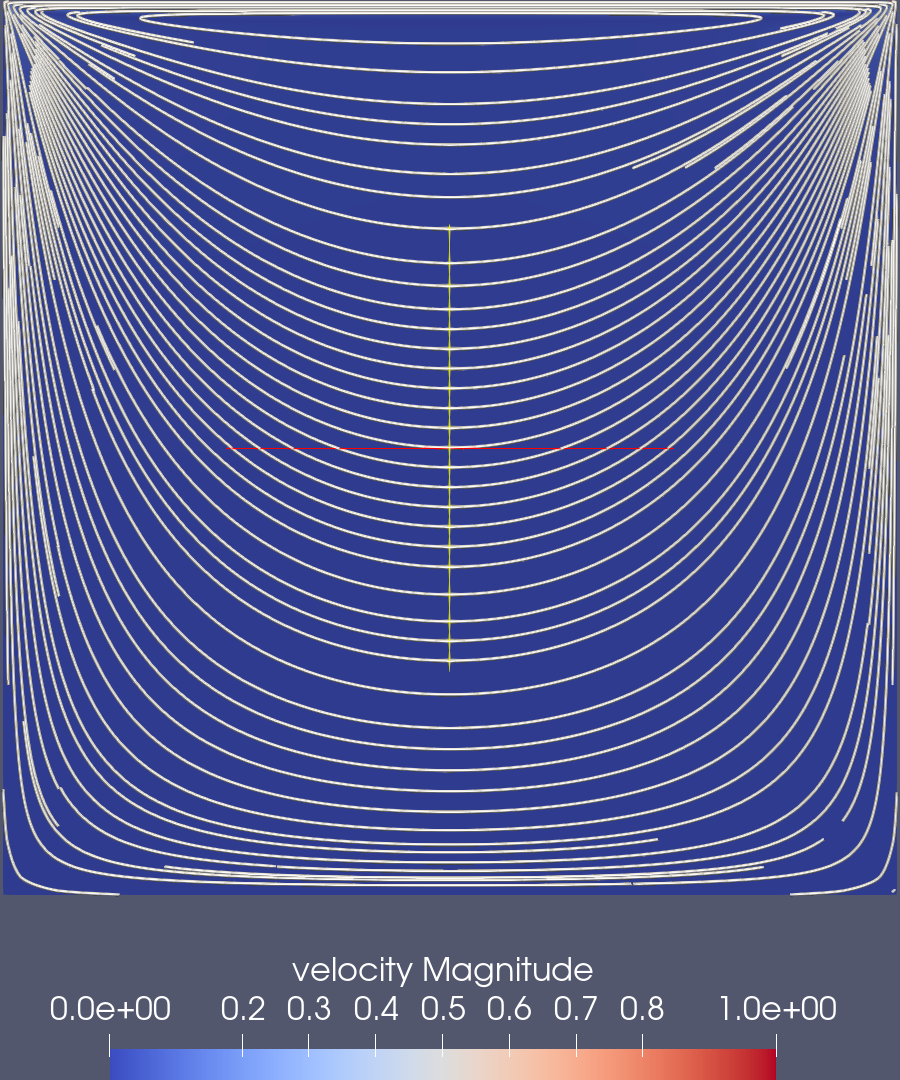}}\hspace*{0.05em}
\subfloat[D-B-F \textbf{\textsl{Test} 5}]
{\includegraphics[width=0.3\textwidth]{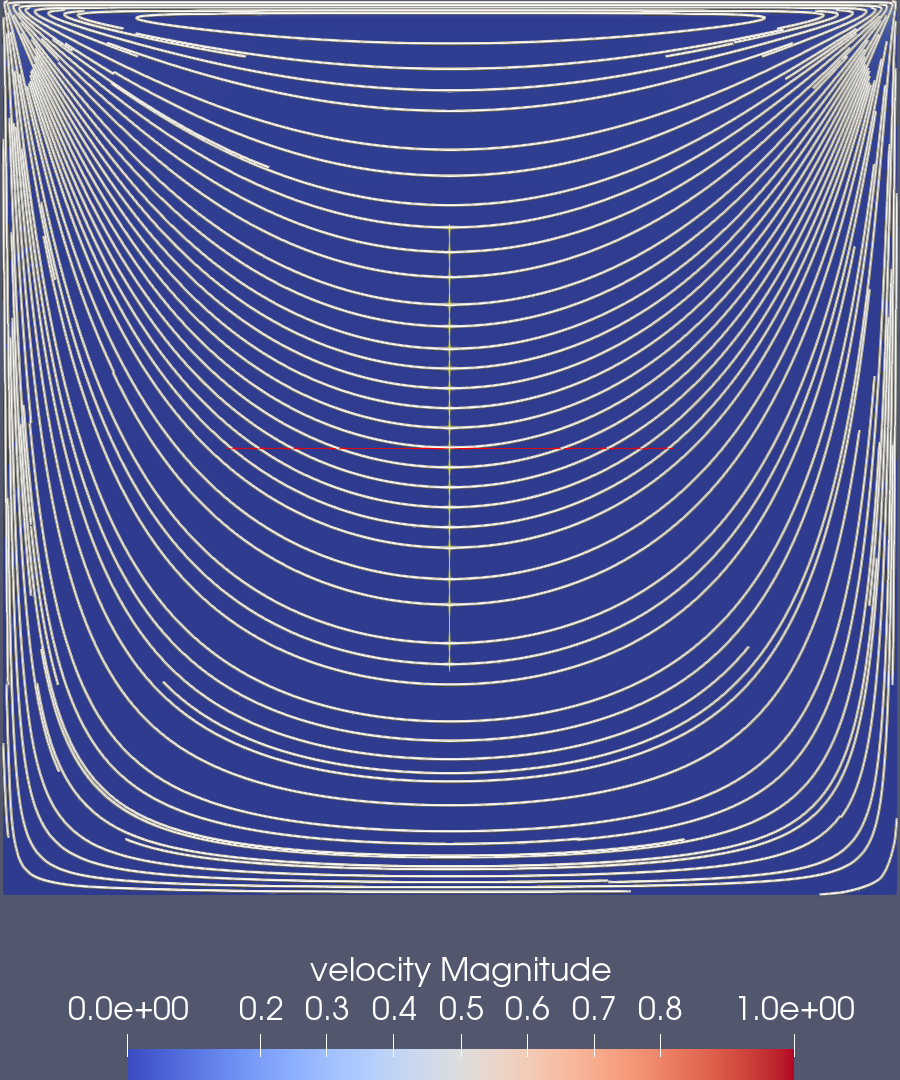}}
\\
\subfloat[Brinkman \textbf{\textsl{Test} 9}]
{\includegraphics[width=0.3\textwidth]{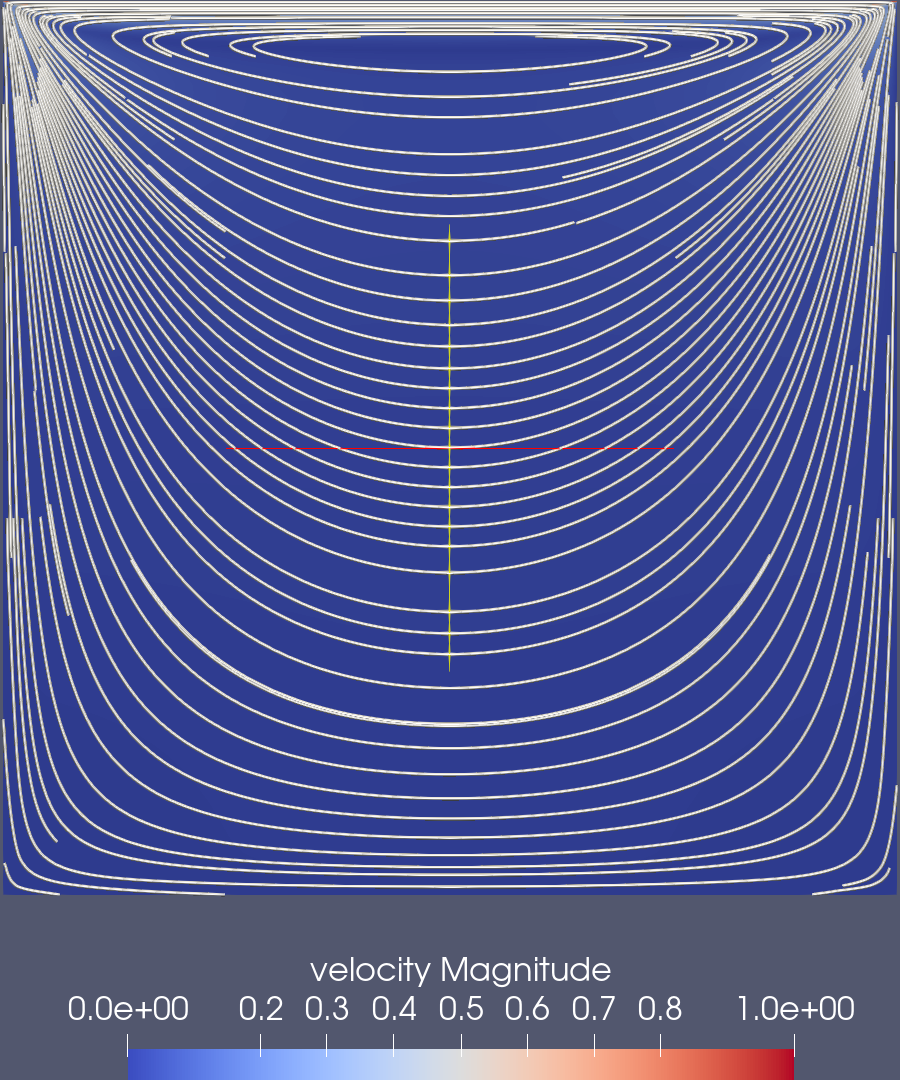}}\hspace*{0.05em}
\subfloat[D-B \textbf{\textsl{Test} 9}]
{\includegraphics[width=0.3\textwidth]{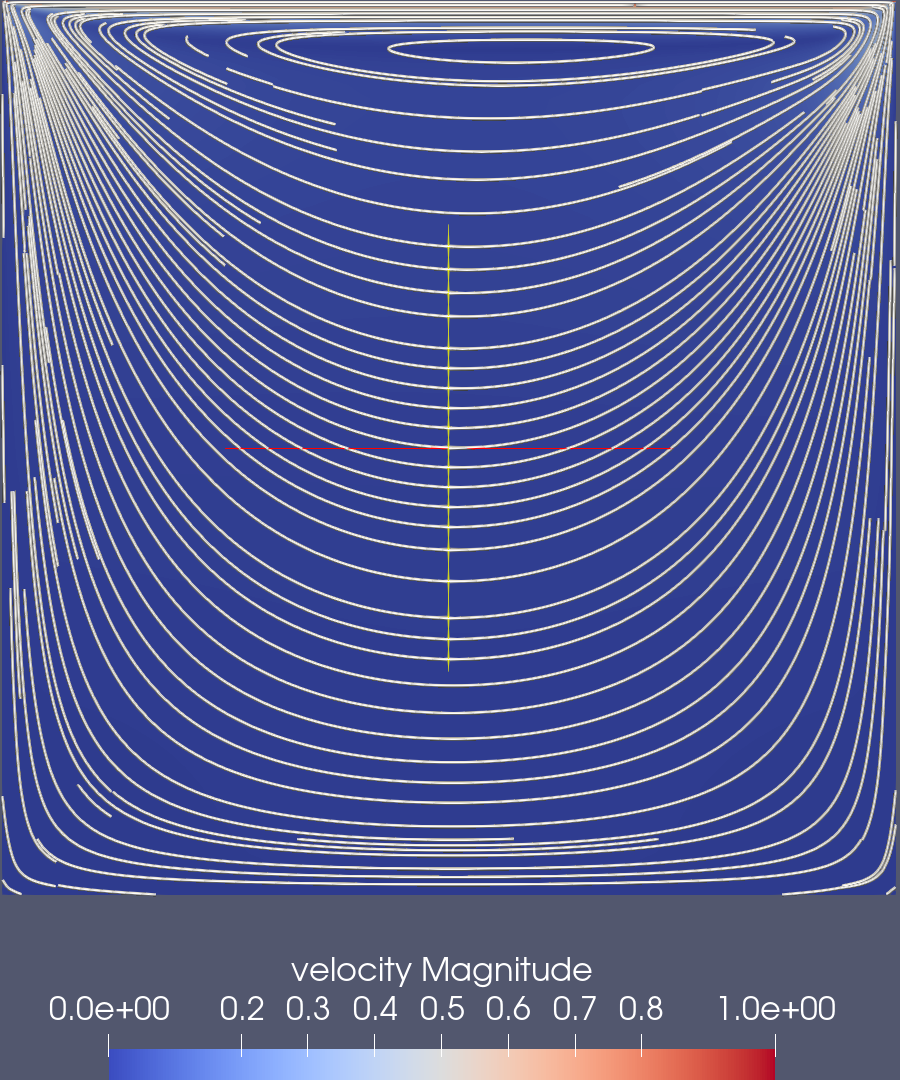}}\hspace*{0.05em}
\subfloat[D-B-F \textbf{\textsl{Test} 9}]
{\includegraphics[width=0.3\textwidth]{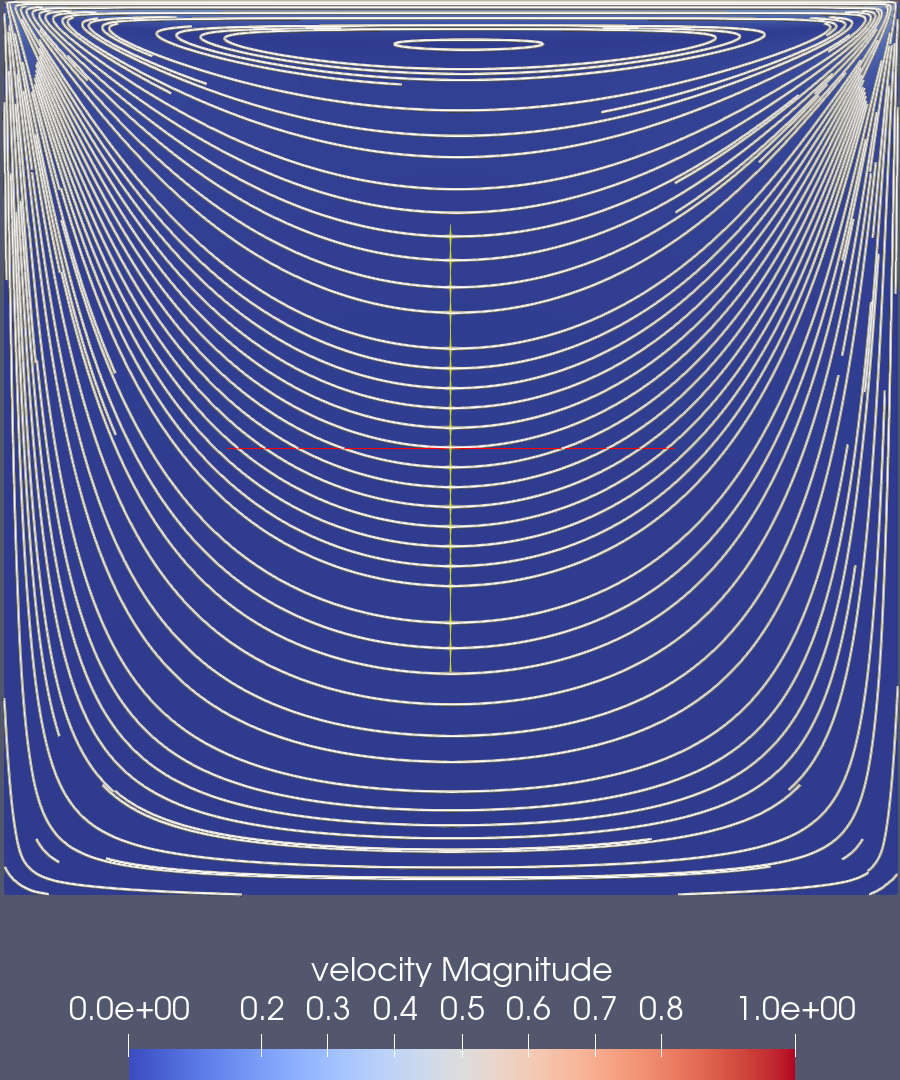}}
\caption{[\texttt{GROUP I}] \textbf{\textsl{Test} 1, 5, and 9}. Velocity streamline: {(left)} Brinkman, {(center)} Darcy-Brinkman, and (right) Darcy-Brinkman-Forchheimer. {\textbf{\textsl{Test} 1, 5, and 9} have dimensionless numbers of 10, 100, 1000 for $Re$, and 2.5e-6, 2,5e-5, and 2.5e-4 for $Da$, respectively.}}
\label{Fig:G1_EX159_CASE_1}
\end{figure}
\begin{figure}[tbhp]
\centering
\includegraphics[width=0.90\textwidth]{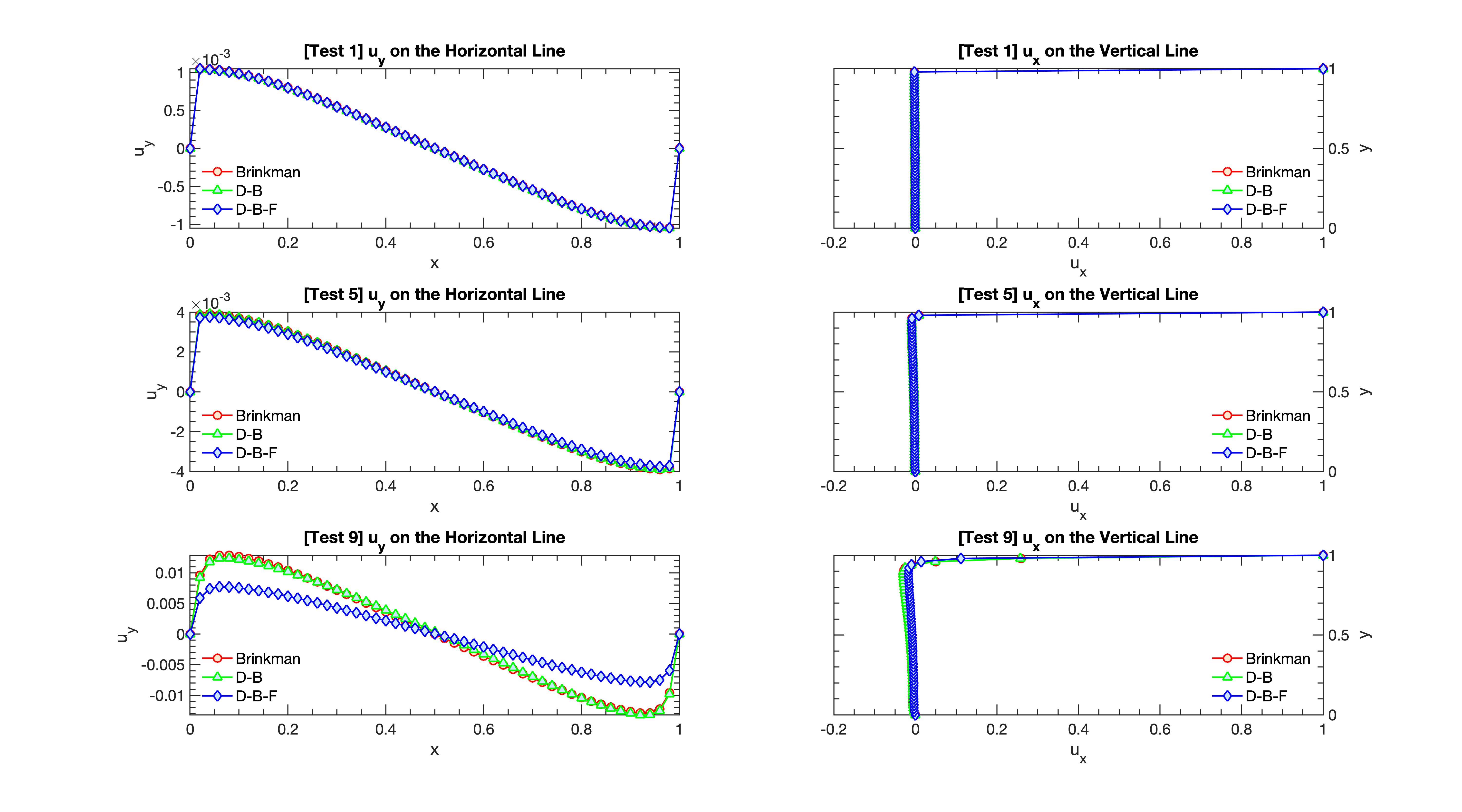}
\caption{[\texttt{GROUP I}] \textbf{\textsl{Test} 1, 5, and 9}. {Velocity} on the reference lines using the Brinkman, Darcy-Brinkman, and Darcy-Brinkman-Forchheimer {models}. {\textbf{\textsl{Test} 1, 5, and 9} have dimensionless numbers of 10, 100, 1000 for $Re$, and 2.5e-6, 2,5e-5, and 2.5e-4 for $Da$, respectively.}}
\label{Fig:G1_EX159_Velocity_Center}
\end{figure}

Table~\ref{tab:Num_Newton_g1_Ex159} presents the {outer} Newton iteration numbers {with flexible preconditioning iterative solver, FGMRES,} in the last adaptive mesh refinement cycle {(i.e., $4^{th}$)} for each model for \textbf{\textsl{Test} 1, 5, and 9}. %three \textbf{\textsl{Tests}}. 
%In this last refinement cycle using the {Kelly error estimator} (Section~\ref{sec:LGR}), 
For both the D-B and D-B-F models, only one more iteration is added  to the iteration number of the Brinkman model. 
{Regarding the stabilization parameter for the Augmented Lagrangian term (\eqref{FLHS-a} and \eqref{FRHS}),} we tested the numerics 
for the D-B-F model with $\gamma$ as $\gamma=0.0$. 
{Although %we omit 
the specifics are omitted, we found that the iteration number for 
FGMRES has been severely increased with inefficiency.} %result is not explicitly shown here}.

On the other hand, {Figure}~\ref{Fig:G1_EX9_MESH} illustrates 
the adaptive refinement using the Kelly error estimator on its last cycle for \textbf{\textsl{Test} 9} with each model. The Brinkman model shows a nearly perfect symmetry due to the absence of the convection term, and the D-B model is the most asymmetric {having the most nodes}. {For the full D-B-F model, it is alleviated by the inertial resistance by the Forchheimer term, which can be also found in Figure~\ref{Fig:G1_EX159_CASE_1} from (g) to (i).} {Thus, it can be inferred that the convection term, still having a small impact, starts to reveal around the corner area compared to the center (see Figure~\ref{Fig:G1_EX159_Velocity_Center})}. We also note that a naive Galerkin finite element method may produce oscillations and numerical instabilities near the corners due to the discontinuous property of the top boundary condition \cite{YooH2018}. 
%In summary for \texttt{GROUP I}, {the flow is found to be controlled mainly by} the viscous resistance, and there are no big differences between the three models: linearized Brinkman, the Darcy-Brinkman, and the Darcy-Brinkman-Forchheimer models. 
\begin{table}[tbhp]
\begin{center}
\caption{[\texttt{GROUP I}] {Total number of Newton iterations} for \textbf{\textbf{\textsl{Test}}} \textbf{1}, \textbf{5}, and \textbf{9}.}
\label{tab:Num_Newton_g1_Ex159}
\begin{tabular}{rccc}
\toprule
\textbf{\textbf{\textsl{Test}}} & Brinkman & D-B & D-B-F \\
\midrule
\textbf{\textbf{\textsl{Test}}} \textbf{1} & $3$ & 4 & {$4$} \\
\textbf{\textbf{\textsl{Test}}} \textbf{5} & $3$ & 4 & $4$ \\
\textbf{\textbf{\textsl{Test}}} \textbf{9} & $3$ & 4 & $4$ \\
\bottomrule
\end{tabular}
\end{center}
\end{table}
\begin{figure}[tbhp]%htbp
\centering
\subfloat[Brinkman \textbf{\textsl{Test} 9}]
{\includegraphics[width=0.3\textwidth]{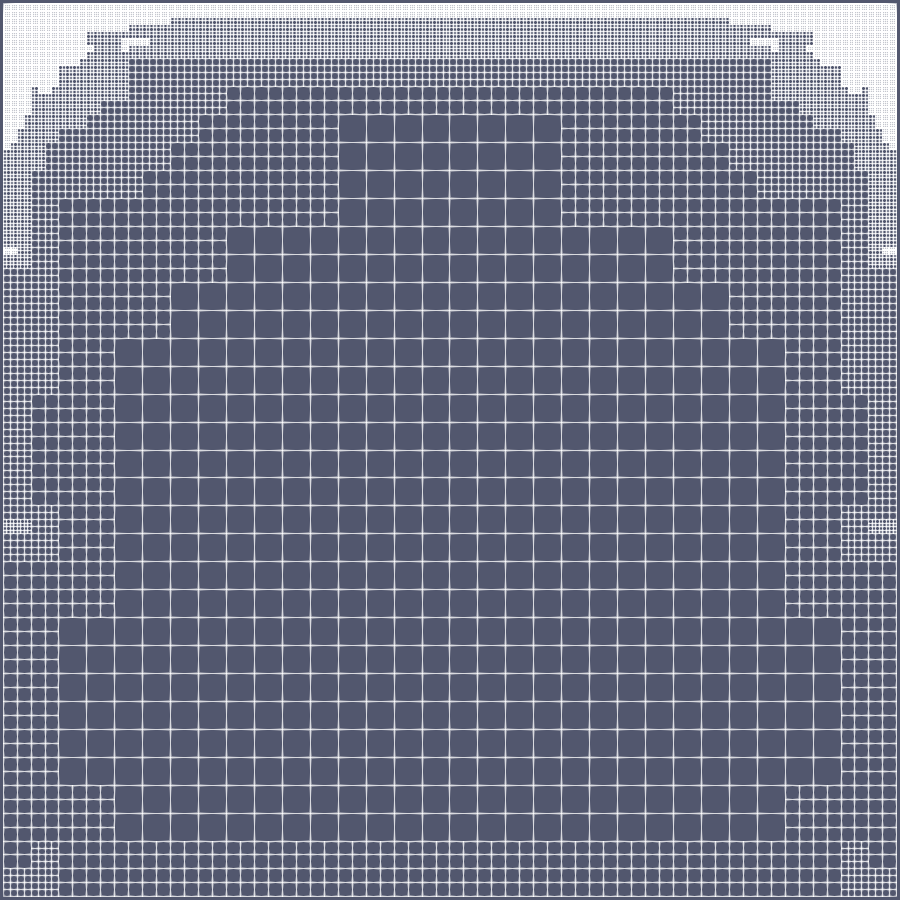}}\hspace*{0.05em}
\subfloat[D-B \textbf{\textsl{Test} 9}]
{\includegraphics[width=0.3\textwidth]{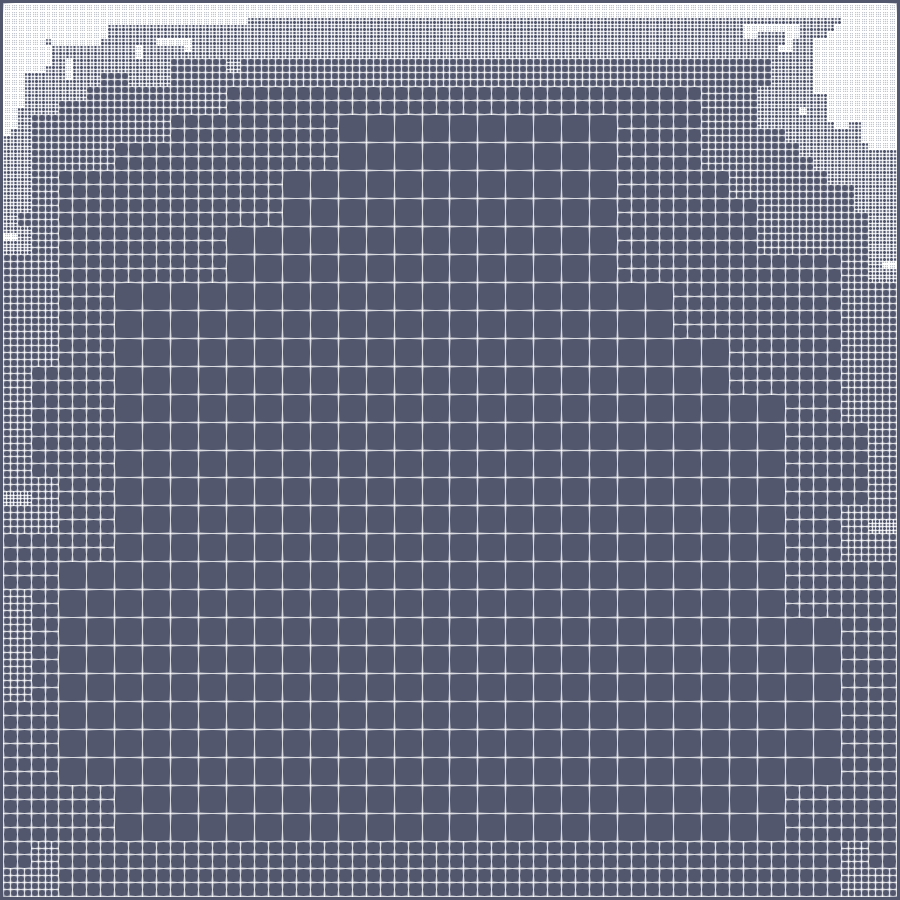}}\hspace*{0.05em}
\subfloat[D-B-F \textbf{\textsl{Test} 9}]
{\includegraphics[width=0.3\textwidth]{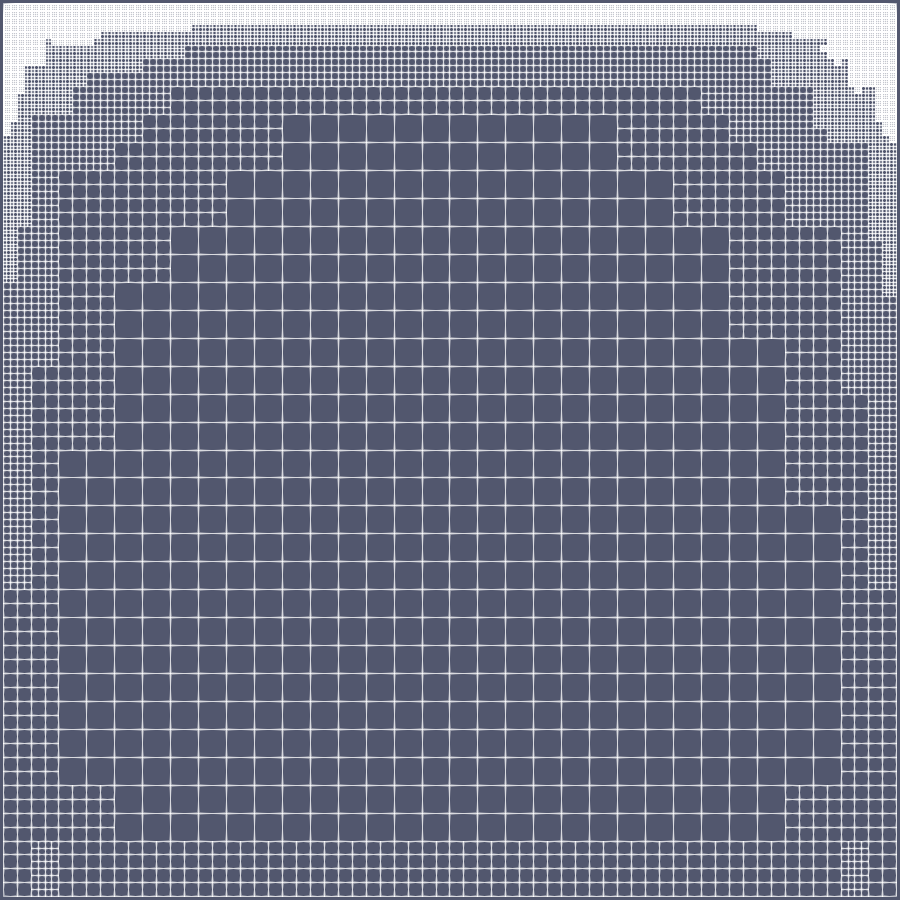}}
\caption{[\texttt{GROUP I}] \textbf{\textsl{Test} 9}. {Adaptive mesh refinements} of the last ($4^{th}$) cycle using the Kelly error estimator: (a) Brinkman, (b) Darcy-Brinkman, and (c) Darcy-Brinkman-Forchheimer. {The total numbers of DoFs are (a) 154099 (DoFs for $\mathbf{u}_{h}$:136762, for $p_h$:17337), (b) 158253 (DoFs for $\mathbf{u}_{h}$:140454, for $p_h$:17799), (c) 144475 (DoFs for $\mathbf{u}_{h}$:128210, for $p_h$:16265)}.}
\label{Fig:G1_EX9_MESH}
\end{figure}

\subsubsection*{{{In-depth study} for \texttt{GROUP I} (with $Re$=10)}}
{In this subsection, we investigate further about \texttt{GROUP I} with $Re$=10, and explore two other $Da$ values within $Re$=10.} 
\subparagraph*{{$\bullet$ \textbf{\textsl{Test 2}} \& \textbf{\textsl{Test 3}}:}} {For in-depth study, we choose  \textbf{\textsl{Test 2}} and \textbf{\textsl{Test 3}} and examine the Newton iterations of each model via probing numbers of FGMRES steps in each iteration. Both \textbf{\textsl{Test 2}} \& \textbf{\textsl{Test 3}} have 3, 4, and 4 iterations in Newton for the Brinkman, D-B, and D-B-F models, respectively, the same as previous \textbf{\textsl{Tests}} in Table~\ref{tab:Num_Newton_g1_Ex159} at the last adaptive mesh refinement cycle, (i.e., $4^{th}$). However in \textbf{\textsl{Test 2}}, reaching the stopping criterion of $\varepsilon=1.0e$-$12$ takes  the total FGRMES steps of 371, 449, and 435 for the Brinkman, D-B, and D-B-F models, respectively. For \textbf{\textsl{Test 3}}, corresponding numbers are 49, 77, 84. Overall, the Brinkman has quite less steps in FGMRES than others, while the D-B, and D-B-F models are close. Note also that \textbf{\textsl{Test 3}} has much less numbers in those steps than those of \textbf{\textsl{Test 2}}, resulting from  the Equation~\eqref{lforms2} as $Da$ number increases while $Re$ fixed as $Re=10$.} 
%EX2\\
%Brink at Refinement 4: 2, 130, 239 = 371\\
%DB at Refinement 4:  2, 118, 187, 142 = 449\\
%DBF at Refinement 4: 2, 118, 182, 133 = 435\\
%\\
%EX3\\
%Brink at Refinement 4: 2, 20, 27 = 49\\
%DB at Refinement 4:  2, 19, 28, 28 = 77\\
%DBF at Refinement 4: 2, 19, 28, 35 = 84\\
\begin{table}[tbhp]
\begin{center}
\caption{{[\texttt{GROUP I}] Number of FGMRES steps with Newton iterations for \textbf{\textbf{\textsl{Test}}} \textbf{2} and \textbf{3}}.}
\label{tab:Num_Newton_g1_Ex23}
\begin{tabular}{ccccc}
\toprule
\textbf{\textbf{\textsl{Test}}} & \textbf{\textbf{\textsl{Newton Iter.}}} & Brinkman & D-B & D-B-F \\
\midrule
\multirow{4}{*}{\textbf{\textsl{Test 2}}} &  \textbf{1} & 2 & 2 & 2 \\
& \textbf{2} & 130 & 118 & 118 \\
& \textbf{3} & 239 & 187 & 182 \\
& \textbf{4} & - & 142 & 133 \\
\midrule
\multirow{4}{*}{\textbf{\textsl{Test 3}}} &  \textbf{1} & 2 & 2 & 2 \\
& \textbf{2} & 20 & 19 & 19 \\
& \textbf{3} & 27 & 28 & 28 \\
& \textbf{4} & - & 28 & 35 \\
\bottomrule
\end{tabular}
\end{center}
\end{table}
\subparagraph*{{$\bullet$ $Da=0.0025$ and $Da=0.025$:}}
{As %we have 
seen in Figure~\ref{Fig:G1_EX159_Velocity_Center}, a clear deviation of D-B-F model from two other models is found for the turbulent regime, i.e., $Re=1000$ (\textbf{\textbf{\textsl{Test}}} \textbf{9}). We further test if any notable deviation starts while fixing $Re=10$. Since \textbf{\textbf{\textsl{Test}}} \textbf{2} and \textbf{3} does not show those deviations (not reported explicitly here), we test cases outside \texttt{GROUP I} (in between  \texttt{GROUP I} and  \texttt{GROUP II}) via increasing $Da$ values to $Da=0.0025$ and $Da=0.025$, while fixing $Re=10$.} 
\begin{figure}[tbhp]
\centering
\includegraphics[width=0.90\textwidth]{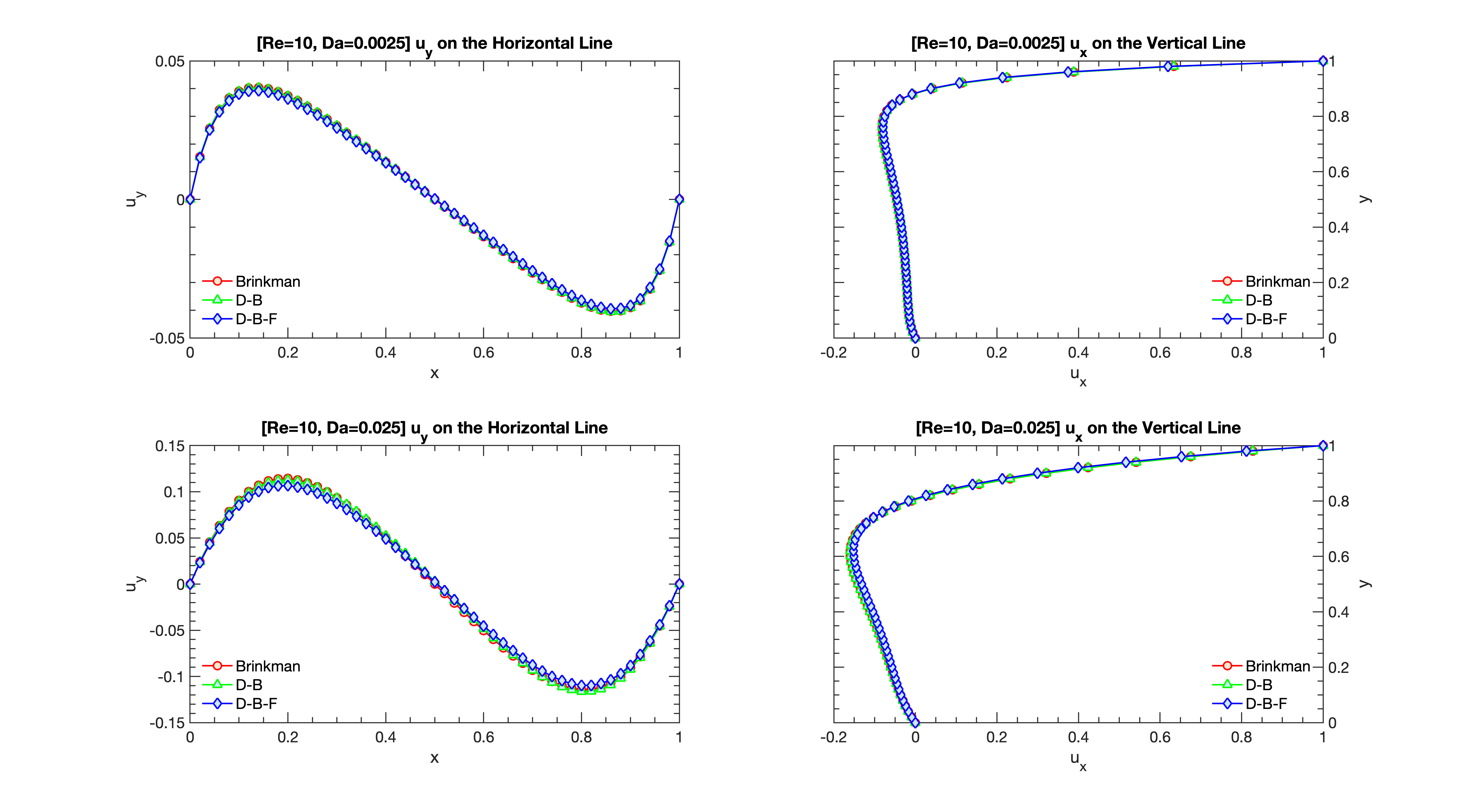}
\caption{[$Re=10$] \textbf{$Da=0.0025$ (first row) and $Da=0.025$ (second row)}. {Velocity} on the reference lines using the Brinkman, Darcy-Brinkman, and Darcy-Brinkman-Forchheimer {models}.}
\label{Fig:G1_EX3p57_Velocity_Center}
\end{figure}

{As demonstrated in Figure~\ref{Fig:G1_EX3p57_Velocity_Center}, we find a slight departure has started when $Da=0.0025$ and more in case of $Da=0.025$. Although $Da$ is not realistic nor any experimental, it seems a very tiny deviation, called \textit{weak inertia} regime \cite{ZimR2004} exists for $Re=10$, which can be negligible from engineering perspective. Note also that \textbf{\textbf{\textsl{Test}}} \textbf{9} with $Re$ variation captures more deviations of D-B-F in $u_y$ than those in $u_x$. Unlike the case, $Da=0.025$ has similar deviations of D-B-F in $u_x$ and $u_y$, even though both \textbf{\textbf{\textsl{Test}}} \textbf{9} and $Da=0.025$ have the same $Re \times Da=0.25$. This is attributed to the property of $Re$, where the inverse of it implies the ratio of viscous force to inertial force (see Equation~\eqref{eq:Re}).  We will see similar trend can be found in the \texttt{GROUP II} \textbf{\textbf{\textsl{Test}}} \textbf{1}.
{Remind} that the D-B-F model is
{convergent} to the N-S system {when} $Da$ goes to greater number, {ultimately being identical if} $Da\rightarrow\infty$.} 

% ----------------------------- % ----------------------------------------------
% -------------------------

\subsubsection{\texttt{GROUP II}: $Re\times Da > 1.0$}
{In this section}, several tests are addressed for \texttt{GROUP II} having $Re\times Da > 1.0$ and the Navier-Stokes model can also be compared with the three models if necessary. 
For brevity as well, \textbf{\textbf{\textsl{Test}}} \textbf{1}, \textbf{5}, and \textbf{9} are selected, 
Some other examples are left for in-depth study in the following subsection. 
%\textbf{\textbf{\textsl{Test}}} \textbf{1} {has}  $Re\times Da = 2.5$, (a) $\dfrac{1}{Re}=0.1$, (b) $\dfrac{1}{ReDa}=0.4$, and (c) $\dfrac{c_F}{\sqrt{Da}}=1.0$. For \textbf{\textbf{\textsl{Test}}} \textbf{5}, the corresponding values are $250$, (a) $0.01$, (b) $0.004$, (c) $0.316$, respectively, and for \textbf{\textbf{\textsl{Test}}} \textbf{9}, those are $25000$, (a) $0.001$, (b) $0.00004$, (c) $0.1$, respectively.

{Unlike} \texttt{GROUP I}, {the secondary and tertiary vortices at the left and right bottom corner appear in \texttt{GROUP II}}, showing the enhanced nonlinearity from increased %ing values as anticipated
of $Re\times Da$. As seen in Figure~\ref{Fig:G2_EX159}, %it is also found that 
both the D-B (center) and the full nonlinear D-B-F (right) models have departed from the Brinkman (left) model in overall shape of streamlines, except for \textbf{\textsl{Test} 1} ((a), (b), (c) in Figure~\ref{Fig:G2_EX159}). %which 
{It results from the relatively high viscous force working  %with small values of $Da$ 
for \textbf{\textsl{Test} 1}, where similar phenomena can apply to other $Re=10$ cases, i.e., \textbf{\textsl{Test} 2} and \textbf{\textsl{Test} 3} as aforementioned in the previous section.}  
%More importantly, 
{As the flow regime changes from the laminar to the turbulent (from \textbf{\textsl{Test} 1} to \textbf{\textsl{Test} 9})}, we find that the center of primary eddy of D-B or D-B-F model  
migrates down toward the center cross due to the  
increased intensity of convection term, as {\texttt{GROUP II} is under weaker viscous resistance} unlike \texttt{GROUP I}. %which is different from \texttt{GROUP I} in that {\texttt{GROUP II} is under weaker viscous resistance}. 
Meanwhile, without the convection term, those of the Brinkman stay at the same height on top of the vertical cross line at the center.
%Compared to the N-S, 
The overall shape of the streamlines of the D-B and D-B-F models appear to be similar except for slight differences in the type and size of curves for the secondary and tertiary vortices at the bottom left and right corners. 
{The size of vortices for D-B-F model is modestly smaller than those of D-B} %and N-S models 
due to the inertial resistance, i.e., (c) $\dfrac{c_F}{\sqrt{Da}}$ in \eqref{eq:dimless-full-nonlinear}. 
\begin{figure}[htbp]
\centering
\subfloat[Brinkman \textbf{\textsl{Test} 1}]
{\includegraphics[width=0.3\textwidth]{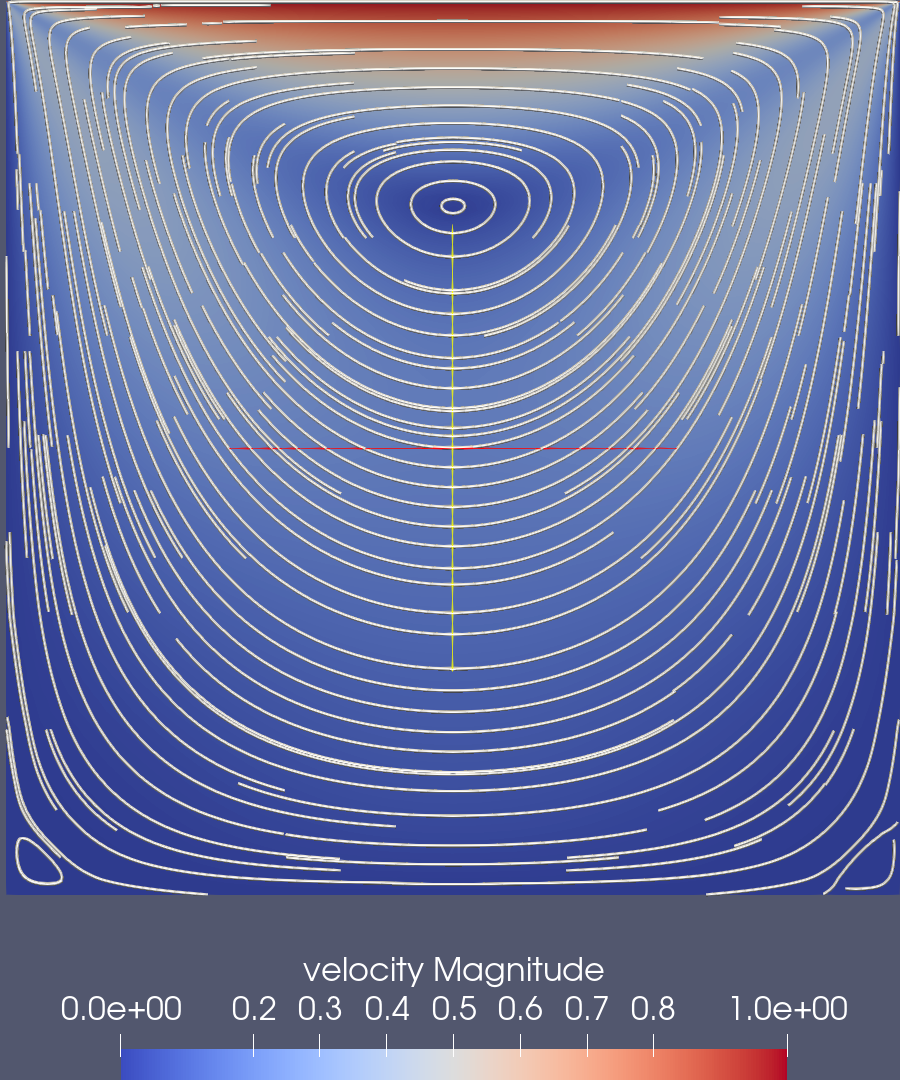}}\hspace*{0.05em}
\subfloat[D-B \textbf{\textsl{Test} 1}]
{\includegraphics[width=0.3\textwidth]{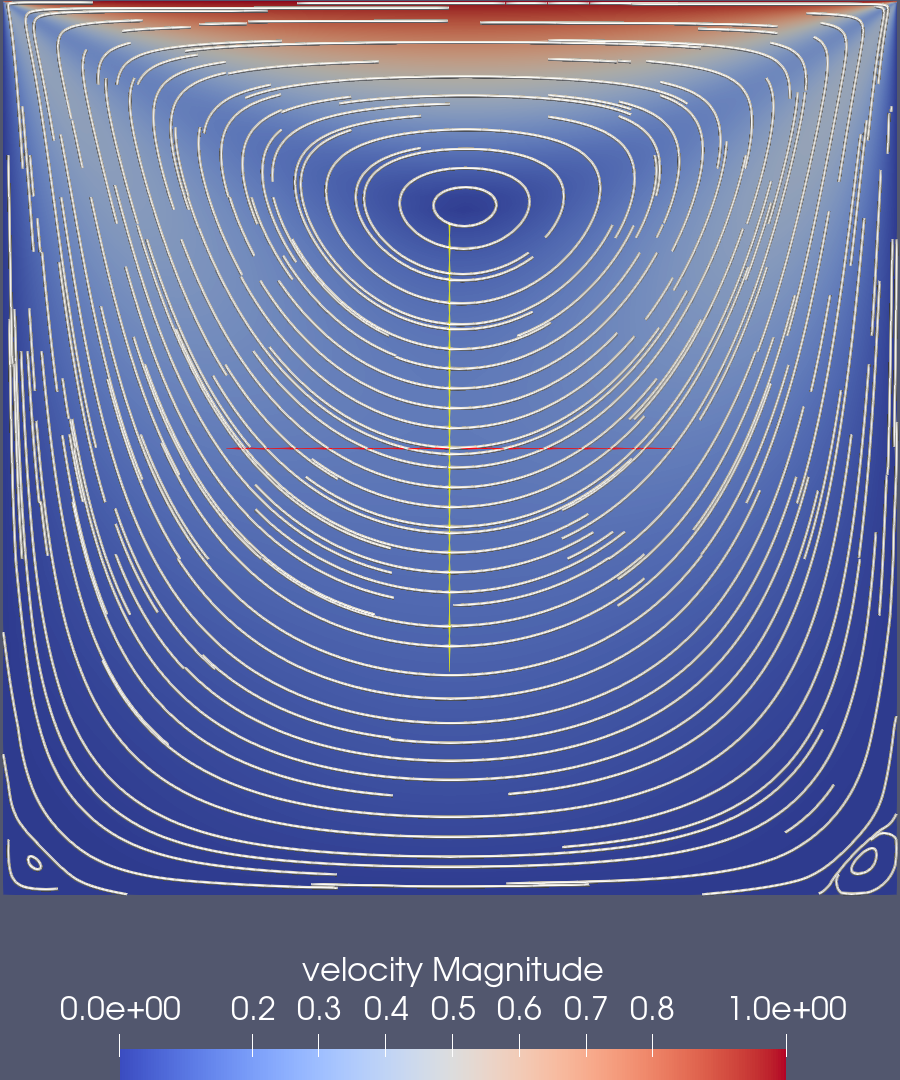}}\hspace*{0.05em}
\subfloat[D-B-F \textbf{\textsl{Test} 1}]
{\includegraphics[width=0.3\textwidth]{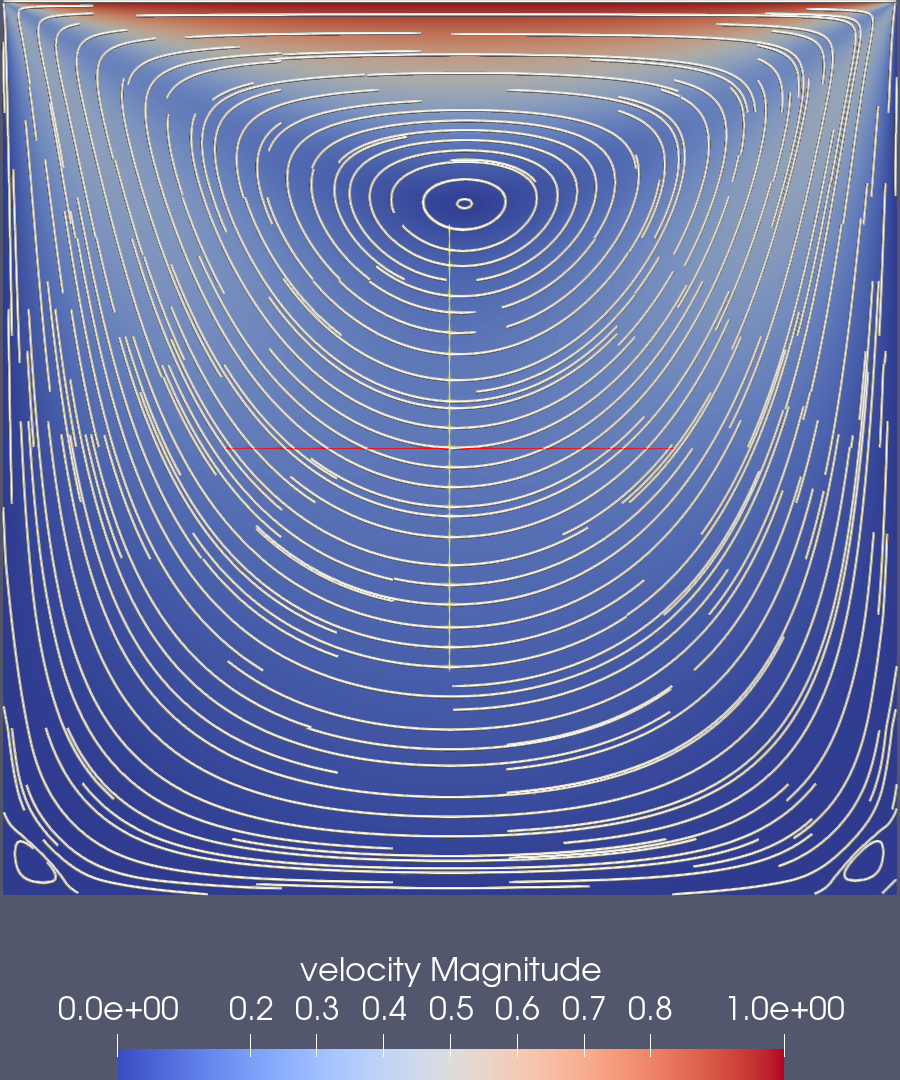}}
\\
\subfloat[Brinkman \textbf{\textsl{Test} 5}]
{\includegraphics[width=0.3\textwidth]{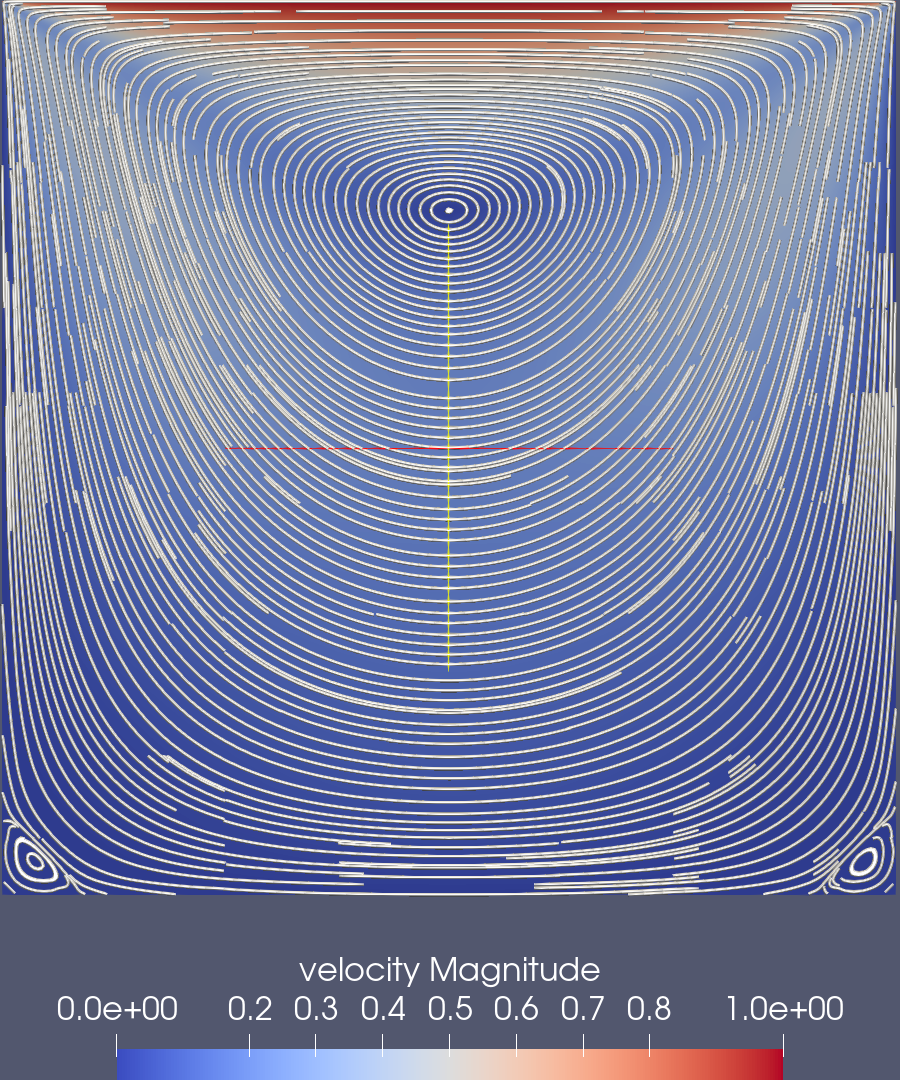}}\hspace*{0.3em}\subfloat[D-B {\textbf{\textsl{Test} 5}}]
{\includegraphics[width=0.3\textwidth]{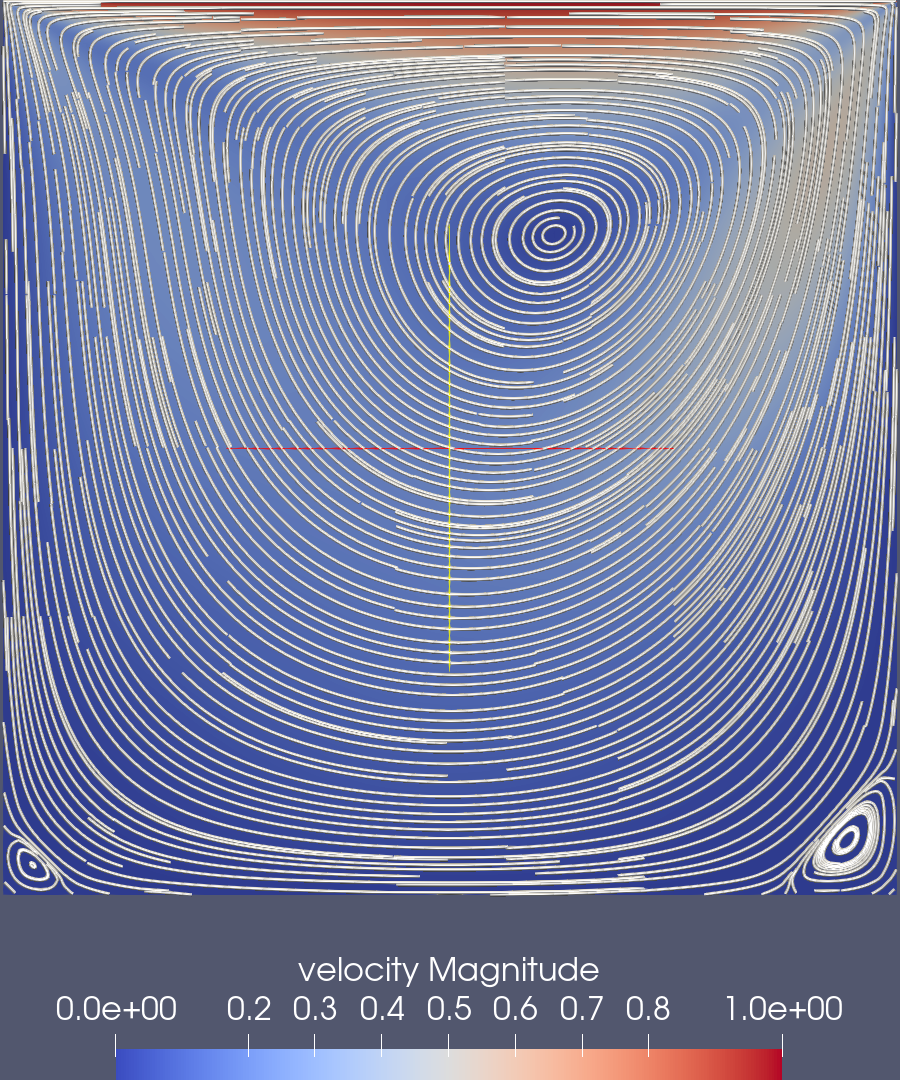}}\hspace*{0.05em}
\subfloat[D-B-F \textbf{\textsl{Test} 5}]
{\includegraphics[width=0.3\textwidth]{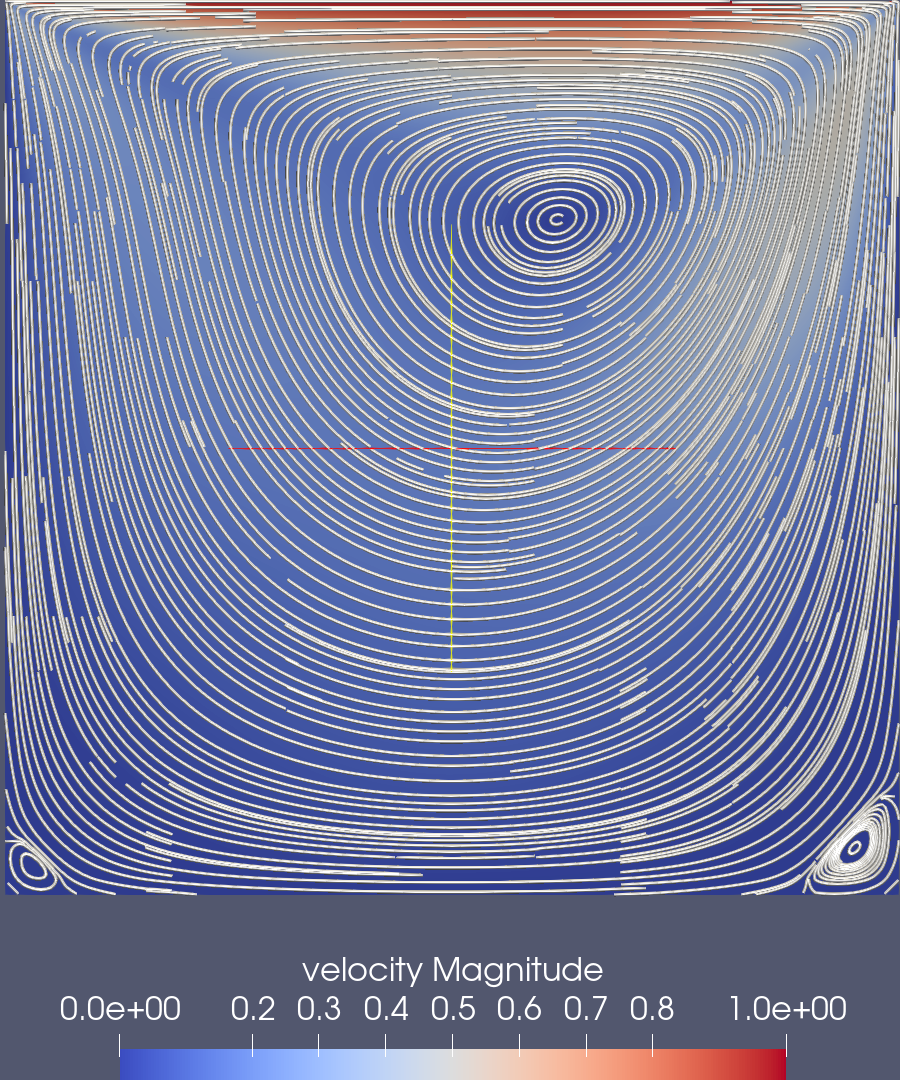}}
\\
\subfloat[Brinkman \textbf{\textsl{Test} 9}]
{\includegraphics[width=0.3\textwidth]{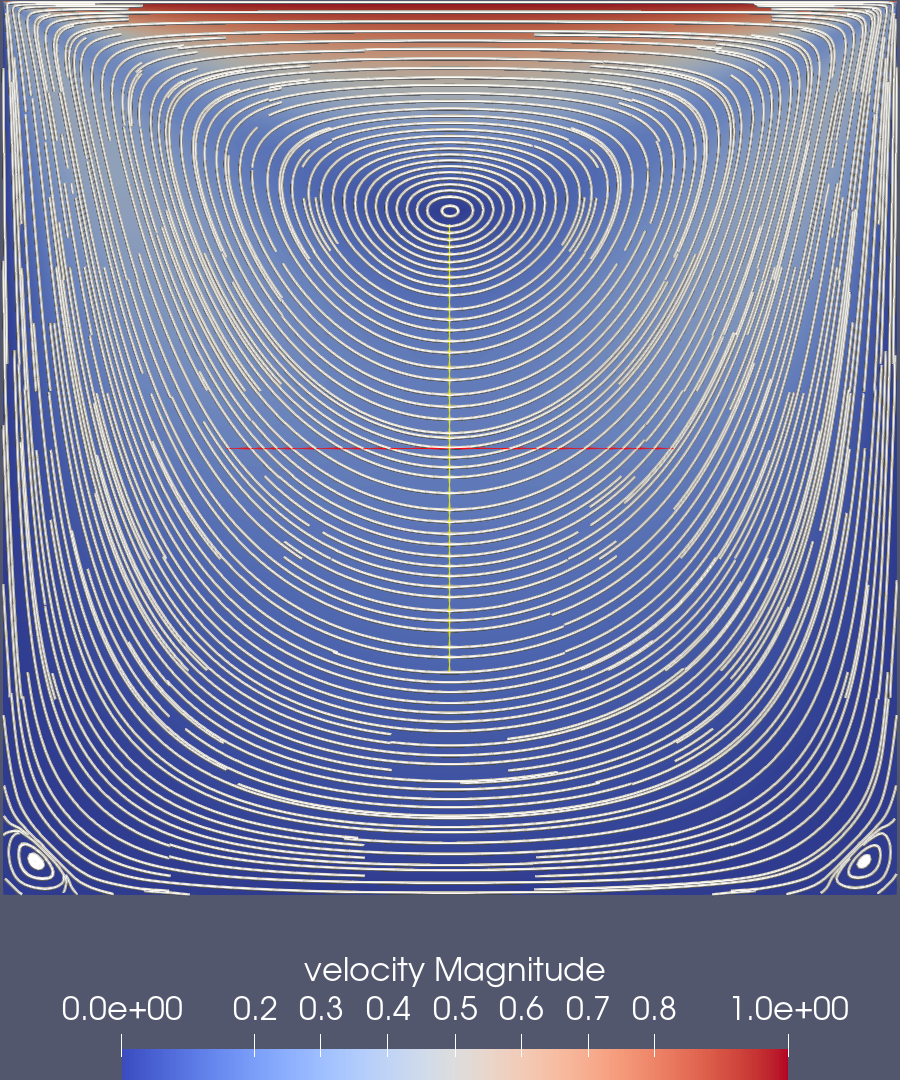}}\hspace*{0.05em}
\subfloat[D-B \textbf{\textsl{Test} 9}]
{\includegraphics[width=0.3\textwidth]{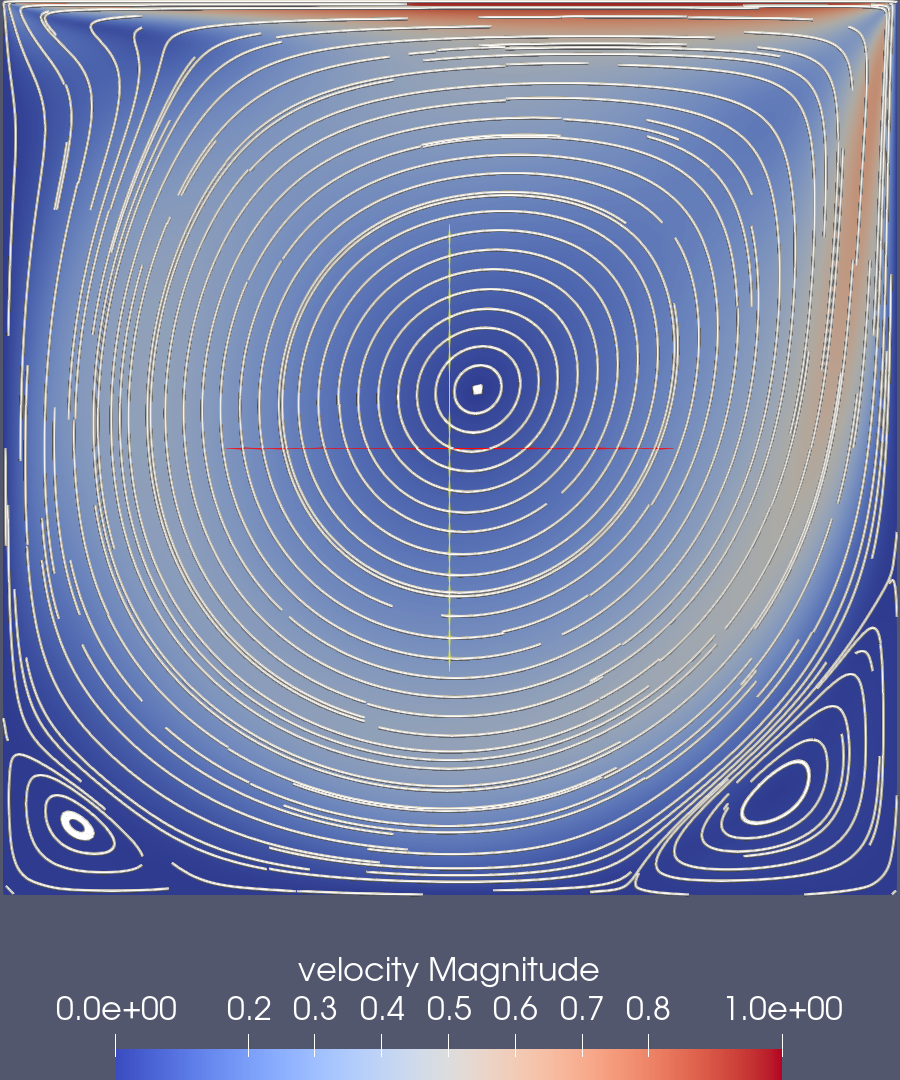}}\hspace*{0.05em}
\subfloat[D-B-F \textbf{\textsl{Test} 9}]
{\includegraphics[width=0.3\textwidth]{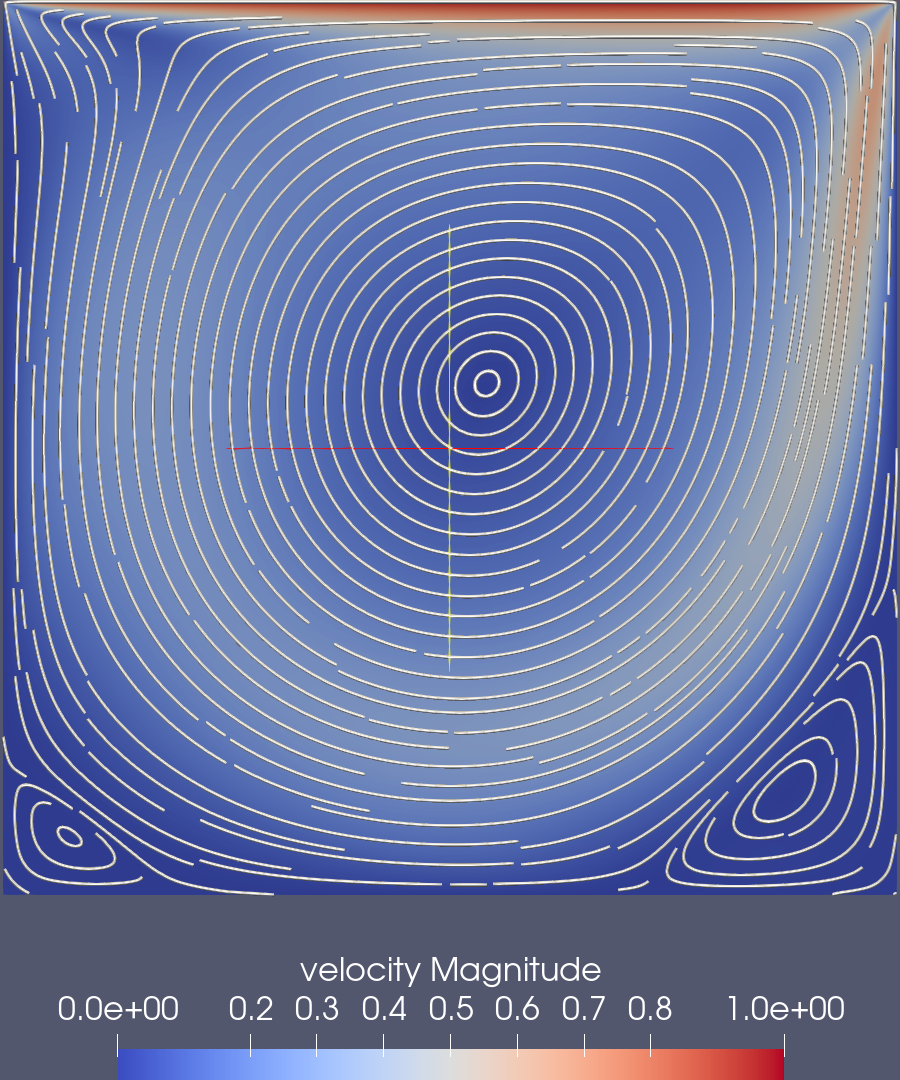}}
\caption{[\texttt{GROUP II}] \textbf{\textsl{Test} 1, 5, and 9}. Velocity streamline: {(left)} Brinkman, {(center)} Darcy-Brinkman, and (right) Darcy-Brinkman-Forchheimer. {\textbf{\textsl{Test} 1, 5, and 9} have dimensionless numbers of 10, 100, 1000 for $Re$, and 2.5e-1, 2,5e+0, and 2.5e+1 for $Da$, respectively.}}
\label{Fig:G2_EX159}
\end{figure}
\begin{figure}[htbp]
\centering
\includegraphics[width=0.9\textwidth]{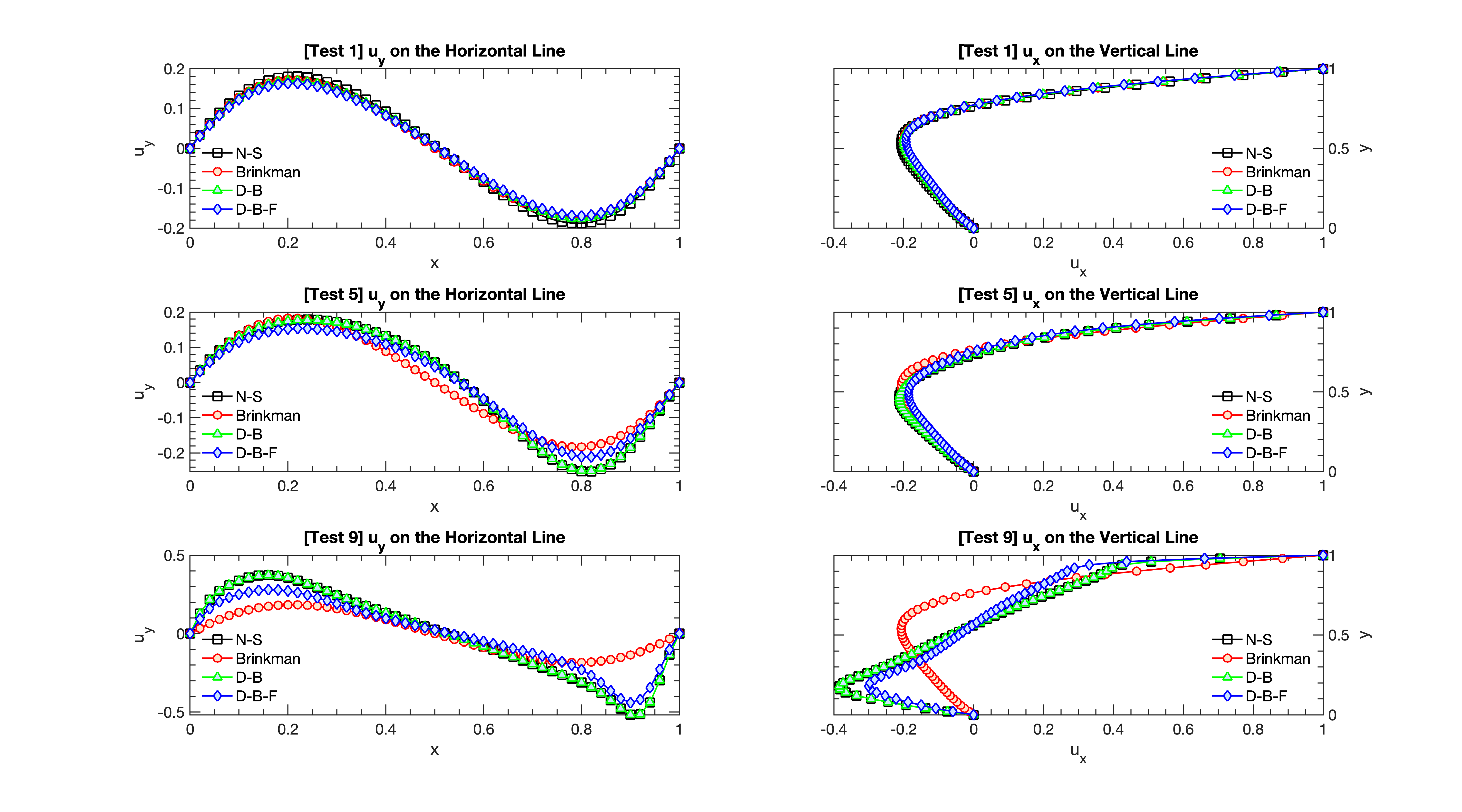}
\caption{[\texttt{GROUP II}] \textbf{\textsl{Test} 1, 5, and 9}. Velocity on the reference lines using the Brinkman, Darcy-Brinkman, and Darcy-Brinkman-Forchheimer {models}, compared with the Navier-Stokes model. {\textbf{\textsl{Test} 1, 5, and 9} have dimensionless numbers of 10, 100, 1000 for $Re$, and 2.5e-1, 2,5e+0, and 2.5e+1 for $Da$, respectively.}}
\label{Fig:G2_EX159_Velocity_Center}
\end{figure}

{The velocity} profiles on the reference lines (Figure~\ref{Fig:G2_EX159_Velocity_Center}) 
present more distinctiveness between the models. 
Regarding $u_y$ profiles on the horizontal reference line 
((left) of Figure~\ref{Fig:G2_EX159_Velocity_Center}), the D-B, D-B-F, and N-S models become more asymmetric compared to more symmetric Brinkman model with increasing values of $Re\times Da$. 
More importantly, unlike \texttt{GROUP I} {where the D-B is overlapping the Brinkman, here} the D-B is almost identical to the N-S  
for three \textbf{\textsl{Tests}}. {As a homogenized model from the N-S \citep{CheZ2001}, the D-B model is the case when $c_F=0$ in (c) $\dfrac{c_F}{\sqrt{Da}}$. %in \eqref{eq:dimless-full-nonlinear}.
} %examples. 
Still, the D-B-F 
is differentiated from the D-B and N-S models, not to mention the Brinkman model. Through the last mesh refinements of \textbf{\textsl{Test} 9} (see Figure~\ref{Fig:G2_EX9_MESH}), we confirm  
the Brinkman demonstrates the symmetry, and there are notable deviations between the three models: the Brinkman, D-B and D-B-F. {We also identify that the maximum adaptive grid refinements occur locally near the top boundary where the discontinuity condition is applied.} {As similar to \texttt{GROUP I}, the D-B model has the most nodes with DoFs.}
\begin{figure}[htbp]
\centering
\subfloat[Brinkman \textbf{\textsl{Test} 9}]
{\includegraphics[width=0.3\textwidth]{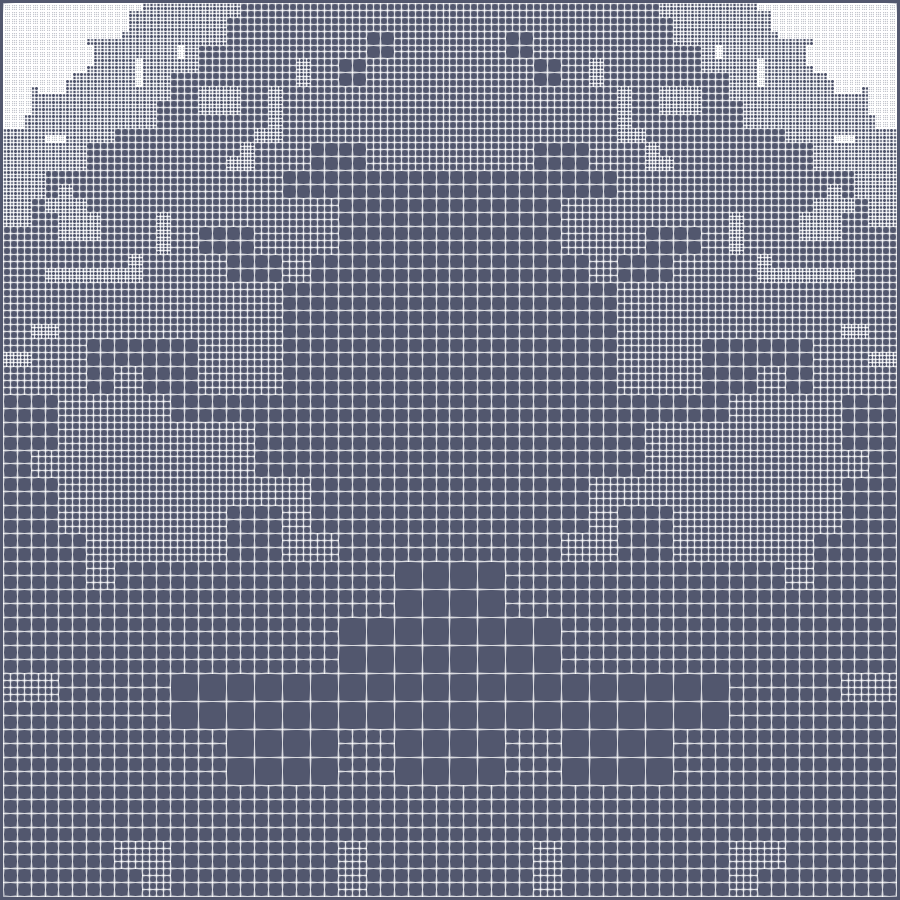}}\hspace*{0.05em}
\subfloat[D-B \textbf{\textsl{Test} 9}]
{\includegraphics[width=0.3\textwidth]{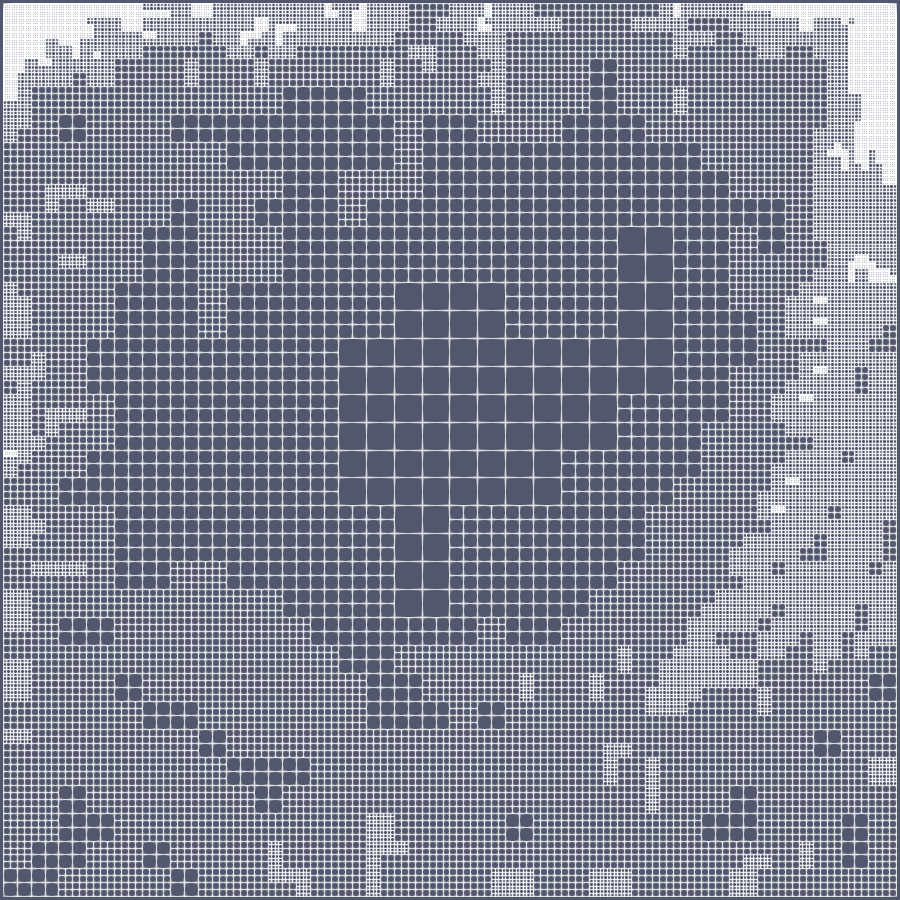}}\hspace*{0.05em}
\subfloat[D-B-F \textbf{\textsl{Test} 9}]
{\includegraphics[width=0.3\textwidth]{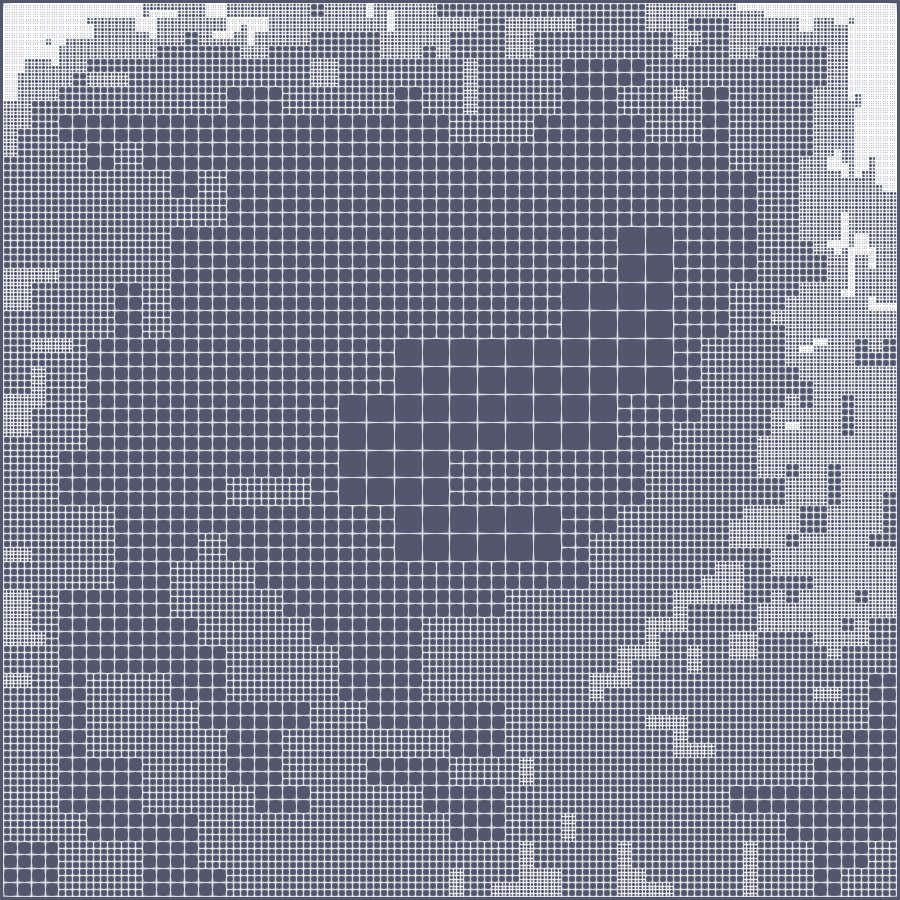}}
\caption{[\texttt{GROUP II}] \textbf{\textsl{Test} 9}.  {Adaptive mesh refinements} of the last ($4^{th}$) cycle using the Kelly error estimator: (a) Brinkman, (b) Darcy-Brinkman, and (c) Darcy-Brinkman-Forchheimer. {The total numbers of DoFs are (a) 175303 (DoFs for $\mathbf{u}_{h}$:155674, for $p_h$:19629), (b) 229747 (DoFs for $\mathbf{u}_{h}$:204006, for $p_h$:25741), (c) 215585 (DoFs for $\mathbf{u}_{h}$:191430, for $p_h$:24155)}.}
\label{Fig:G2_EX9_MESH}
\end{figure}

%Note in the same manner that 
Table~\ref{tab:Num_Newton_g2_Ex159} shows the Newton iteration numbers, which are identical to the ones in \texttt{GROUP I}. Unlike  \texttt{GROUP I}, when the Augmented Lagrangian term is set with $\gamma=0.0$ for the D-B-F model, {the iterative solver of FGMRES fails to converge. This implies that the increased nonlinearity requires the balanced Augmented Lagrangian term, and an appropriate value of $\gamma$ must be provided.} 
\begin{table}[!h]
\begin{center}
\caption{[\texttt{GROUP II}] Total number of Newton iterations for \textbf{\textbf{\textsl{Test}}} \textbf{1}, \textbf{5}, and \textbf{9}.}
\label{tab:Num_Newton_g2_Ex159}
\begin{tabular}{rccc}
\toprule
\textbf{\textbf{\textsl{Test}}} & Brinkman & D-B & D-B-F \\
\midrule
\textbf{\textbf{\textsl{Test}}} \textbf{1} & $3$ & 4 & {$4$} \\
\textbf{\textbf{\textsl{Test}}} \textbf{5} & $3$ & 4 & $4$ \\
\textbf{\textbf{\textsl{Test}}} \textbf{9} & $3$ & 4 & $4$ \\
\bottomrule
\end{tabular}
\end{center}
\end{table}

\subsubsection*{{In-depth study} for \texttt{GROUP II}}

{In this subsection}, the remaining \textbf{\textbf{\textsl{Tests}}} in \texttt{GROUP II} are investigated with in-depth study on the departure of the full Darcy-Brinkman-Forchheimer model from other models, including the Navier-Stokes, 
highlighting the Forchheimer term for inertial resistance. 

\subparagraph*{{$\bullet$ \textbf{\textsl{Test 4}} \& \textbf{\textsl{Test 6}}:}} 
\textbf{\textbf{\textsl{Test 4}}} and \textbf{\textbf{\textsl{Test 6}}} are in the transient regime for the porous medium with fixed $Re=100$. As can be seen in Figure~\ref{Fig:G2_EX46}, particularly from the Brinkman model, 
the increased $Da$ intensifies the vorticity, which can be found in other models as well. \textbf{\textbf{\textsl{Test 4}}} has a slight more difference in the location of primary eddy between the D-B and the D-B-F models (see (b) and (c) in Figure~\ref{Fig:G2_EX46}) than  the that of \textbf{\textbf{\textsl{Test 6}}} (see (e) and (f) in Figure~\ref{Fig:G2_EX46}), exhibiting more effect of inertial resistance compared to the viscous resistance, which is via the Forchheimer term. %In the same context, it is clear that the inertial resistance of \textbf{\textsl{Test 4}} is greater than that of \textbf{\textsl{Test 6}}.
\begin{figure}[htbp]
\centering
\subfloat[Brinkman \textbf{\textsl{Test} 4}]
{\includegraphics[width=0.3\textwidth]{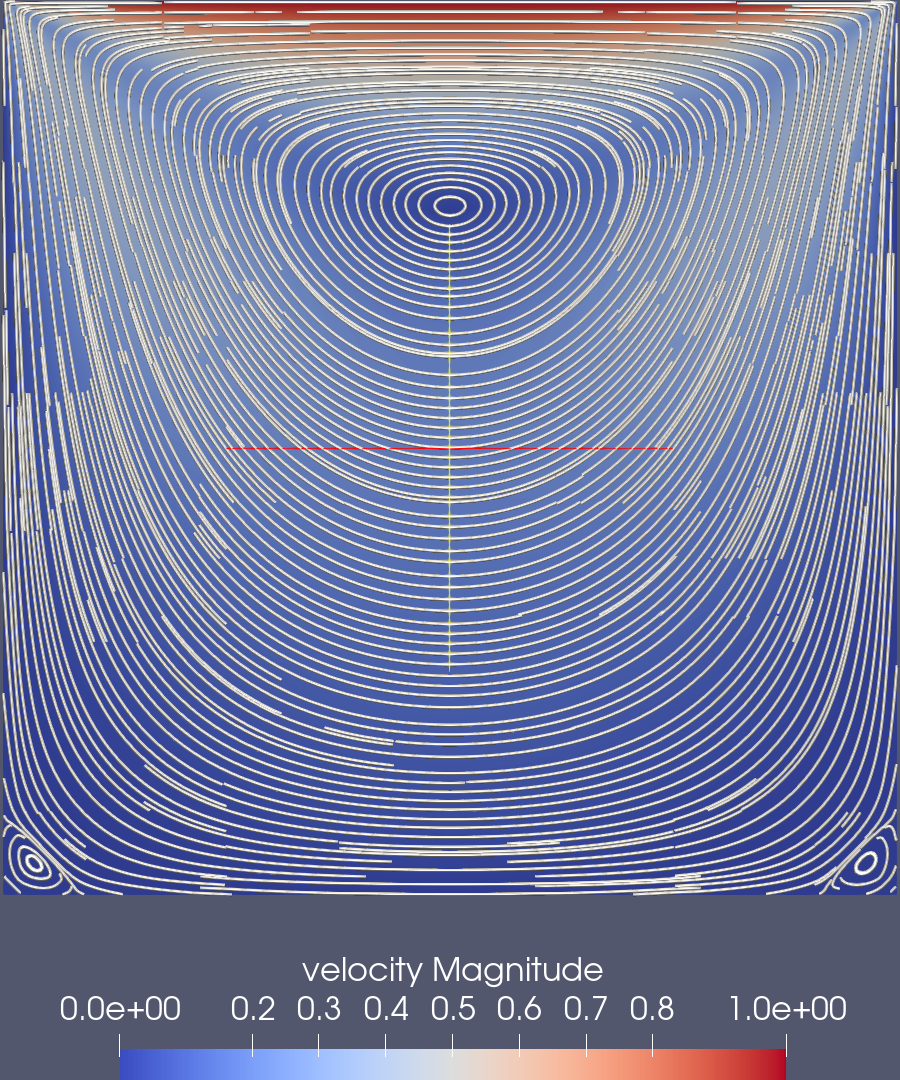}}\hspace*{0.05em}
\subfloat[D-B \textbf{\textsl{Test} 4}]
{\includegraphics[width=0.3\textwidth]{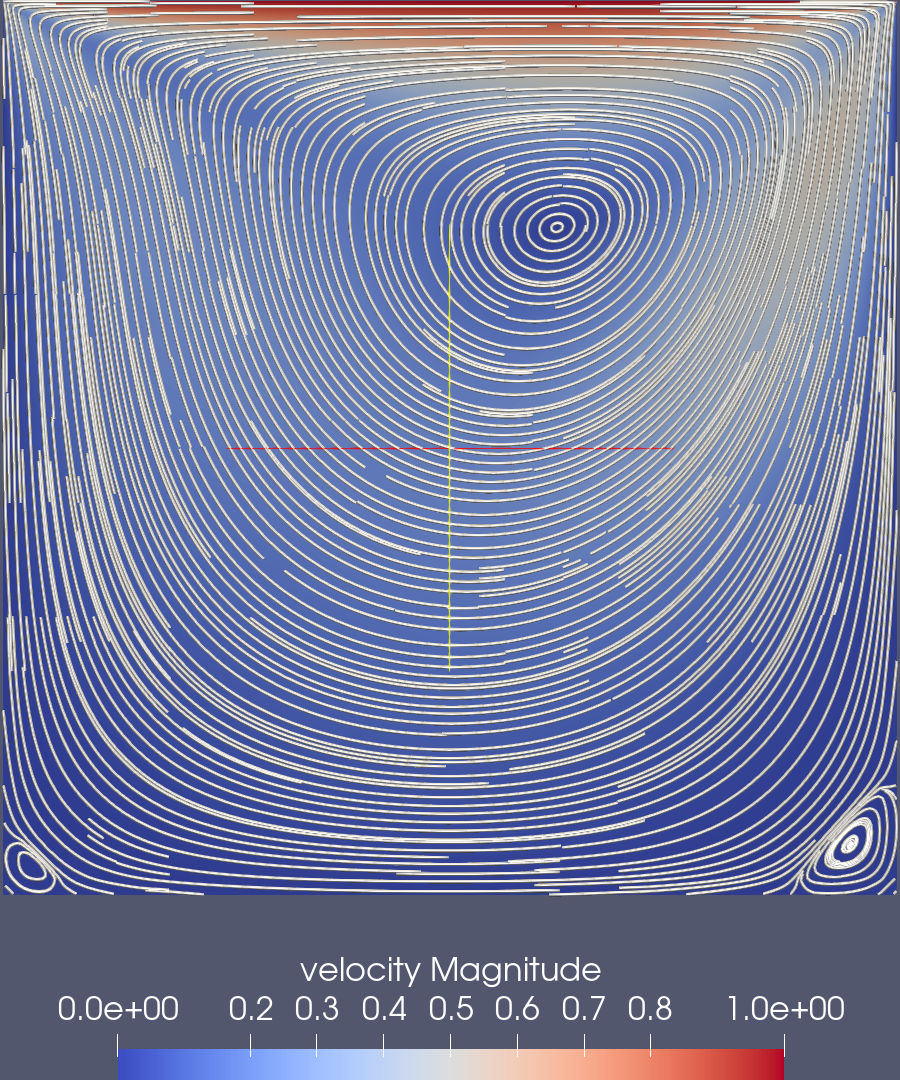}}\hspace*{0.05em}
\subfloat[D-B-F \textbf{\textsl{Test} 4}]
{\includegraphics[width=0.3\textwidth]{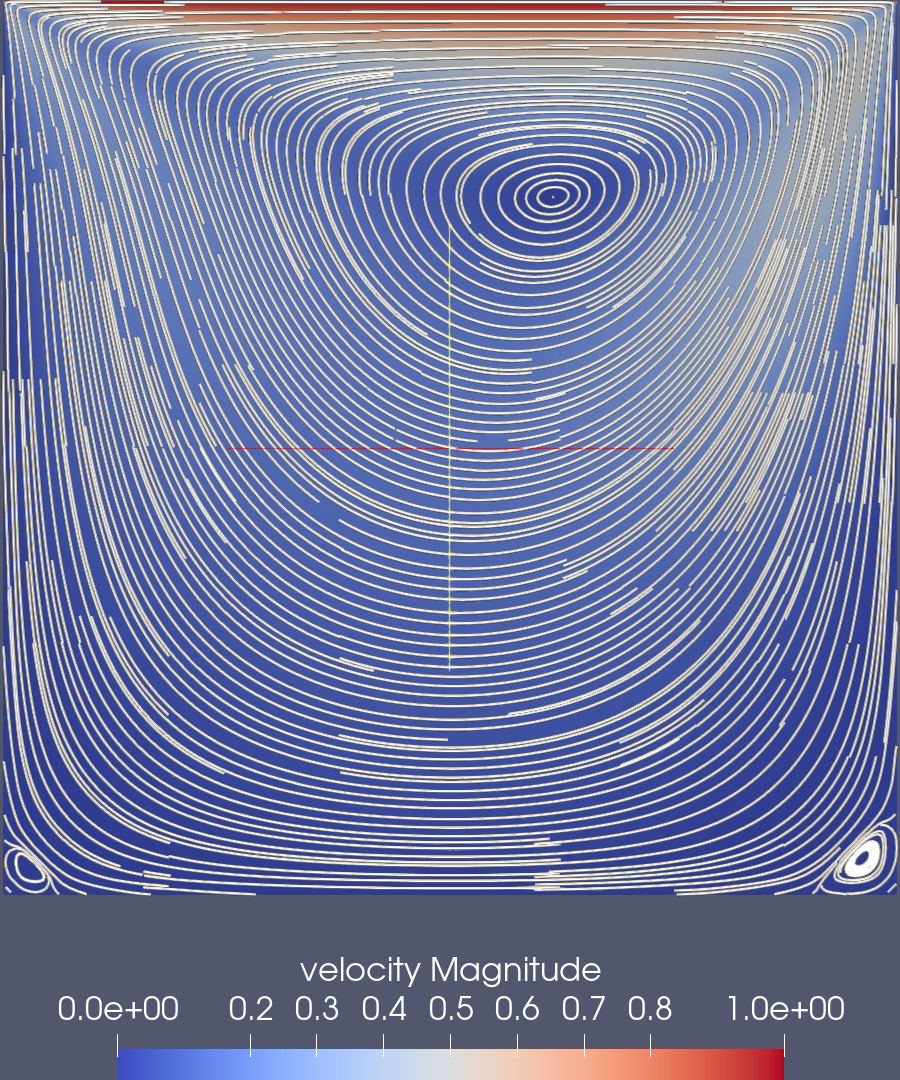}}\\
\subfloat[Brinkman \textbf{\textsl{Test} 6}]
{\includegraphics[width=0.3\textwidth]{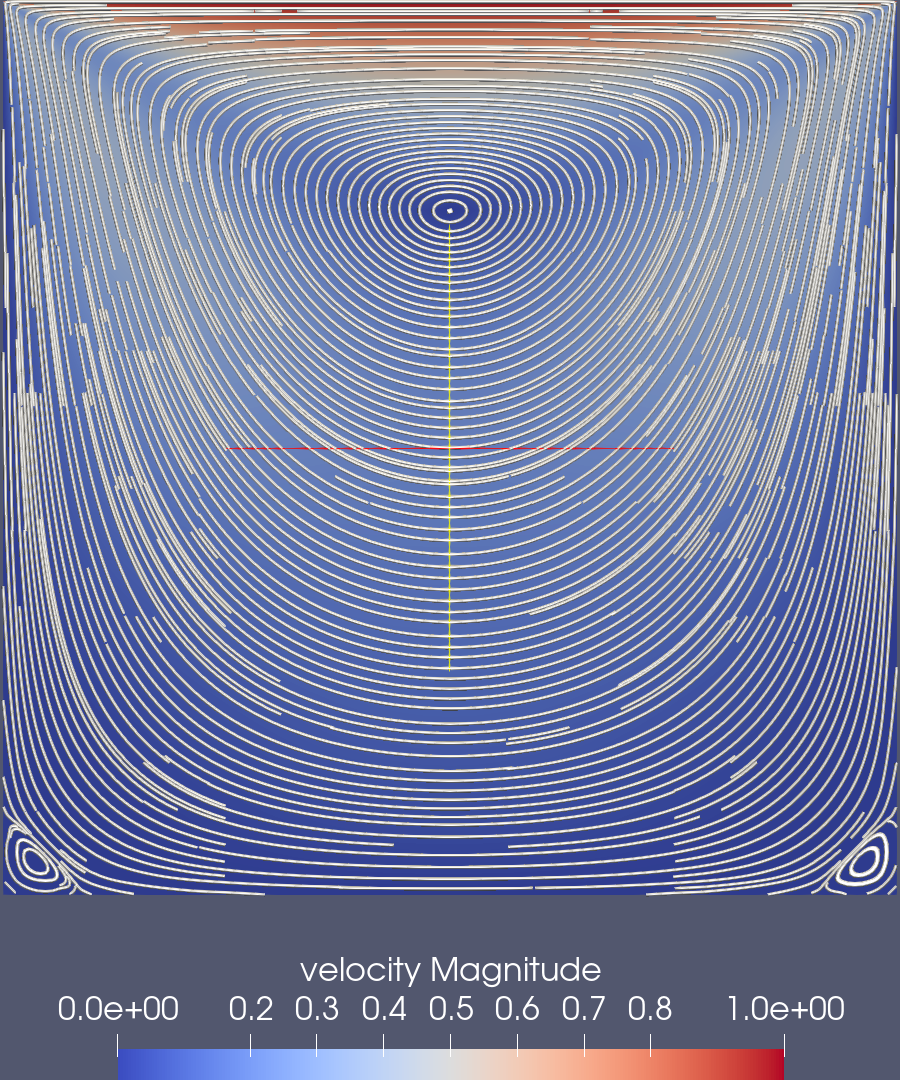}}\hspace*{0.05em}
\subfloat[D-B \textbf{\textsl{Test} 6}]
{\includegraphics[width=0.3\textwidth]{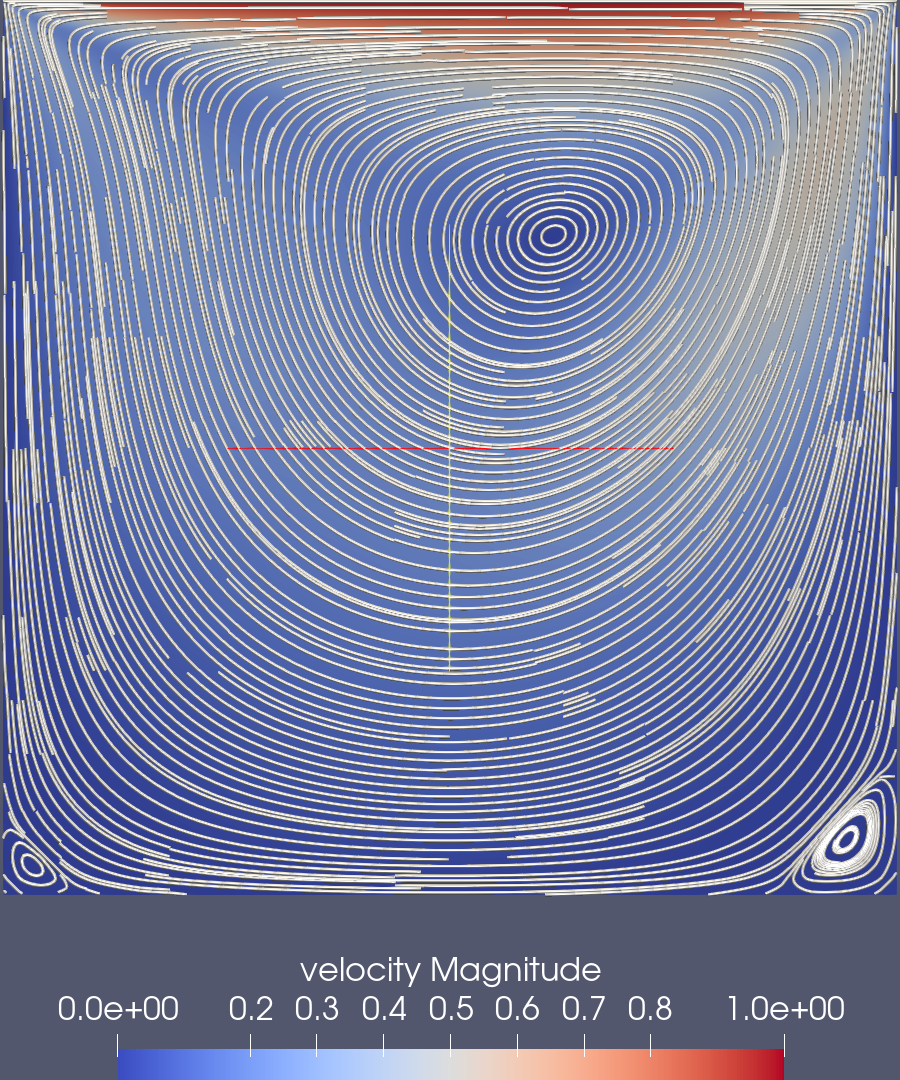}}\hspace*{0.05em}
\subfloat[D-B-F \textbf{\textsl{Test} 6}]
{\includegraphics[width=0.3\textwidth]{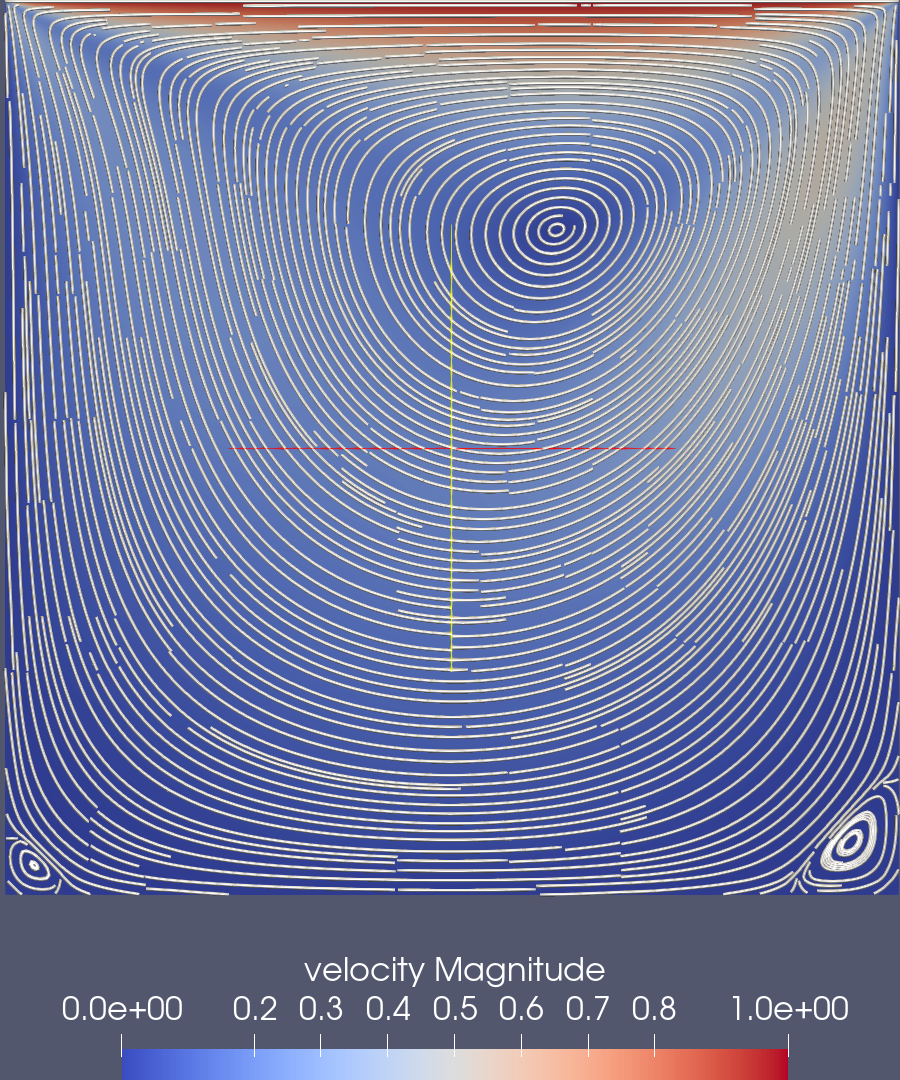}}
\caption{[\texttt{GROUP II}] \textbf{\textsl{Test} 4 and 6}. Velocity streamline: {(left)} Brinkman, {(center)} Darcy-Brinkman, and (right) Darcy-Brinkman-Forchheimer. {\textbf{\textsl{Test} 4 and 6} have dimensionless numbers of 100 for $Re$ and 2.5e-1 and 2.5e+1 for $Da$, respectively.}}
\label{Fig:G2_EX46}
\end{figure}
\begin{figure}[htbp]
\centering
\subfloat[Brinkman \textbf{\textsl{Test} 2}]
{\includegraphics[width=0.3\textwidth]{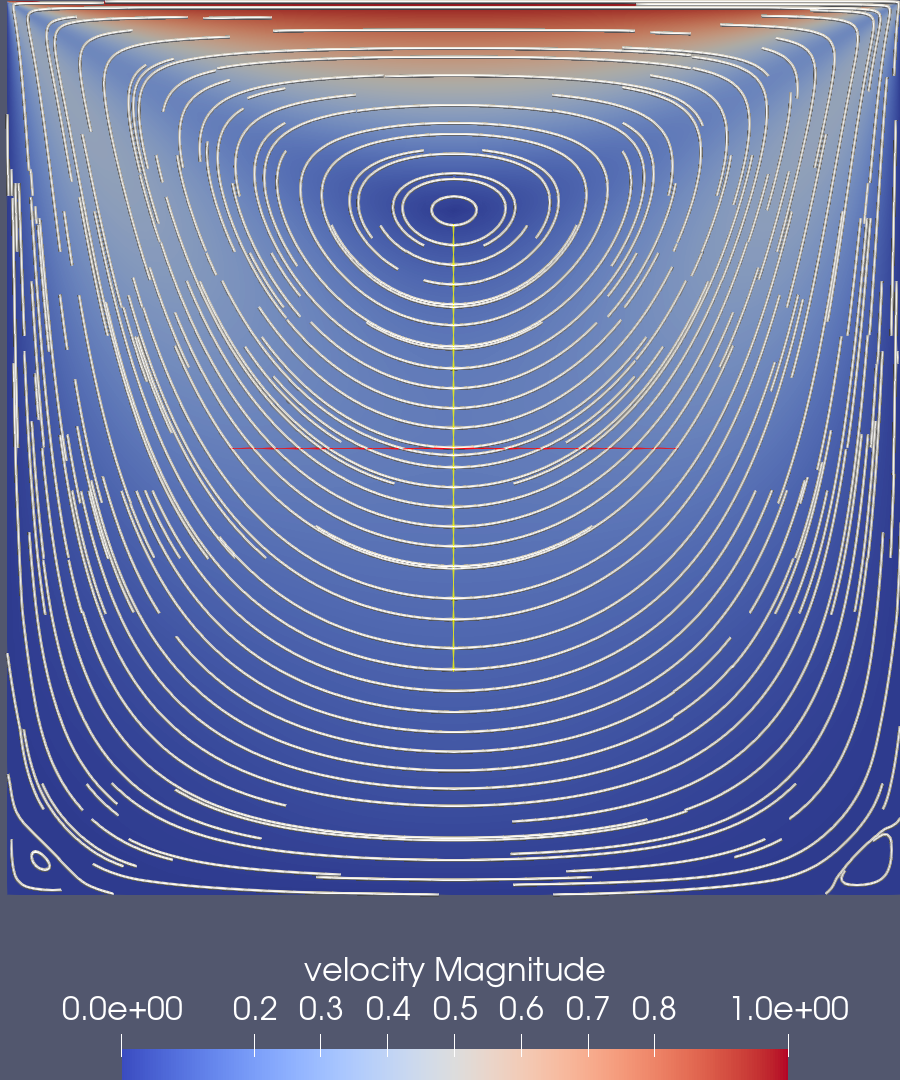}}\hspace*{0.05em}
\subfloat[D-B \textbf{\textsl{Test} 2}]
{\includegraphics[width=0.3\textwidth]{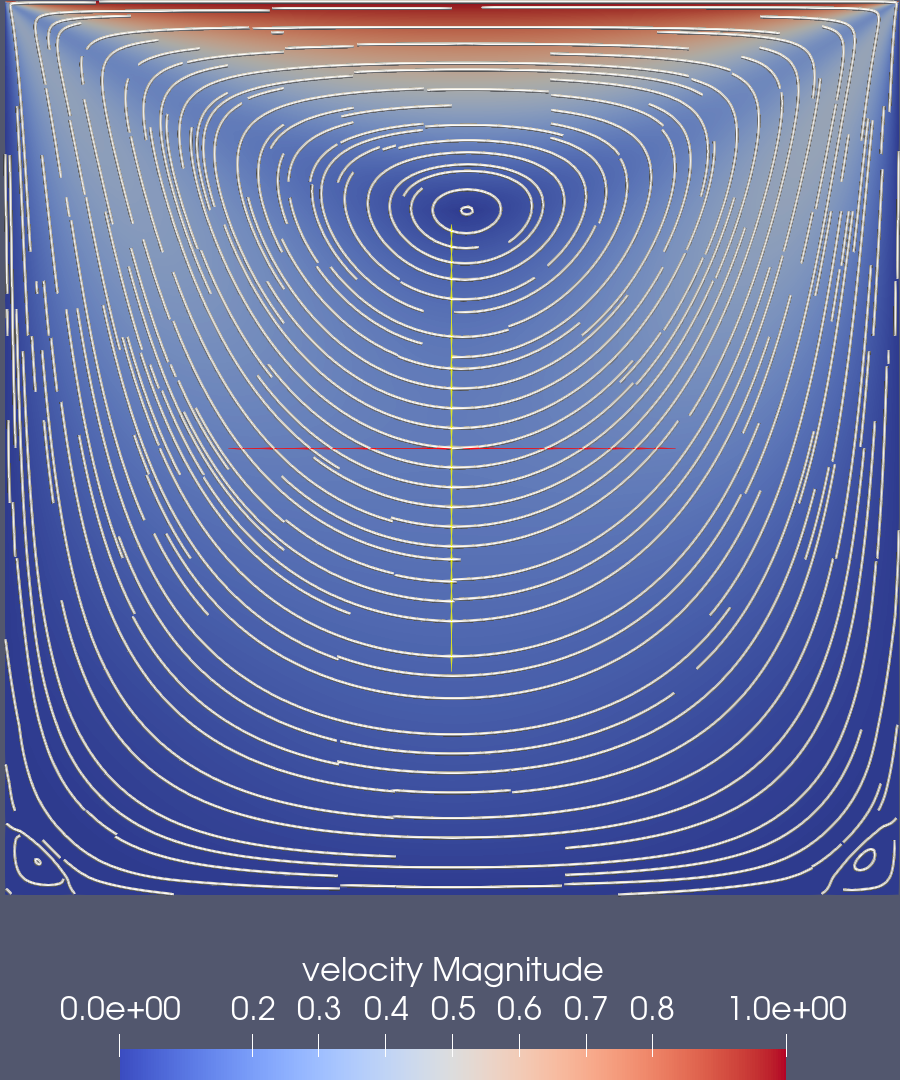}}\hspace*{0.05em}
\subfloat[D-B-F \textbf{\textsl{Test} 2}]
{\includegraphics[width=0.3\textwidth]{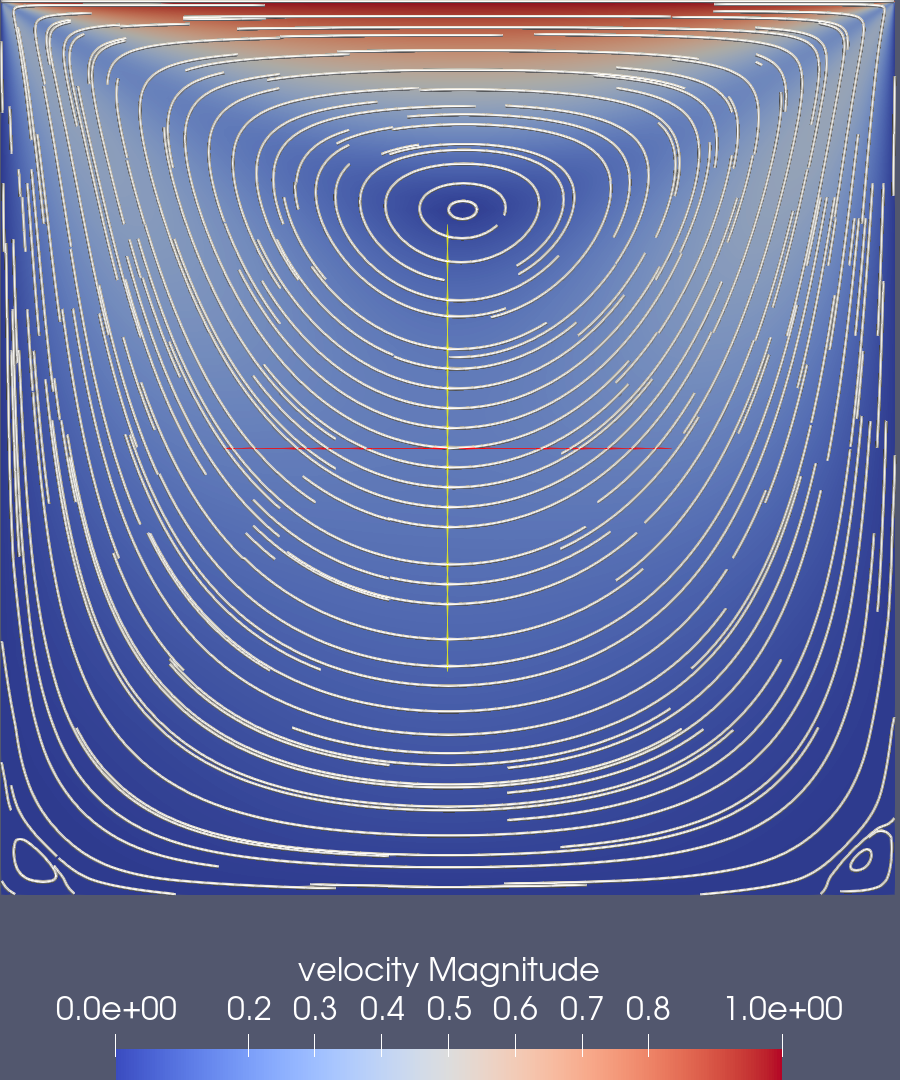}}\\
\subfloat[Brinkman \textbf{\textsl{Test} 8}]
{\includegraphics[width=0.3\textwidth]{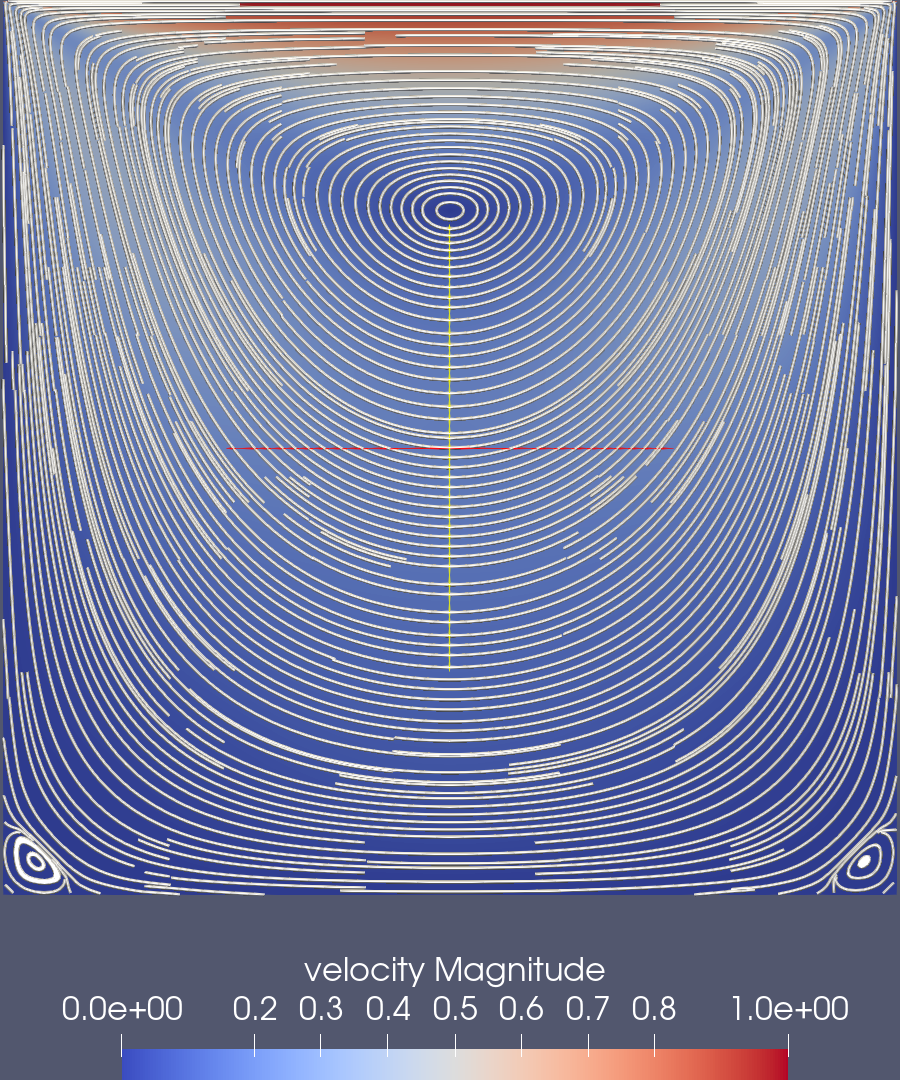}}\hspace*{0.05em}
\subfloat[D-B \textbf{\textsl{Test} 8}]
{\includegraphics[width=0.3\textwidth]{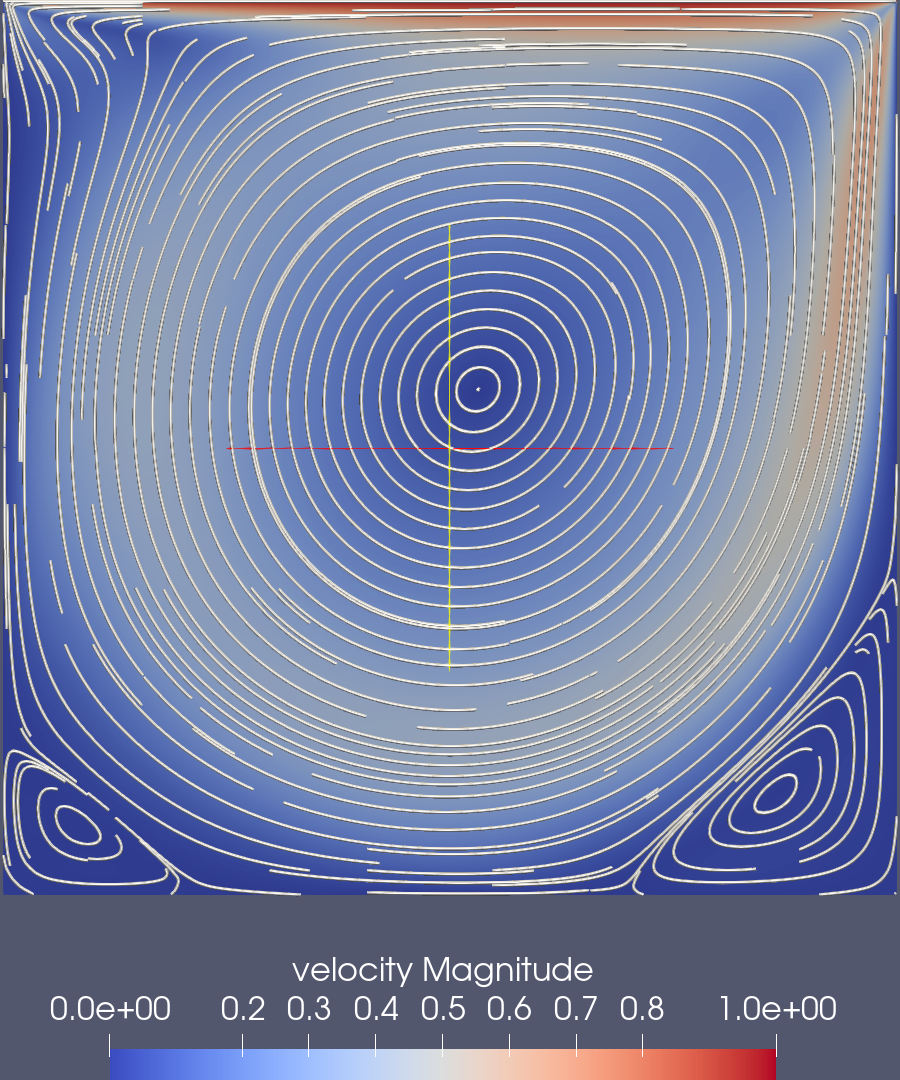}}\hspace*{0.05em}
\subfloat[D-B-F \textbf{\textsl{Test} 8}]
{\includegraphics[width=0.3\textwidth]{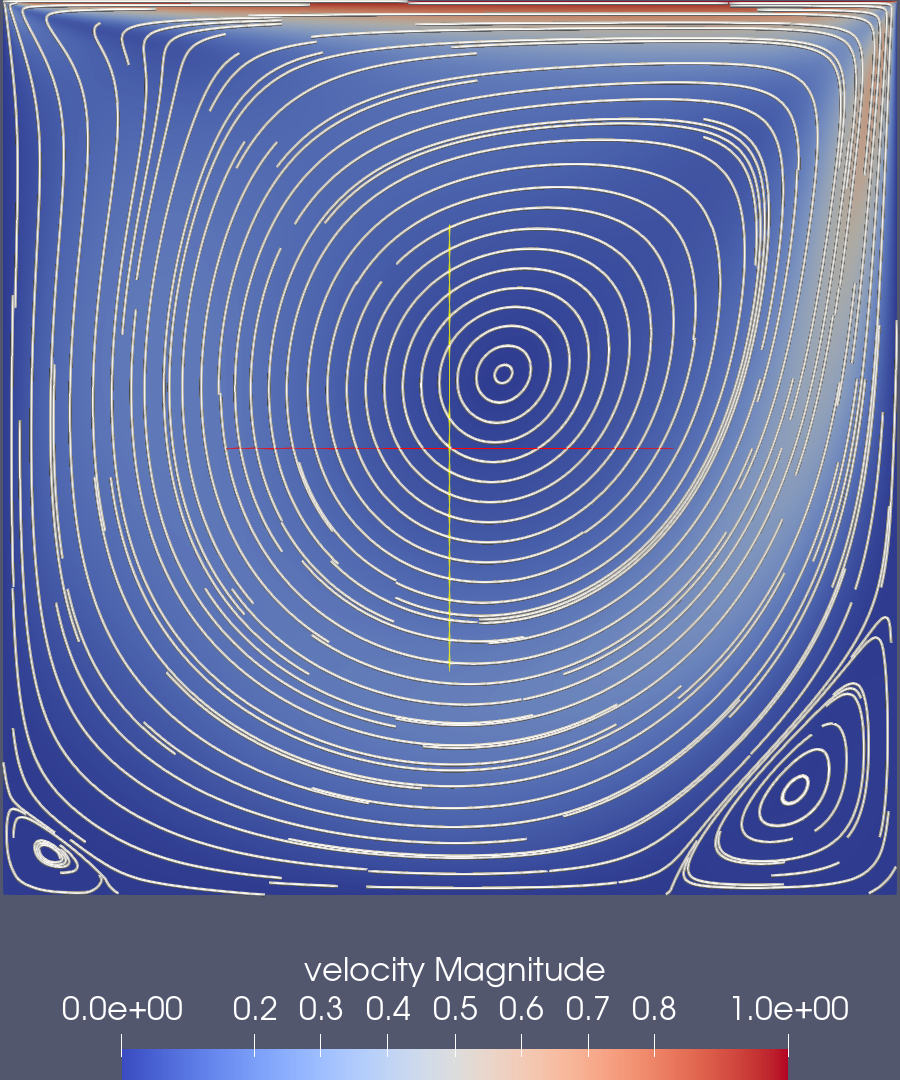}}
\caption{[\texttt{GROUP II}] \textbf{\textsl{Test} 2 and 8}.  Velocity streamline: {(left)} Brinkman, {(center)} Darcy-Brinkman, and (right) Darcy-Brinkman-Forchheimer. {\textbf{\textsl{Test} 2 and 8} have dimensionless numbers of 10 and 1000 for $Re$, respectively, and 2.5e+0 for $Da$.}}
\label{Fig:G2_EX28}
\end{figure}

\subparagraph*{{$\bullet$~\textbf{\textsl{Test 2}} \& \textbf{\textsl{Test 8}}:}}
{From} \textbf{\textsl{Test 2}} and \textbf{\textsl{Test 8}}, we identify that $Re$ variation with fixed $Da$ has more effects on the flow than $Da$ variation with fixed $Re$ does (see Figure~\ref{Fig:G2_EX28} compared to Figure~\ref{Fig:G2_EX46}). 
Much difference in streamlines between examples, i.e., \textbf{\textbf{\textsl{Test 2}}} and \textbf{\textbf{\textsl{Test 8}}} within each model, is due to different $Re$ values, resulting in the change of flow regimes; \textbf{\textbf{\textsl{Test 2}}} is in the laminar, whereas \textbf{\textbf{\textsl{Test 8}}} is in the turbulence for a typical porous medium.
{Albeit the same increasing ratio of $100$ applied for $Re$ as for $Da$ in the previous case (see Table~\ref{tab:Examples}), %this is due to 
we find that $Re$ that has more influence (see \eqref{eq:dimless-full-nonlinear}) for this lid-driven cavity problem under the unit velocity (i.e., the maximum velocity), %and which also implies 
implying that the flow regime change is more abrupt and instantaneous~\cite{BearJ1972}}.

Note also that when the intensity of vorticity is increased by either $Da$ or $Re$, the growth of vortices is somewhat restricted for the D-B-F model compared to the D-B model, where the only difference is the existence of inertial resistance. 
As seen in Figure~\ref{Fig:G2_EX46} and \ref{Fig:G2_EX28}, the size of vortices from the D-B-F is smaller than those from the D-B, which is distinctively shown in the turbulent regime. 
\subparagraph*{{$\bullet$~\textbf{\textbf{\textsl{Test 7}}}:}}
\textbf{\textbf{\textsl{Test 7}}} is selected to investigate further the effect of the inertial resistance to the flow. 
Note that {the ratio of} (c) $\dfrac{c_F}{\sqrt{Da}}$ in \eqref{eq:dimless-full-nonlinear} compared to $Re$, i.e., the value of $\dfrac{c_F}{\sqrt{Da}}$, \textbf{\textbf{\textsl{Test 7}}} has the largest for all \textbf{\textbf{\textsl{Tests}}} in \texttt{GROUP I} and \texttt{GROUP II} (see Table~\ref{tab:Examples}). In addition, %it has the \hl{farthest} 
the distance from 
the Navier-Stokes system, i.e.,  $Da\rightarrow\infty$, in terms of  (b) $\dfrac{1}{ReDa}$ and (c) $\dfrac{c_F}{\sqrt{Da}}$ in \eqref{eq:dimless-full-nonlinear}, \textbf{\textbf{\textsl{Test 7}}} is the farthest among all \textbf{\textbf{\textsl{Tests}}} in \texttt{GROUP II}.  
\begin{figure}[htbp]
\centering
\subfloat[Brinkman {\textbf{\textsl{Test} 7}}] 
{\includegraphics[width=0.3\textwidth]{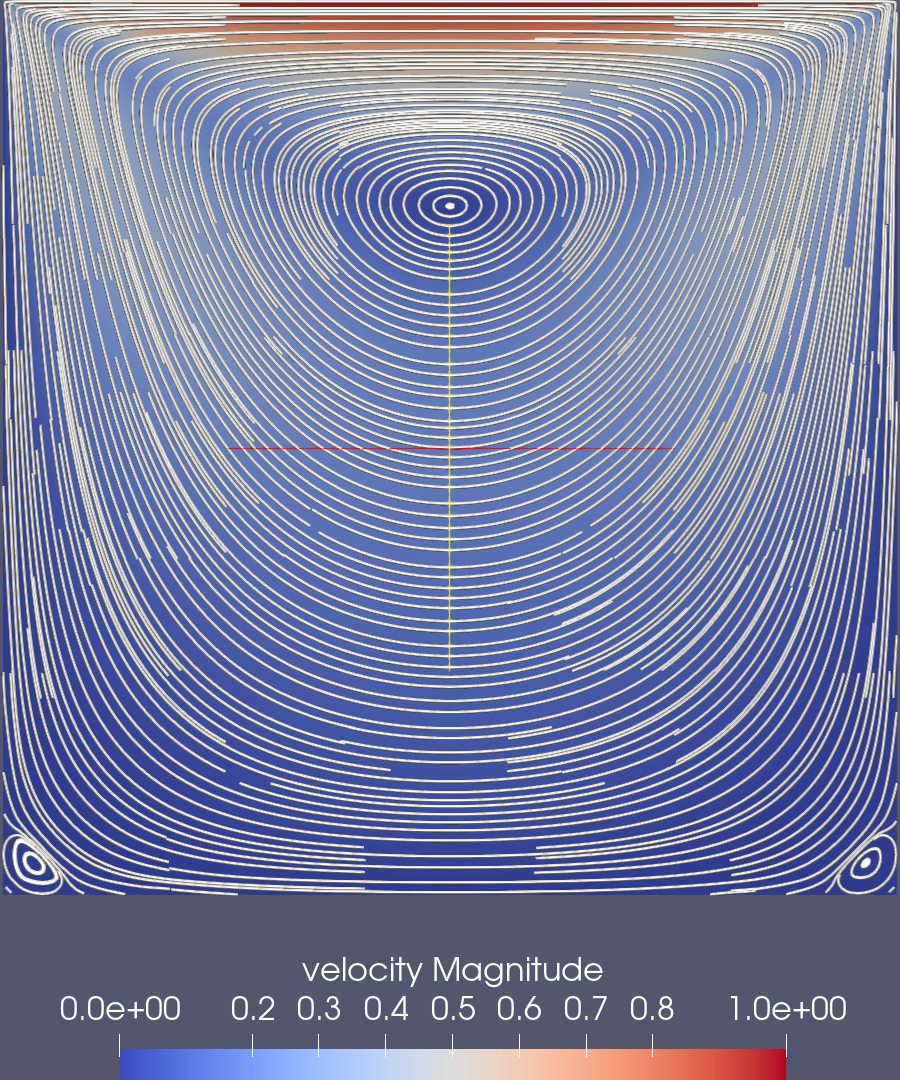}}\hspace*{0.05em}
\subfloat[D-B {\textbf{\textsl{Test} 7}}] 
{\includegraphics[width=0.3\textwidth]{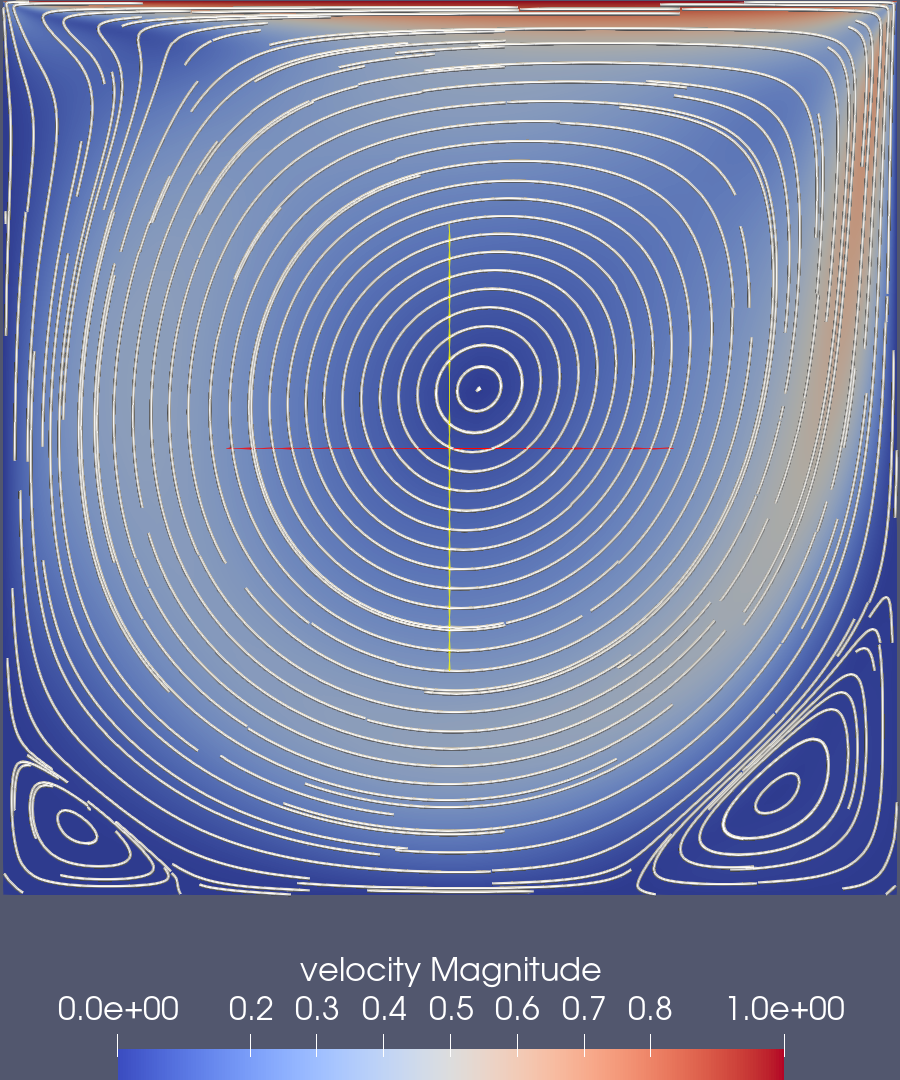}}\hspace*{0.05em}
\subfloat[D-B-F {\textbf{\textsl{Test} 7}}] 
{\includegraphics[width=0.3\textwidth]{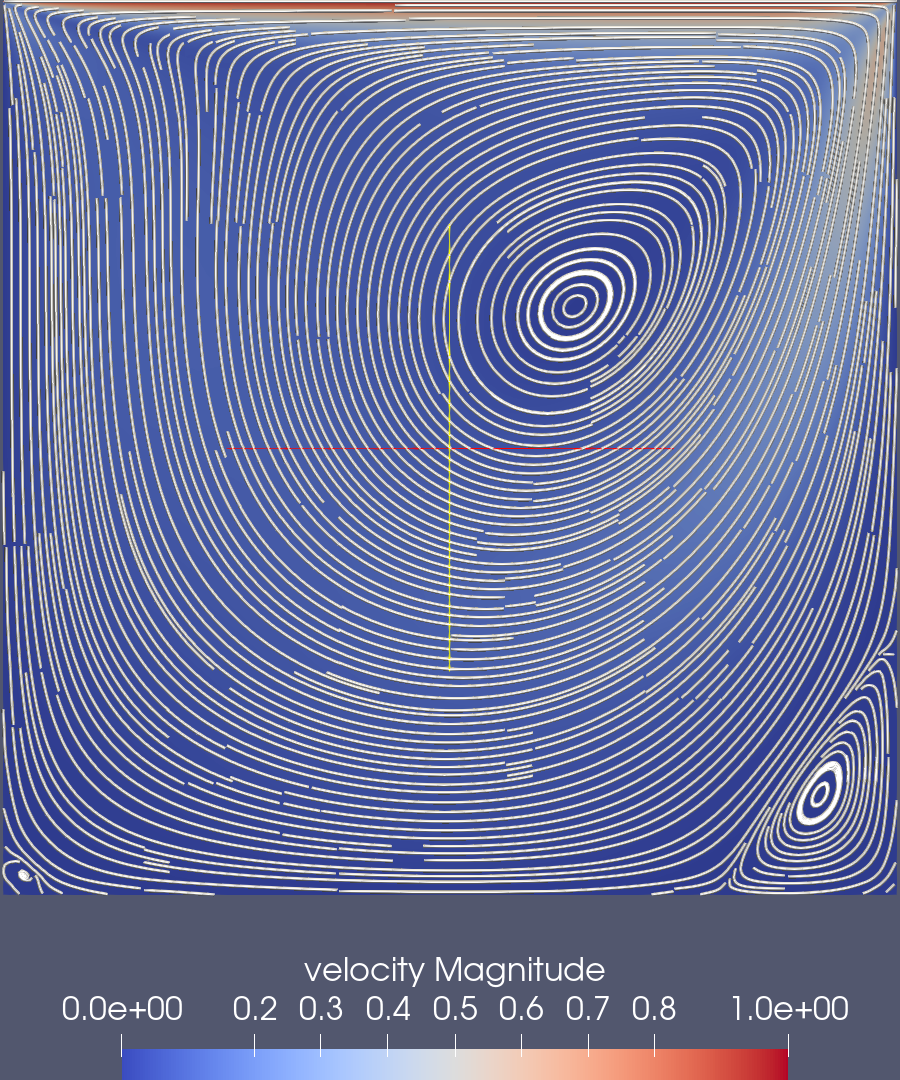}}\\
\subfloat[Velocity on the center line]
{\includegraphics[width=0.75\textwidth]{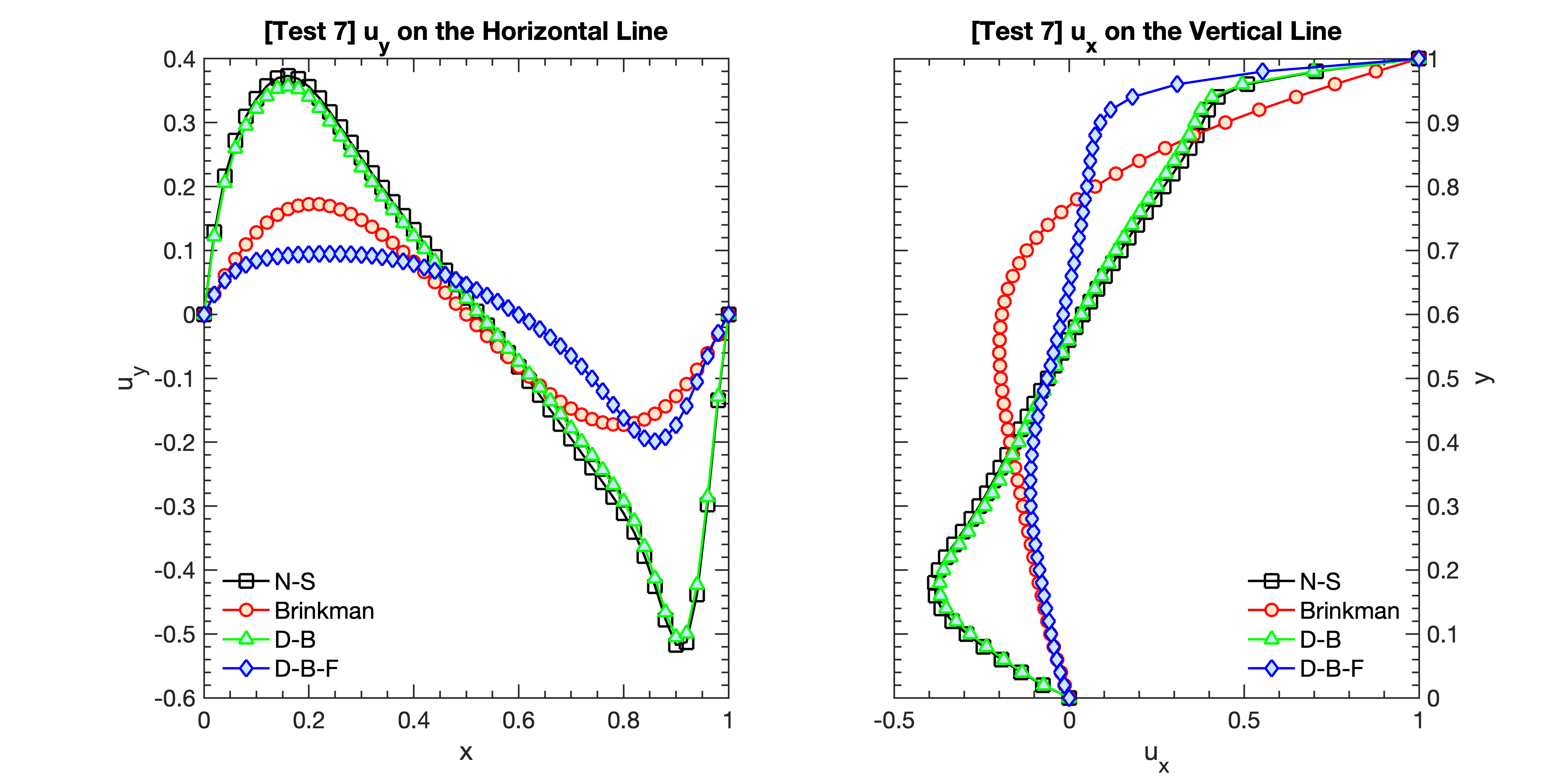}}
\caption{[\texttt{GROUP II}] \textbf{\textsl{Test} 7}. {Velocity} profiles: (a) Brinkman, (b) Darcy-Brinkman, (c) Darcy-Brinkman-Forchheimer, and (d) velocity on the reference lines using the Brinkman, Darcy-Brinkman, and Darcy-Brinkman-Forchheimer models, compared with the Navier-Stokes model. {\textbf{\textsl{Test} 7} has dimensionless numbers of 1000 for $Re$ and 2.5e-1 for $Da$.}}
\label{Fig:G2_EX7}
\end{figure}

Figure~\ref{Fig:G2_EX7} illustrates 
the distinguishing effect of the inertial resistance. Similarly, as depicted in Figure~\ref{Fig:G2_EX159_Velocity_Center}, the deviation of the full nonlinear D-B-F model is even more distinctive than the previous ones. 
When we compare the D-B with D-B-F with the streamline results and  velocities on the reference lines in \textbf{\textsl{Test~7} and \textbf{\textsl{Test}~8}}, %In contrast, 
%the velocities of the D-B model on the reference lines 
%compared to \textbf{\textsl{Test} 7}, 
 the D-B-F model in \textbf{\textsl{Test} 8} with the increased $Da$ number illustrates more similar streamlines to the D-B, possessing larger area for the first eddy near the cross %as the N-S model 
 (see (e), (f) in Figure~\ref{Fig:G2_EX28} and (b), (c) in Figure~\ref{Fig:G2_EX7}). 
 
\subparagraph*{{$\bullet$~\textit{Forchheimer coefficient, $c_F=0.4$ vs. $c_F=0.6$}:}}
We further test with two different values of $c_F$ for the D-B-F model: $c_F=0.4$ and $c_F=0.6$. Although it is not severe, the effect of the Forchheimer term can be found, 
using the same $Re$ and $Da$ used in \textbf{\textsl{Test} 7} and \textbf{\textsl{Test} 9} in \texttt{GROUP II}. From (a) and (b) in Figure~\ref{Fig:G2_EX79_DBF_CF}, the delay effect resulting from the inertial resistance is also confirmed; 
The location of the primary eddy center of (b) is higher than that of (a), occupying a smaller area.  Note that all the vortices are smaller in (b) and (d) than those of (a) and (c) due to the increased value of $c_F=0.6$. From Figure~\ref{Fig:G2_EX79_DBF_CF}, when we compare \textbf{\textsl{Test} 7} with \textbf{\textsl{Test} 9} with the same $c_F$, i.e., $c_F=0.4$ or $c_F=0.6$, we also confirm that the inertial resistance compared to $Da$ and $Re$ for \textbf{\textsl{Test} 7} is relatively larger than that for \textbf{\textsl{Test} 9}.  

\begin{figure}[htbp]
\centering
\subfloat[D-B-F \textbf{\textsl{Test} 7}: $c_F=0.4$]
{\includegraphics[width=0.35\textwidth]{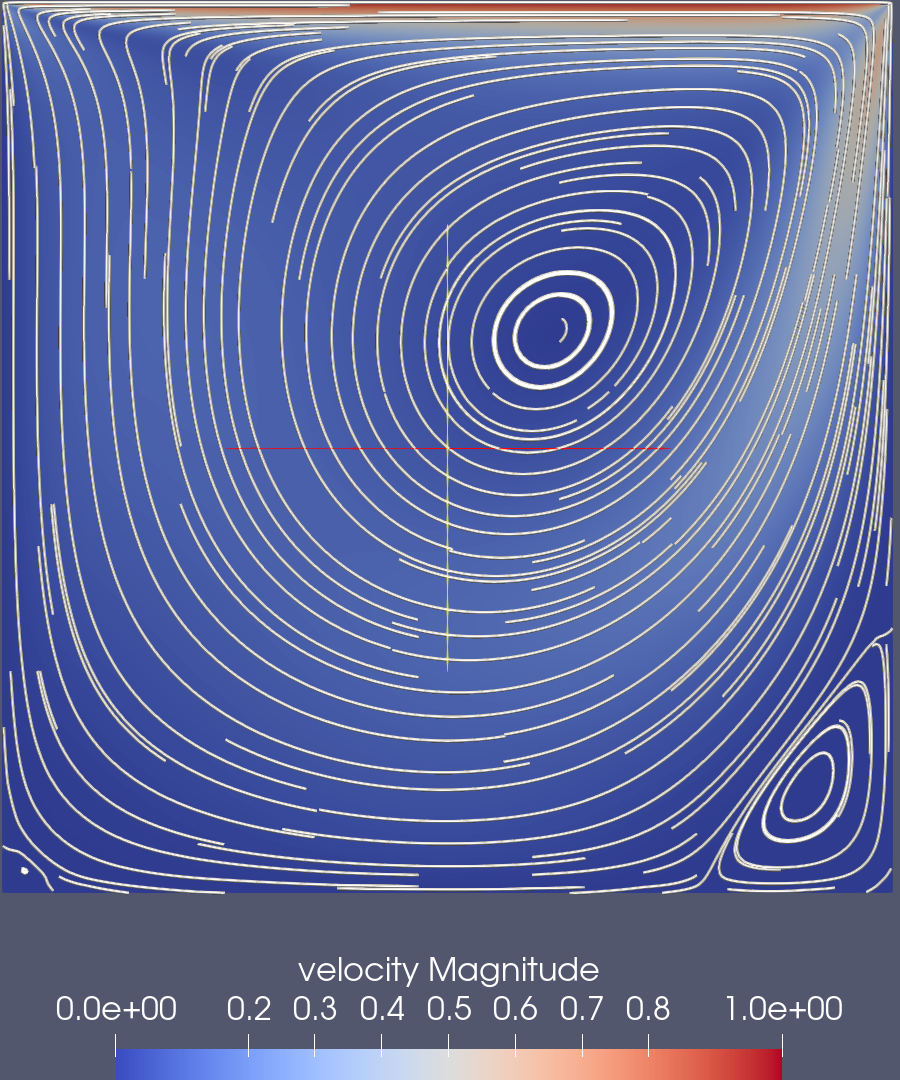}}\hspace*{0.05em}
\subfloat[D-B-F \textbf{\textsl{Test} 7}: $c_F=0.6$]
{\includegraphics[width=0.35\textwidth]{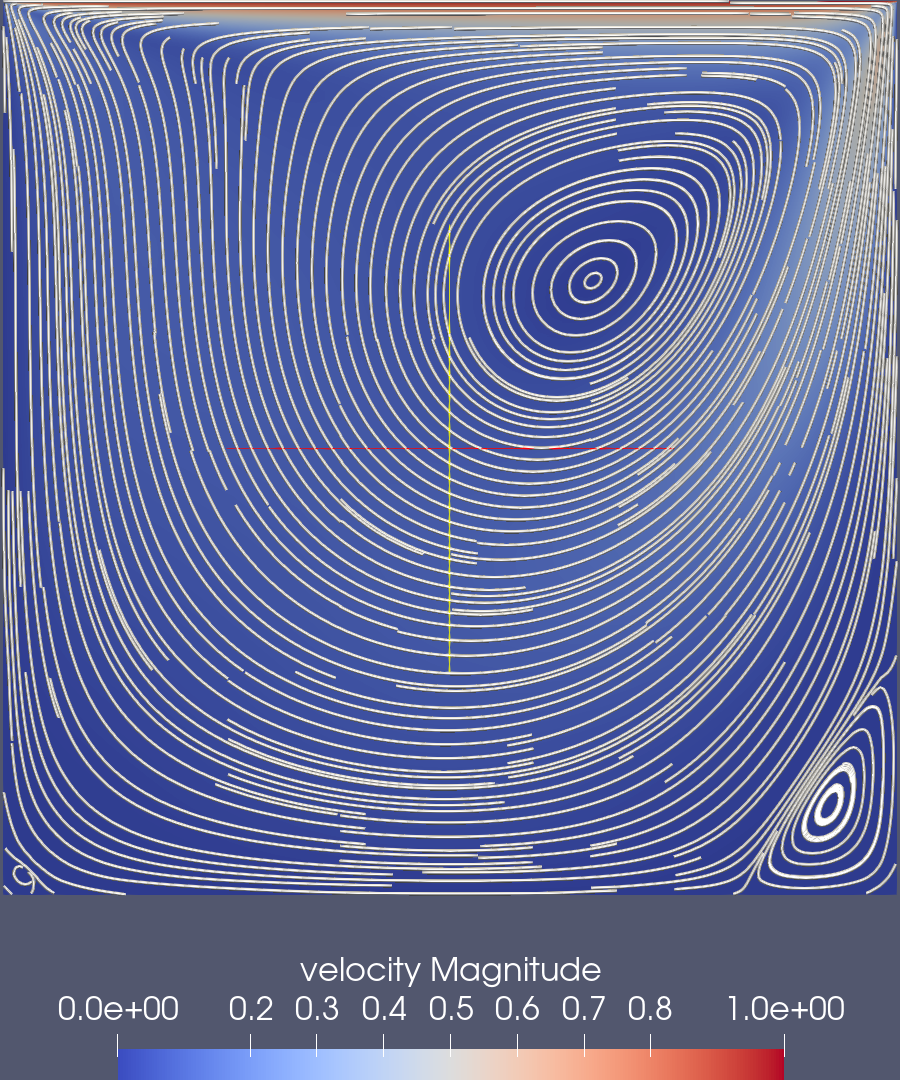}}\\
\subfloat[D-B-F \textbf{\textsl{Test} 9}: $c_F=0.4$]
{\includegraphics[width=0.35\textwidth]{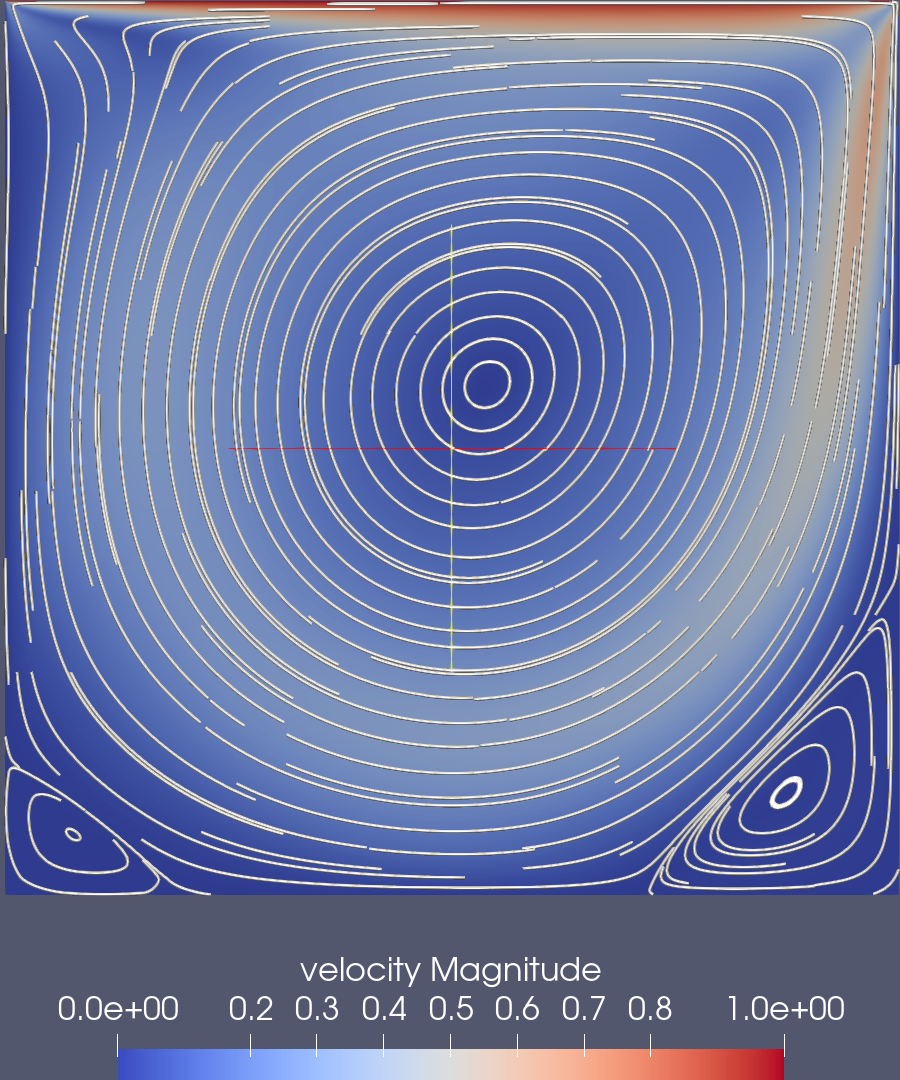}}\hspace*{0.05em}
\subfloat[D-B-F \textbf{\textsl{Test} 9}: $c_F=0.6$]
{\includegraphics[width=0.35\textwidth]{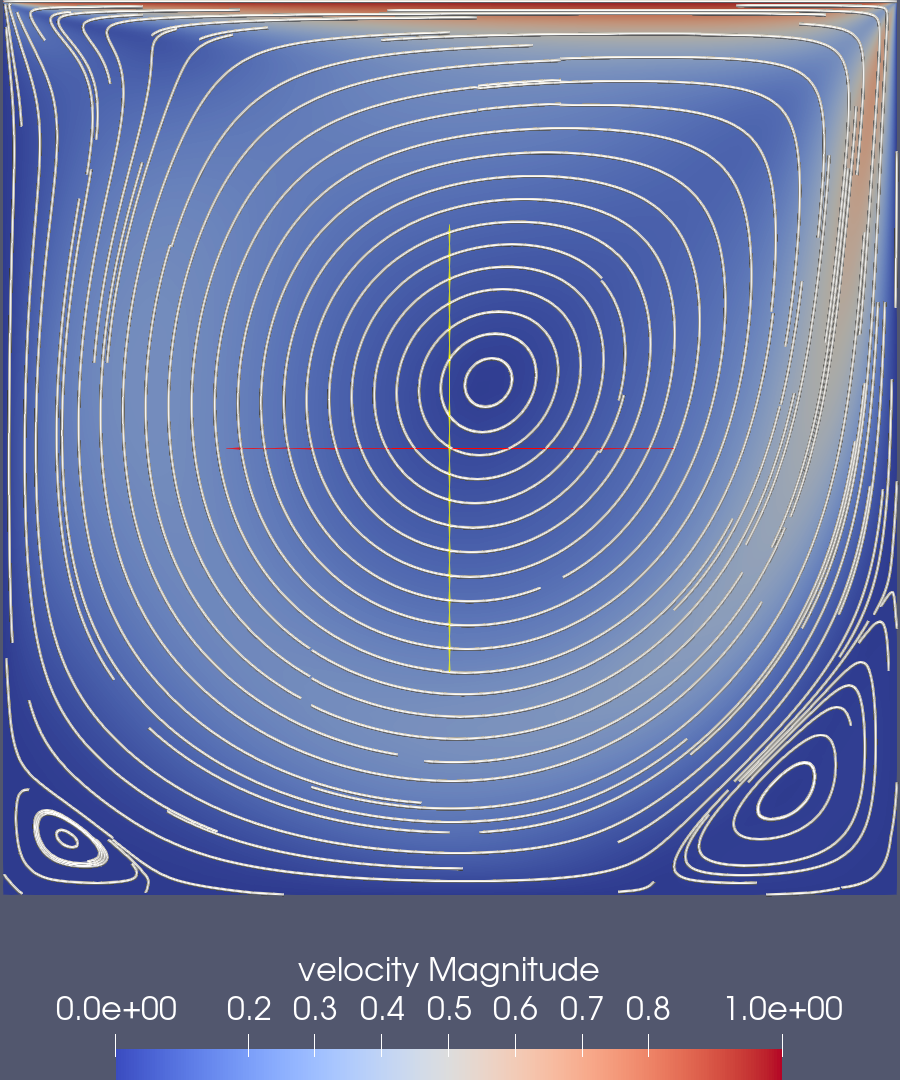}}\\
\caption{[\texttt{GROUP II}] \textbf{\textsl{Test} 7 and 9}. Velocity streamlines $c_F=0.4$ (left) vs. $c_F=0.6$ (right). {\textbf{\textsl{Test} 7 and 9} have dimensionless numbers of 1000 for $Re$ and 2.5e-1 and 2.5e+1 for $Da$, respectively.}}
\label{Fig:G2_EX79_DBF_CF}
\end{figure}

\subsubsection{Summary on \texttt{GROUP I} and \texttt{GROUP II} results}
%\hl{In summary for} \texttt{GROUP I}, \hl{the flow is found to be controlled mainly by the viscous resistance, and there are no big differences between the three models: linearized Brinkman, the Darcy-Brinkman, and the Darcy-Brinkman-Forchheimer models.} 
{
Depending on the parameters, three models, i.e., the linearized Brinkman, Darcy-Brinkman, and Darcy-Brinkman-Forchheimer, sometimes behave similarly or depart from one another. For example, the Darcy-Brinkman is close to the Brinkman when $Re\times Da<1.0$ (\texttt{GROUP I}), but it is similar to the Navier-Stokes when $Re\times Da>1.0$ (\texttt{GROUP II}). The Brinkman can model correctly within the small $Da$ values (or $Re\times Da < 1.0$), but not in the large $Da$ (or $Re\times Da > 1.0$).  
{The detailed findings for each group are what follows. For $Re\times Da<1.0$, 
\begin{itemize}
\item Few distinctions are found in velocity profiles between the three models, i.e., the linearized Brinkman, Darcy-Brinkman, and Darcy-Brinkman-Forchheimer. Within these models, we identify the flow is controlled mainly by the viscous resistance with small $Da$ values.
\item The linearized Brinkman and Darcy-Brinkman (without the inertial resistance term) models are closer to each other as little convection exists and they have the viscous resistance only.
\item {The linearized Brinkman can also be replaced with the Darcy-Brinkman-Forchheimer in \texttt{GROUP I} without much error. %However, 
Little deviations that can be discerned by the Darcy-Brinkman-Forchheimer model have started, implying the intermediate regime or weak inertial regime \cite{ZimR2004}.}
\end{itemize}
Meanwhile, for $Re\times Da>1.0$, 
\begin{itemize}
\item The distinctions of each model are clearly shown  
with its different location and size of primary, secondary, and tertiary vortices. The Forchheimer term reduces the effect and size of the vortices with the inertial resistance coefficient, which is known to depend on hydraulic properties of the porous medium. 
\item A rather abrupt change in flow patterns for the Darcy-Brinkman model is identified during the flow regime change, i.e., the laminar to turbulent, when $Re$ values are increased. And the Darcy-Brinkman model is the closest to the Navier-Stokes as it is homogenized from it \cite{CheZ2001}. 
\item %\hl{However,} 
Since the Darcy-Brinkman-Forchheimer model has the inertial resistance, %accompanies the corresponding change of inertial resistance in its effect , 
%it can delay 
the change of flow regime is delayed under the equivalent inertial and viscous forces. 
\end{itemize}
}

%%%%%%%%%%%%%%%%%%%%%%%%%%%%%%%%%%%%%%%%%%%%%%%%%%%
\section{Conclusion}\label{conclusion}
{This work} presents a stabilized finite element discretization of the full nonlinear steady Darcy-Brinkman-Forchheimer model of the porous medium, {also highlighting the effect of dimensionless constants and the role of the Forchheimer term 
in the steady lid-driven cavity flow}. We  
{employ} a beneficial {Grad-Div stabilization technique} and 
{implement} a Schur complement-based preconditioner for the established linear solver. 
%{The strengths of our algorithm are} 
%the state-of-art preconditioning strategy that significantly reduces the computational cost, 

{The efficient preconditioning strategy results in a few Newton's iteration number even with the full nonlinear model, {yielding the stable and convergent numerical solutions against the spurious oscillations that may occur from the \textit{inf-sup} and discontinuous boundary conditions for the lid-driven cavity flow} with {relatively} high Reynolds number problem {for the porous media}.
  %and 
Moreover, local mesh adaptivity works effectively for the large scale sparse problem via efficient and independent procedure.}}
{With all these approaches, different viscous and inertial resistances are applied, and} we find that the full nonlinear Darcy-Brinkman-Forchheimer model can position itself in between %the 
%different models 
the linearized Brinkman and Darcy-Brinkman/Navier-Stokes
in the turbulent flow regime {due to the role of Forchheimer term}. {When the $Da$ or $Re$ is significantly low, the model can replace the Brinkman, but it also can capture a very minute regime change as the $Da$ gets increased.}

{%needs further in-depth investigations. 
Meantime, the study still needs to be further investigated relating to (but not limited to) more realistic fluid properties and physics such as porosity/permeability variations and Reynolds number based on the experimental reference values (e.g., the reference distance and discharge), and multi-phase flow regime with fluid pressure and capillarity. Moreover, the relation between the dimensionless Forchheimer coefficient and real phenomena or the physical implication of coefficient is  still unknown. In terms of computational performance, the study %has its limitation on the 
still lacks the in-depth and comprehensive knowledge on appropriate selection of the Augmented Lagrangian parameter for its full efficiency.} %\hl{Limitation of the study: pressure, porosity, realistic, cf parameter, computational comparisons}}.

{The proposed solver and algorithm for the Darcy-Brinkman-Forchheimer model, with its robustness in the stability and convergence of 
computation on incompressible turbulent flow, can be expanded and utilized straightforwardly to study other flow states or physics, such as the heat flow in the porous medium.} 
%On the far side,  there are many interesting things to look {further} into, but some need to fit into the scope of the current contribution. 
%{Particularly when computing the Schur complement, globally constant viscosity $\nu$ and the Augmented Lagrangian parameter $\gamma$ {have been chosen in this work.} However,} these parameters could be made variable, and meaningful modified pressure mass matrix elements can be computed via:
%\[
%\left( \mathbf{M}_p(\nu, \; \gamma)\right)_{ij} = \int_{\Omega} \left( \nu(\mathbf{x}) + \gamma(\mathbf{x}) \right)^{-1} \phi_i  \, \phi_j \; d\Omega.
%\]
A more challenging problem is to {expand the model investigated in this study} and develop a {coupling model of the nonlinear flow and mechanics of porous media (e.g., the subsurface flow and geomechanics), i.e.,} the Forchheimer model for the non-Darcy flow and 
nonlinear mechanics such as the strain-limiting theories of {elasticity} \cite{lee2020nonlinear,yoon2021quasi,yoon2021finite,yoon2022preferential,yoon2022finite,yoon2024finite} and stuying flow in porous media containing network of fractures modeled by surface mechanics \cite{ferguson2015numerical}. Devising a numerical method for these multi-physics models would also be a significant challenge. 

\section*{Acknowledgements} {Authors appreciate the feedback received from both reviewers, which ultimately contributed to the enhancement of the paper.  The authors acknowledge the support of the College of Science \& Engineering, Texas A\&M University-Corpus Christi, for this research. {First author,} Hyun C. Yoon, would like to appreciate the support of the Korea Institute of Geoscience and Mineral Resources  {(KIGAM; project number {GP2020-006}/{GP2023-005})}. The authors are also indebted to Dr. D. Palaniappan (Professor of Applied Mathematics, Texas A\&M University-Corpus Christi) for several insightful discussions.  }

\bibliographystyle{unsrt}
\bibliography{fluids}

\end{document}